\documentclass[12pt]{article}

\usepackage{latexsym, amssymb,amsmath,amscd, amsfonts, amsthm, ifpdf, fancybox,float,  geometry, epsfig, graphicx, color, colordvi, setspace,caption,subcaption,indentfirst,soul,bm}
\usepackage[numbers,sort&compress]{natbib}
\usepackage[enableskew,vcentermath]{youngtab}
\usepackage{tikz}
\usetikzlibrary {arrows.meta}
\usetikzlibrary{decorations.pathreplacing}
\usepackage{caption}
\allowdisplaybreaks[4]
\newtheorem{thm}{Theorem}[section]
\newtheorem{pro}[thm]{Proposition}

\newtheorem{rem}[thm]{\it Remark}
\newtheorem{defi}[thm]{Definition}
\newtheorem{lem}[thm]{Lemma}
\newtheorem{cor}[thm]{Corollary}

\numberwithin{equation}{section}
\numberwithin{figure}{section}

\sethlcolor{yellow}
\soulregister\cite7
\soulregister\ref7 

\begin{document}
\begin{center}
	
	{\Large \bf New combinatorial proof of Gaussian polynomial and the monotonicity of the Garvan's $k$-rank}
\end{center}

\begin{center}
	{  Wenxia Qu}$^{1}$  and
	{Wenston J.T. Zang}$^{2}$ \vskip 2mm
	
    	$^{1,2}$ School of Mathematics and Statistics, Northwestern Polytechnical University, Xi'an 710072, P.R. China\\[6pt]
	$^{1,2}$ MOE Key Laboratory for Complexity Science in Aerospace, Northwestern Polytechnical University, Xi'an 710072, P.R. China\\[6pt]
	$^{1,2}$ Xi'an-Budapest Joint Research Center for Combinatorics, Northwestern Polytechnical University, Xi'an 710072, P.R. China\\[6pt]
	\vskip 2mm
	
	$^1$quwenxia0710@mail.nwpu.edu.cn, $^2$zang@nwpu.edu.cn
\end{center}

\vskip 6mm \noindent {\bf Abstract.}  Gaussian polynomial, which is also known as $q$-binomial coefficient, is one of the fundamental concepts in the theory of partitions. Zeilberger provided a  combinatorial proof of Gaussian polynomial, which is called Algorithm Z by Andrews and Bressoud. In this paper, we provide a new bijection on Gaussian polynomial, which leads to a refinement of Algorithm Z. Moreover, using this bijection, we provide an alternative proof of generalized Rogers-Ramanujan identity, which was first proved by Bressoud and Zeilberger. Furthermore, we give a combinatorial proof of the monotonicity property of Garvan's $k$-rank, which is a generalization of Dyson's rank and Andrews-Garvan's crank.



\section{Introduction}
Gaussian polynomials, also known as $q$-binomial coefficients, constitute one of the fundamental objects in the theory of integer partitions. Moreover, Pak, Panova, and Vallejo~\cite{Pak-1,Pak-2,Pak2017Bounds,Pak-3,Pak-4} have established deep connections between Gaussian polynomials and Kronecker coefficients in representation theory. Gaussian polynomials also share a strong relationship with the second-order Reed–Muller code in coding theory~\cite{MacWilliams1977}.

Given $M,N\ge 0$,  let $p_{M,N}(n)$ denote the number of partitions of $n$ with at most $M$ parts, each part not exceeding $N$. Then the generating function for $p_{M,N}(n)$ is given by the Gaussian polynomial ${N+M\brack M}$, expressed as follows:
\begin{equation}\label{equ-gauss-binomial}
       \sum_{n=0}^\infty p_{M,N}(n)q^n={N+M\brack M} = 
           \frac{(q;q)_{M+N}}{(q;q)_{M}(q;q)_{N}}.
\end{equation}
Here we use the standard $q$-series notation
\[
(a;q)_n=\prod_{i=1}^{n}(1-aq^{i-1}),\quad
(a;q)_\infty=\prod_{i=1}^{\infty}(1-aq^{i-1}).
\]
There are several ways to prove \eqref{equ-gauss-binomial}. For example, Andrews \cite{andrews1998theory} established this identity by showing that both sides satisfy the same recurrence relation.  Another approach is to interpret the Gaussian polynomial as enumerating the number of $k$-dimensional subspaces in an $n$-dimensional vector space over $\mathbb{F}_q$, see Stanley \cite{Stanley-Enumerative-book}. Zeilberger \cite{Bressound-1989-Generlaized} also provided a combinatorial proof, which is  outlined below. 

\begin{thm}\label{thm-ZA}(Algorithm Z)
    There is a bijection $\Gamma$ between the set of pairs of partitions $(\alpha,\gamma)$ where $\alpha$ is a partition with each part not exceeding $M+N$, $\gamma$ has $N$ parts with each part not exceeding $M$ and the set of pairs of partitions $(\xi, \delta)$ where $\xi$ has at most $N$ parts , $\delta$ is a partition with each part not exceeding $M$.
\end{thm}

It is clear that Theorem \ref{thm-ZA} gives a combinatorial proof of the following identity:
\[\frac{1}{(q;q)_{M+N}}{M+N\brack N}=\frac{1}{(q;q)_M}\frac{1}{(q;q)_N}.\]

This theorem was originally introduced by Zeilberger \cite{Bressound-1989-Generlaized}, who together with Bressoud \cite{Bressound-1989-Generlaized} used the algorithm to provide a combinatorial proof of the generalized Rogers–Ramanujan identity. Andrews and Bressoud \cite{Andrews Bressoud-1984-identities} referred to this method as Algorithm Z and gave a constructive combinatorial proof of the $q$-analog of the Pfaff–Saalschütz summation with its aid.

There are numerous further applications of Algorithm Z. For instance, Bessenrodt \cite{Bessenrodt1995} employed it to give a bijective proof of a strong refinement of the Alladi–Gordon theorem. Joichi and Stanton \cite{Joichi}  noted that Algorithm Z could be applied to the $q$-binomial theorem. Moreover, Algorithm Z was used by Chen, Chen, Fu, and Zang \cite{ChenChenFuZang2011} to provide a new combinatorial proof of Ramanujan’s $_1 \psi _1$ summation, and by Fu \cite{Fu2007a} to offer a combinatorial interpretation of the Lebesgue identity. Guo and Zeng \cite{GuoZeng2010} also applied Algorithm Z to establish a combinatorial proof of a curious $q$-binomial coefficient identity. By combining novel combinatorial bijections with Algorithm Z, Berndt, Kim, and Yee \cite{Berndt} obtained the first complete combinatorial proofs for a family of identities from Ramanujan’s lost notebook arising from Heine’s transformation and partial theta functions.

Our first main result is to give a refinement of Algorithm Z. For fixed integers $M$, $N$, let $A_{M,N}(n)$ denote the set of partition pairs $(\alpha, \beta)$ such that $\alpha$ is a partition with at most $N$ parts, each part not exceeding $M$, and $\beta$ is a partition with each part lying in $[M+1,M+N]$ satisfying $|\alpha|+|\beta|=n$. Let $B_{N}(n)$ denote the set of partitions $\gamma $ of $n$ with at most $N$ parts. Then we have the following theorem.
	
\begin{thm}\label{thm-suanfa}
   Given positive integer $N$ and for any integer $M\ge 0$, there exists a bijection $\phi_M$ between $A_{M,N}(n)$ and $B_{N}(n)$.
\end{thm}

 Theorem \ref{thm-suanfa} can be viewed as a combinatorial proof of the following identity.
\begin{equation}\label{eq-suanfaeq}
\frac{1}{(q;q)_{N}}=\frac{1}{(q^{M+1};q)_{N}}{M+N\brack N}.
\end{equation}
Note that \eqref{eq-suanfaeq} is a transformation of \eqref{equ-gauss-binomial}. Moreover, Theorem \ref{thm-suanfa} implies the following refinement of Theorem \ref{thm-ZA}.

\begin{cor}
   Preserve the notation of the Theorem \ref{thm-ZA}, there is a bijection between $(\alpha,\gamma)$ and $(\xi,\delta)$. Moreover,  the subpartition of $\alpha$ consisting of all parts not exceeding $M$ coincides with $\delta$. 
\end{cor}

It is worth noting that combinatorial proofs utilizing Algorithm Z may be simplified with the aid of Theorem \ref{thm-suanfa}. For example, we can provide an alternate combinatorial proof of the following lemma using Theorem \ref{thm-suanfa}; this lemma is a key step in proving the generalized Rogers–Ramanujan identity in \cite{Bressound-1989-Generlaized}. To state the lemma, we first introduce two definitions: $R_{k,m}(n)$ and $S_{k,m}(n)$.

Let $R_{k,m}(n)$ denote the set of partition pairs $(\lambda,\delta)$ such that $\lambda$ is an ordinary partition, and $\delta$ is an explicit partition with exactly $|m|$ parts as described below.
    \begin{equation}\label{eq-defofdelta}
    \begin{cases}
        ((2k+1)(m-1)+(k+1), (2k+1)(m-2)+(k+1), \ldots,  k+1)&\text{if }m\ge0;\\
        ((2k+1)(-m-1)+k, (2k+1)(-m-2)+k, \ldots,  k)&\text{if }m<0.
    \end{cases}
\end{equation}
Moreover, $|\lambda|+|\delta|=n$.

Let $S_{k,m}(n)$ denote the set of $4$-tuple partition $\Omega=(s_1,\ldots,s_k;\alpha, \beta, \gamma, \xi)$ satisfies the following restrictions.
\begin{itemize}
    \item[(1)] $s_1\ge s_2\ge\cdots\ge s_k\ge 0$ are nonnegative integers;
    \item[(2)] $\alpha$ is a partition with all parts not less than $s_k$, moreover, the length of the $i$-th Durfee square in $\alpha$ equals $s_i$, where $1\le i\le k$;
    \item[(3)] $\beta$ is a partition with all parts strictly larger than $s_k$ and less than or equal to $2s_k$;
    \item[(4)] $\gamma$ is an explicit partition with the following form:
\begin{equation*}
    \begin{cases}
        (m,m-1, \ldots, 2,1), & \text{if }m\ge 0;\\
        (-m-1, -m-2, \ldots,1), & \text{if }m<0.
    \end{cases}
\end{equation*}
\item[(5)] $\xi$ is a partition with less than or equal to  $s_k+m$ parts and all parts less than or equal to $s_k-m$;
\item[(6)] $|\alpha|+|\beta|+|\gamma|+|\xi|=n$.
\end{itemize}

Then \cite[Lemma 3.4]{Bressound-1989-Generlaized} can be restated as follows:

\begin{thm}{\rm (\cite[Lemma 3.4]{Bressound-1989-Generlaized})}\label{thm-B1989-lem3.4}
Given integral $m$ and positive integral $k$ and $n$, there exists a one-to-one correspondence $\chi$ between the set $R_{k,m}(n)$ and the set $S_{k,m}(n)$. 
\end{thm}

\begin{rem}
   The original proof of \cite[Lemma 3.4]{Bressound-1989-Generlaized} applies Algorithm Z successively 
$2k-1$ times. In Section \ref{sec-thm-B1989-lem3.4},  by employing the bijection $\phi_M$ in Theorem \ref{thm-suanfa}, we achieve the same result in only three steps.
\end{rem}

Our second main result in this paper concerns Garvan's $k$-rank. 
Recall that the rank of an ordinary partition was introduced by Dyson \cite{Dyson1944rank} as the largest part minus the number of parts. The crank of an ordinary partition was defined by Andrews and Garvan \cite{Andrews and Garvan 1988 crank} as the largest part if the partition contains no ones, otherwise as the difference between the number of parts larger than the number of ones and the number of ones. It should be noted that Dyson \cite{Dyson1944} conjectured that rank can provide combinatorial interpretations of the first two Ramanujan congruences, which was confirmed by Atkin and Swinnerton-Dyer in \cite{Atkin1954some}. Andrews and Garvan show that crank can give combinatorial interpretations of all the three Ramanujan congruences in \cite{Andrews and Garvan 1988 crank}. For more details about rank and crank, see \cite{AndrewsGarvan1989Ramanujan,AndrewsOno2005Ramanujan,BringmannDousse2016,BringmannOno2010the,BringmannOno2010dyson,Garvan1990the,Lewis1911on}.

If we let $N(m,n)$ denote the number of partitions of $n$ with rank $m$ and $M(m,n)$ denote the number of partitions of $n$ with crank $m$, the generating functions of $N(m,n)$ and $M(m,n)$ were given by \cite[Eq.(2.13)]{Atkin1954some} and \cite[Eq.(7.20)]{Garvan1988new} as follows.
\begin{equation}
    \sum_{n=0}^{\infty}N(m,n)q^n=\frac{1}{(q;q)_{\infty}}\sum_{n=1}^{\infty}(-1)^{n-1}q^{n(3n-1)/2+mn}(1-q^n)
\end{equation}
and 
\begin{equation}
   \sum_{n=0}^{\infty}M(m,n)q^n=\frac{1}{(q;q)_{\infty}}\sum_{n=1}^{\infty}(-1)^{n-1}q^{n(n-1)/2+mn}(1-q^n).
\end{equation}
The monotonicity of $N(m,n)$ and $M(m,n)$ have also been investigated. In \cite{chan-2014ineq}, Chan and Mao gave the monotonicity of $N(m,n)$ as follows.
\begin{thm}{\rm (\cite[Theorem 4]{chan-2014ineq})}\label{thm-chan-2014ineq}
    For all nonnegative integers $m$ and positive integers $n$, 
    \begin{equation*}
        N(m,n)\ge N(m,n-1)
    \end{equation*}
    except when $(m,n)=(\pm 1, 7),(0,8),(\pm 3, 11)$ and when $n=m+2, m\ge0$.
\end{thm}
In \cite{JiZang-2021-uni}, Ji and Zang investigated the unimodality of $M(m,n)$.
\begin{thm}{\rm (\cite[Theorem 1.6]{JiZang-2021-uni})}\label{thm-JiZang-2021-uni}
    For $n\ge 14$ and $0\le m\le n-2$,
    \begin{equation*}
        M(m,n)\ge M(m,n-1).
    \end{equation*}
\end{thm}

In \cite{Garvan1994generalizations}, Garvan introduced a generalized Dyson's rank. He defined $N_{k}(m,n)$ by
\begin{equation}
    \sum_{n=0}^{\infty}N_{k}(m,n)q^n=\frac{1}{(q;q)_{\infty}}\sum_{n=1}^{\infty}(-1)^{n-1}q^{n((2k-1)n-1)/2+|m|n}(1-q^n),
\end{equation}
for any positive integer $k$, in which the case $k=2$ coincides with Dyson's rank and the case $k=1$ corresponds to the Andrews-Garvan crank. There are several studies on the property of $k$-rank, see \cite{BringmannMahlburg2013,DixitYee2013,Garvan1994generalizations,Mao2014asymptotic,Waldherr2013Asymptotics} for example. 

The second main result of this paper is to give a combinatorial proof of the following inequality on $N_k(m,n)$ with the aid of Theorem \ref{thm-suanfa}.

\begin{thm}\label{thm-yingyong}
For $k\ge 3$, $n\ge k-1$, $m\in \mathbb{Z}$, we have
\[
N_{k}(m,n+1)\ge N_{k}(m,n)
\]
except for the cases $n= |m|+k-1$ and $(k,m,n)= (0,3,8)$.
\end{thm}

Combining Theorem \ref{thm-chan-2014ineq} and Theorem \ref{thm-JiZang-2021-uni}, we obtain the monotonicity of $N_k(m,n)$ as stated below.

\begin{cor}
    For any integer $k, m, n$ satisfying $k\ge 1$, $m\ge 0$ and $n\ge k-1$, we have
    \begin{equation}
       N_{k}(m,n)\ge N_{k}(m,n-1) 
    \end{equation}
    except when either $n=m+k$ or 
    $$(k,m,n)\in \{(1,2,5), (1,3,10), (1,4,9), (1,6,13), (2,1,7), (2,0,8), (2,3,11), (3,0,9)\}.$$
\end{cor}

It should be denoted that the bijection $\phi_M$ in Theorem \ref{thm-suanfa} relates to two new combinatorial structures, namely $C_{M,N}(n)$ and $D_{M,N}(n)$. For  positive integers $M,N$, let $C_{M,N}(n)$ denote the set of partition $\delta$ of $n$, where  $\delta=(\delta_1,\delta_2,\ldots,\delta_{N})$ is a partition with $N$ non-negative parts and $0\le \delta_{i}-\delta_{i+1}\le M+N-i$, here we use the convention that $\delta_{N+1}=0$. We use $D_{M,N}(n)$ to denote the set of partition pairs $(\pi, \mu)$ in which $\pi$ is a partition with all parts lying in $[M+1, M+N-1]$ and the number of appearances of $M+i$ does not exceed $N-i$, and $\mu$ is the partition with at most $N$ non-negative parts and each part less than or equal to $M$. Moreover, $|\pi|+|\mu|=n$. Then we have the following result, which will play a crucial role in the proof of Theorem \ref{thm-suanfa}.

\begin{thm}\label{thm-CandD}
    There exists a bijection $\psi$ between $C_{M,N}(n)$ and $D_{M,N}(n)$
\end{thm}

This paper is organized as follows. In Section \ref{sec-thm-CandD}, we give a proof of Theorem \ref{thm-CandD}. Explicit constructions of $\psi$ and the inverse map $\psi^{-1}$ will be given. Some properties of $\psi$ will also be discussed in this section. In Section \ref{sec-thm-suanfa}, we will prove Theorem \ref{thm-suanfa} with the aid of Theorem \ref{thm-CandD}. Section \ref{sec-thm-B1989-lem3.4} is devoted to providing an alternative combinatorial proof of Theorem \ref{thm-B1989-lem3.4} via the bijection $\phi_M$ introduced in Theorem \ref{thm-suanfa}.  The combinatorial proof of Theorem \ref{thm-yingyong} will be given in Section \ref{sec-yingyongNkmn}. 

\section{Proof of Theorem \ref{thm-CandD}}\label{sec-thm-CandD}
This section aims to prove Theorem \ref{thm-CandD}. We will present our proof in four subsections. In Subsection \ref{subs-1}, we will construct the map $\psi$ from $C_{M,N}(n)$ to $D_{M,N}(n)$ and then introduce its inverse map $\psi^{-1}$. Although the maps $\psi$ and $\psi^{-1}$ are explicitly described, it remains non-trivial to verify that $\psi$ is the desired map in Theorem \ref{thm-CandD} and $\psi^{-1}$ is the inverse map of $\psi$. In Subsection \ref{subs-2} we will first establish several key properties  of $\psi$, and then use these properties to show that the image of $\psi$ lies in $D_{M,N}(n)$. Subsection \ref{subs-3} is devoted to proving that $\psi^{-1}$ is indeed a map from $D_{M,N}(n)$ to $C_{M,N}(n)$. Finally, in Subsection \ref{subsec-2.4}, we will show that $\psi^{-1}$ acts as the inverse of $\psi$ and this completes the proof of Theorem \ref{thm-CandD}.

\subsection{The description of $\psi$ and $\psi^{-1}$}\label{subs-1}
The main purpose of this subsection is to give explicit descriptions of the bijection $\psi$ and its inverse map $\psi^{-1}$. To this end, we first introduce some notations on integer partitions. Then we sketch Algorithm Z, which will be used in the construction of $\psi$. Next we provide the descriptions of the map $\psi$ and its inverse  $\psi^{-1}$. 

Here and throughout this paper, an integer partition of $n$ is a finite sequence of non-increasing positive integers $\lambda=(\lambda_1,\ldots,\lambda_\ell)$ such that $\lambda_1+\cdots+\lambda_{\ell}=n$. We use $\ell(\lambda)$ to denote the  length of $\lambda$, and let $|\lambda|$ denote the total sum of all parts of $\lambda$.
Moreover, for convenience, we also write a partition $\lambda$ of $n$ as $(1^{f_1(\lambda)}, 2^{f_2(\lambda)},\ldots,n^{f_{n}(\lambda)})$, where $f_i(\lambda)$ denotes the number of appearances of $i$ in $\lambda$. We use the convention that we may omit the term $i^0$ when $f_i(\lambda)=0$, and we may write $i^1$ as $i$ for short when $f_i(\lambda)=1$. Given $\lambda,\mu$, let $\lambda\cup \mu$ be the union of $\lambda$ and $\mu$. In other words,  $f_i(\lambda\cup\mu)=f_i(\lambda)+f_i(\mu)$ for any $i$. Furthermore, the conjugation of $\lambda$, which is denoted by $\lambda'$, is another partition defined as follows:
\[\lambda'_j=\#\{i\colon \lambda_i\ge j\},\]
where $1\le j\le \ell(\lambda)$.
For example, let $\lambda=(5,5,3,2,1,1,1)$. Then $\ell(\lambda)=7$, $|\lambda|=18$, and we may write $\lambda$ as $(1^3,2,3,5^2)$. Moreover, $\lambda'=(7,4,3,2,2)$ and $\lambda\cup\lambda'=( 1^3,2^3,3^2,4,5^2,7)$.
Now we sketch Algorithm Z. We only describe the map $\Gamma$, which will be used in the construction of $\phi$. For the inverse map $\Gamma^{-1}$, we refer the reader to \cite{Bressound-1989-Generlaized} for details.
    
\noindent{\textit{Description of Algorithm Z}. }
    Given a partition $\xi=(\xi_1,\ldots, \xi_N)$ and $\delta=(\delta_1,\ldots,\delta_M)$. For any $1\le i\le M$, let $\gamma_i$ ($0\le \gamma_i\le N$) be the unique integer such that $\xi_{N-\gamma_i}\ge \delta_i-\gamma_i\ge \xi_{N-\gamma_i+1}$. Here we use the convention that $\xi_0=+\infty$ and $\xi_{N+1}=0$. Let 
        \begin{equation}\label{eq-az-alpha}
            \alpha_j=\begin{cases}
                \xi_{t},  &\text{if }j=t+\gamma'_{N-t+1} \text{ for some } 1\le t\le N;\\
                \delta_i-\gamma_i,  &\text{if }j=N-\gamma_i+i \text{ for some }1\le i\le M,
            \end{cases}
        \end{equation}
where $1\le j\le M+N$. Then we define $\Gamma(\xi,\delta)=(\alpha,\gamma)$.

\begin{rem}
    It is worth mentioning that the statement in Theorem \ref{thm-ZA} is slightly different from the above algorithm.  To be specific, the partitions $\alpha$, $\gamma$ and $\delta$ mentioned above correspond to the conjugates of these partitions in Theorem \ref{thm-ZA}. The partition $\xi$ above remains the same as in Theorem \ref{thm-ZA}.
\end{rem}

\begin{rem}\label{rem20}
    It should be noted that if for some $t$, $\delta_i-t\le \xi_{N-t}$, then $\gamma_i\le t$. Assume the contrary, if $\gamma_i>t$, then $\delta_i-\gamma_i<\delta_i-t\le \xi_{N-t}\le \xi_{N-\gamma_i+1}$, a contradiction. We emphasize this fact since it will be used in the proof of Lemma \ref{lem-513-1}.
\end{rem}
	
{\noindent\bf Description of the map $\psi$. }
Given $\delta\in C_{M,N}(n)$, define $\delta^0=\delta$ and initialize $\Pi^0=(\Pi^{0,1},\ldots,\Pi^{0,N})$ to be the $N$-tuple of partitions $(\emptyset,\ldots,\emptyset)$. If $\delta^0_1=\delta_1\le M$, set $\psi(\delta)=(\emptyset,\delta)$.  For $i\ge 0$, let $\delta^{i}=(\delta^i_1, \ldots, \delta^i_N)$ and $\Pi^{i}=(\Pi^{i,1}, \ldots, \Pi^{i,N})$. If $\delta^i_1>M$, we apply the following operation to obtain $\delta^{i+1}$ and $\Pi^{i+1}$.  
\begin{itemize}
    \item[\textbf{Step 1:}]   Let $k_i$ be the maximum integer such that $\delta_{k_i}^i\ge M+1$. Since $\delta^i_1>M$,  such $k_i$ exists. Define
    \begin{equation}\label{equ-delta-i}
        \overline{\delta^{i}}:=(\delta_{k_i+1}^i, \delta_{k_i+2}^i, \ldots, \delta_{N}^i)\quad\text{and}\quad \tilde{\delta^i}=( \delta_{1}^i-(M+1), \delta_{2}^i-(M+1), \ldots,\delta_{k_i}^i-(M+1)).
    \end{equation}
    \item[\textbf{Step 2:}]  Perform Algorithm Z on $\overline{\delta^i}$ and $\tilde{\delta^i}$ to obtain a partition $\delta^{i+1}$ with $N$ nonnegative parts and a partition $\gamma^{i+1}$ with $k_i$ nonnegative parts and each part not exceeding $N-k_i$.
\item[\textbf{Step 3:}] Define 
\begin{equation}\label{eq-def-fi}
    f^i(j)=\begin{cases}
        N-k_i-\gamma^{i+1}_j+j,&\text{if }1\le j\le k_i;\\
        j-k_i+\gamma^{i+1'}_{N-j+1},&\text{if }k_i+1\le j\le N.
    \end{cases}
\end{equation}
 Then by \cite[(2.7)]{Dyson1989mappings}, $f^i(j)$ is a bijection from $\{1,2,\ldots,N\}$ to itself. 
\item[\textbf{Step 4:}] Let $\overline{\Pi}^{i,j}=(\overline{\Pi}^{i,j}_1, \ldots, \overline{\Pi}^{i,j}_i, \overline{\Pi}^{i,j}_{i+1})$ be a partition of length $i+1$ defined by
\begin{equation}\label{eq-def-overlinepi}
    \overline{\Pi}^{i,j}_s=\begin{cases}
        \Pi^{i,j}_s+\#\{N-\gamma^{i+1}_j+1\le t\le N\colon \ell(\Pi^{i,t})=i-s+1\},&\text{if }1\le s\le i;\\
        M+1+\#\{N-\gamma^{i+1}_j+1\le t\le N\colon \ell(\Pi^{i,t})=0\},&\text{if }s=i+1.
    \end{cases}
\end{equation}
Now  define $\Pi^{i+1}$ as follows:
\begin{equation}\label{equ-pi-i+1-w}
    \Pi^{i+1,w}=\begin{cases}
        \Pi^{i,t},&\text{if } w=f^i(t) \text{ for } k_i+1\le t\le N;\\
        \overline{\Pi}^{i,t},&\text{if }w=f^i(t) \text{ for } 1\le t\le k_i.
    \end{cases}
\end{equation}
 \end{itemize}

We apply the above operation iteratively until $\delta^{q_\delta}_1\le M$ holds for some $q_\delta\ge 0$ and define $\psi(\delta)=(\bigcup_{t=1}^N \Pi^{q_\delta,t},\delta^{q_\delta})$.

\begin{rem}\label{rem21}
    In Step 4, the fact that $\overline{\Pi}^{i,j}$ is a partition is not immediately obvious; it is an immediate consequence of Lemma \ref{lem-pialpha}.
\end{rem}
 For example, let $M=2$, $N=10$, $\delta^{0}=(28,26,20,12,6,6,5,3,1,1)\in C_{2,10}(108)$ and $\Pi^0=(\emptyset, \emptyset, \emptyset, \emptyset, \emptyset, \emptyset, \emptyset, \emptyset, \emptyset, \emptyset)$. 

 First, we have $k_0=8$, we perform \textbf{Step 1} to get
 \begin{equation}\nonumber
	     \overline{\delta^{0}}=(1,1),\quad \tilde{\delta}^{0}=(25,23,17,9,3,3,2,0).
	 \end{equation}

After \textbf{Step 2, 3, 4}, we can get 
 \begin{equation*}
	     \delta^{1}=(23, 21, 15,7, 1, 1, 1, 1, 1, 0)
	 \end{equation*}
 and
 \begin{equation*}
	     \Pi^1=((5), (5), (5), (5), (5), (5), \emptyset, (4), \emptyset, (3)).
	 \end{equation*}

Now $k_1=4$. Performing \textbf{Step 1, 2, 3, 4} we can get 
\begin{equation*}
	     \delta^{2}=(14,12,6,1,1,1,1,1,1,0)
	 \end{equation*}
 and 
 \begin{equation*}
	     \Pi^2=((9,5), (9,5), (9,5), (5), (5), \emptyset, (7,4), (4), \emptyset, (3)).
\end{equation*}
Now $k_2=3$. By \textbf{Step 1, 2, 3, 4} we get
 \begin{equation*}
	     \delta^{3}=(4,2,1,1,1,1,1,1,1,0)
	 \end{equation*}
 and 
 \begin{equation*}
	     \Pi^3=((10,9,5), (10,9,5), (5), (5), \emptyset, (7,4), (4), (9,6,4), \emptyset, (3)).
	 \end{equation*}
Now $k_3=1$. Performing \textbf{Step 1, 2, 3, 4} again we deduce
 \begin{equation*}
	     \delta^{4}=(2,1,1,1,1,1,1,1,0,0)
	 \end{equation*}
 and 
 \begin{equation*}
	     \Pi^4=((10,9,5), (5), (5), \emptyset, (7,4), (4), (9,6,4), \emptyset, (10,9,6,3), (3)).
	 \end{equation*}
We now have $\delta^4_1=2\le M$, so the iteration terminates. Thus we get $\mu=\delta^4=(2,1,1,1,1,1,1,1,0,0)$ and $\pi=\bigcup_{t=1}^{N}\Pi^{4,t}=(3^2, 4^3, 5^3, 6^2, 7^1, 9^3, 10^2)$.

We begin by analyzing the function $f^i(t)$ and introduce some notations that will be used frequently in this section.
The function $f^i(t)$, defined in \eqref{eq-def-fi} is a bijection from $\{1,2,\ldots,N\}$ to itself. In fact, from \eqref{eq-az-alpha}, we see that if $f^i(t)=j$ then 
\begin{equation}\label{eq-fi-delta}
    \delta^{i+1}_j=\begin{cases}
        \delta^i_t&\text{if }  k_i+1\le t\le N;\\
        \delta^i_t-(M+1)-\gamma^{i+1}_t&\text{if } 1\le t\le k_i.
    \end{cases}
\end{equation} 
In other words, \eqref{equ-pi-i+1-w} and \eqref{eq-fi-delta} imply that $\delta^{i+1}_{j}$ is generated from $\delta^i_t$ and $\Pi^{i+1}_j$ is generated from $\Pi^i_t$.
We also introduce the notation $f_j^i(t)$ to denote the index $r$ such that $\delta^{i+1}_r$ is generated from $\delta^j_t$. That is, 
\[f_j^i(t)=f^i(f^{i-1}(\cdots f^j(t)\cdots)).\]
Here we adopt the convention that $f^{i}_{j}(t)=t$ for all $1\le t\le N$ and $j>i$. Moreover, given $\delta\in C_{M,N}(n)$, let $q_\delta$ denote the total times of iterations to obtain $\psi(\lambda)$. More precisely, let $\psi(\delta)=(\pi,\mu)$, then $\delta^{q_\delta}=\mu$. Furthermore, throughout the remainder of this section, the symbol $k_i$ denotes the maximum integer such that $\delta^i_{k_i} \geq M + 1$, as defined in \textbf{Step 1}, we adopt the convention that $k_{-1}=N$, $k_{q_{\delta}}=0$.

We now describe the construction of the inverse map of $\psi$.

{\noindent\bf Description of the map $\psi^{-1}$. }
Let $(\pi,\mu)\in D_{M,N}(n)$ where 
$$\pi=((M+1)^{f1}, (M+2)^{f_2}, \ldots, (M+N)^{f_{N}})$$
with $0\le f_{i}\le N-i$ and $$M\ge \mu_{1}\ge \mu_{2}\ge \cdots \ge \mu_{N}\ge 0.$$

Consider a table \bm{$T$} with $N+1$ rows and $N$ columns. Label the cell in the $i$-th row and $j$-th column by $(i,j)$, where $0\le i\le N$ and $1\le j\le N$. Define a linear order $\preceq$ on $(i,j)$ as follows:
\[(i,j)\preceq (h,k) \text{ iff either } i<h \text{ or both } i=h \text{ and } j\ge k. \]
First, we put $\mu_i$ at coordinate $(0,i)$ in the $0$th row of the table. Then we start from $(1,N)$ coordinate and put all parts (from the smallest part to the largest part) of $\pi$ into this table under the following operations:
\begin{itemize} 
    \item We fill the coordinate $(1,N),(1,N-1),\ldots,(1,N+1-f_1)$ with  $(M+1)$  and mark the coordinate $(1,N-f_1)$ with \textit{``$F_1$"}, we further delete the entire column under \textit{``$F_1$"};
    \item For $i=2,\ldots, N$, let $(a,b)$ denote the coordinate of the last \textit{``$F_{i-1}$"},  we fill the next $f_i$ undeleted coordinates $(a_1,b_1),\ldots,(a_{f_i},b_{f_i})$ by means of $\preceq$ with $M+i$. In other words, $(a_1,b_1),\ldots,(a_{f_i},b_{f_i})$ are the next minimum $f_i$ coordinates by means $\preceq$ which are undeleted after $(a,b)$. 
    Then we mark the next undeleted coordinate with \textit{``$F_i$"}, we also delete  the entire column under \textit{``$F_i$"}.
\end{itemize}
 Finally, let $\delta_{N+1-i}$ be the sum of the number in the column where \textit{``$F_i$"} is in and define
\begin{equation}
\delta=\psi^{-1}(\mu,\pi)=(\delta_1,\delta_2,\ldots,\delta_N).
\end{equation}

For example, given $M=2$, $N=10$, 
\[
\mu=(2,1,1,1,1,1,1,1,0,0) \quad\text{ and } \pi=(3^2,4^3,5^3,6^2,7,9^3,10^2)
\]
where $f_1=2$, $f_2=3$, $f_3=3$, $f_4=2$, $f_5=1$, $f_6=0$, $f_7=3$, $f_8=2$, $f_9=f_{10}=0$. We consider a table \bm{$T$} with $11$ rows and $10$ columns. First, we add $\mu_{1},\mu_{2}, \ldots, \mu_{10}$ to the $0$th row of the table \bm{$T$}. Then we apply above operations on $\pi$ and fill the coordinates with the parts of $\pi$ as shown in Table \ref{tab-T2}.
\begin{table}[h]
    \centering
    \begin{tabular}{c|cccccccccc}
        &\textbf{C 1} & \textbf{C 2} & \textbf{C 3} & \textbf{C 4} & \textbf{C 5 } & \textbf{C 6} & \textbf{C 7} & \textbf{C 8} & \textbf{C 9} & \textbf{C 10} \\ \hline
        R 0 & 2 & 1 & 1 & 1 & 1 & 1 & 1 & 1 & 0 & 0 \\
        R 1 & 5 & 5 & 5 &\textit{$F_2$} & 4 & 4 & 4 & \textit{$F_1$} &  3 & 3 \\
        R 2 & 9 & \textit{$F_6$} & \textit{$F_5$}& /&7& \textit{$F_4$} & 6 & / & 6 & \textit{$F_3$}\\
        R 3  &10 & / & / & / & \textit{$F_7$} &/ & 9 & / & 9 & / \\
        R 4  &\textit{$F_9$} & / & / & / & / & / &\textit{$F_8$}& / & 10 & /   \\
        R 5  & / & / & / & / & / & / & / & / & \textit{$F_{10}$} & / \\
        R 6  & / & / & / & / & / & / & / & / & / & / \\
        R 7  & / & / & / & / & / & / & / & / & / & / \\
        R 8  & / & / & / & / & / & / & / & / & / & / \\
        R 9  & / & / & / & / & / & / & / & / & / & / \\
        R 10  & / & / & / & / & / & / & / & / & / & / \\
        \hline
        &26&6&6&1&12&5&20&1&28&3
    \end{tabular}
    \caption{An example of table \bm{$T$} in $\psi^{-1}$.}
    \label{tab-T2}
\end{table}

By summing the entries in each column, we can easily deduce that $\delta_{10}=1$, $\delta_{9}=1$, $\delta_{8}=3$, $\delta_{7}=5$, $\delta_{6}=6$, $\delta_{5}=6$, $\delta_{4}=12$, $\delta_{3}=20$, $\delta_{2}=26$, $\delta_{1}=28$.

\subsection{The image set of $\psi$}\label{subs-2}

Although the maps $\psi$ are defined as above, it is not clear how to explicitly describe the set of images of $\psi$. The main purpose of this section is to show that for any $\delta\in C_{M,N}(n)$, we have $\psi(\delta)\in D_{M,N}(n)$. To this end, we need to further analyze the properties on $\psi$. We begin by introducing the following two propositions on $\psi$ which will be frequently used throughout this section.

\begin{pro}\label{prop-1}
        For $0\le i\le q_\delta-1$ and $1\le j\le k_i$, we have $\ell(\Pi^{i+1,f^{i}_0(j)})=i+1$. Moreover,
        \begin{equation}\label{eq-f00j-f01j}
            f_0^0(j)=f_0^1(j)=\cdots=f_0^{i-1}(j)=j.
        \end{equation}
    \end{pro}
    \begin{proof}
    When $i=0$, for $1\le j\le k_0$, we have $\ell(\Pi^{1,f^{0}(j)})=\ell(\overline{\Pi}^{0,j})=1$ from \eqref{eq-def-overlinepi} and \eqref{equ-pi-i+1-w}. 
    
    Suppose that for $0\le i\le t<q_{\delta}-1$ and $1\le j\le k_i$, we have $\ell(\Pi^{i+1,f^{i}_0(j)})=i+1$ and $f^{i-1}_0(j)=j$. Now we consider the case $i=t+1$. By the definition of $k_{t+1}$ and the construction of $\psi$, we know that $\delta^{t+1}_j\ge M+1$ for $1\le j\le k_{t+1}$. Since $1\le j\le k_{t+1}< k_t$, the induction hypothesis implies that $f_0^{t-1}(j)=j$. Moreover, note that 
    \[
    \delta^{t+1}_j\ge M+1>M\ge \delta^t_{k_t+1}=\overline{\delta^t_1}.
    \] Thus, by the construction of Algorithm Z \eqref{eq-az-alpha}, we deduce that the number $\gamma^{t+1}_j$ in Step 2 is equal to $N-k_t$. Therefore, from \eqref{eq-def-fi} we obtain 
    \begin{equation}\label{eq-f0tj}
        f_0^t(j)=f^t(f_0^{t-1}(j))=f^t(j)=j.
    \end{equation}
    By the induction hypothesis and \eqref{eq-f0tj}, we have $\ell(\Pi^{t+1,j})=t+1$. Then by \eqref{eq-def-overlinepi} and \eqref{equ-pi-i+1-w}, we deduce that $\ell(\Pi^{t+2,f^{t+1}_0(j)})=\ell(\overline{\Pi}^{t+1,j})=\ell(\Pi^{t+1,j})+1=t+2$. This completes the proof.
    \end{proof}

\begin{pro}\label{prop-2}
    For $0\le s\le q_{\delta}$, $k_s+1\le t\le k_{s-1}$ and $0\le s\le i\le q_{\delta}$,  we have $\ell(\Pi^{i,f^{i-1}_0(t)})=s$.
\end{pro}
\begin{proof}
Fix $t$ with $ k_s+1\le t\le k_{s-1}$, by Proposition \ref{prop-1} we have  $\ell(\Pi^{s,f^{s-1}_0(t)})=s$. Moreover, if $\delta^{s}_{f^{s-1}_0(t)}\ge M+1$, then by the definition of $k_{s}$, we know $f^{s-1}_0(t)\le k_s$. From \eqref{eq-f00j-f01j} we see that $f^{s-1}_0(f^{s-1}_0(t))=f^{s-1}_0(t)$, which is contradict to $f^{s-1}_0$ is a bijection. Hence $\delta^{s}_{f^{s-1}_0(t)}\le M$, which implies $\delta^{s}_{f^{s-1}_0(t)}=\overline{\delta}^{s}_{f^{s-1}_0(t)-k_s}$.
By the construction of $\psi$, we see that for any $i\ge s$, 
\[\delta^{i}_{f^{i-1}_0(t)}=\delta^{s}_{f^{s-1}_0(t)}\le M\quad\text{and}\quad\Pi^{i,f^{i-1}_0(t)}=\Pi^{s,f^{s-1}_0(t)}.\]
Thus $\ell(\Pi^{i,f^{i-1}_0(t)})=\ell(\Pi^{s,f^{s-1}_0(t)})=s$.
\end{proof}

We then show $|\pi|+|\mu|=|\delta|$  as given below.

\begin{lem}\label{lem-21-weight}
    Given $\delta\in C_{M,N}(n)$ and $\psi(\delta)=(\pi,\mu)$, we have
    \begin{equation}
        |\pi|+|\mu|=|\delta|.
    \end{equation}
\end{lem}

\begin{proof}
It suffices to show that for any $0 \le i\le q_\delta$, 
\begin{equation}\label{eq-delta-deltai}   |\delta|=|\delta^i|+\sum_{j=1}^N|\Pi^{i,j}|. 
\end{equation}
When $i=0$, \eqref{eq-delta-deltai} follows from $\delta^0=\delta$ and $\Pi^0=(\emptyset,\ldots,\emptyset)$. Assume $f^i(t)=j$ for some $1\le t,j\le N$ and $0\le i\le q_\delta-1$. We proceed to show that 
\begin{equation}\label{equ-deli-sum}
    \delta^i_t+|\Pi^{i,t}|=\delta^{i+1}_j+|\Pi^{i+1,j}|.
\end{equation}
Then clearly  \eqref{eq-delta-deltai}  holds.

There are two cases.
\begin{itemize}
    \item[Case 1.] If $t>k_i$. In this case, by \eqref{eq-fi-delta} and \eqref{equ-pi-i+1-w}, 
we find that $\delta^i_t=\delta^{i+1}_j$ and $\Pi^{i,t}=\Pi^{i+1,j}$. This yields \eqref{equ-deli-sum}.
    \item[Case 2.] If $t\le k_i$. In this case, again by \eqref{eq-fi-delta} and \eqref{equ-pi-i+1-w}, we deduce that
    \begin{equation}\label{eq-delta-i+1-j}
        \delta^{i+1}_j=\delta^i_t-(M+1)-\gamma^{i+1}_t,
    \end{equation}
    and
    \begin{equation}\label{eq-pii+1j=overline}
        |\Pi^{i+1,j}|=|\overline{\Pi}^{i,t}|=|\Pi^{i,t}|+M+1+\#\{N-\gamma^{i+1}_t+1\le s\le N\colon \ell(\Pi^{i,s})\le i\}.
    \end{equation}
\end{itemize}  
    
    From the analysis in Proposition \ref{prop-1}, it is easy to see that $\ell(\Pi^{i,s})\le i$ for any $1\le s\le N$. Thus 
    \begin{equation}\label{eq-0912-1}
        \#\{N-\gamma^{i+1}_t+1\le s\le N\colon \ell(\Pi^{i,s})\le i\}=\gamma^{i+1}_t
    \end{equation}
    and \eqref{equ-deli-sum} follows from \eqref{eq-delta-i+1-j}, \eqref{eq-pii+1j=overline} and \eqref{eq-0912-1}.

 Thus in both cases \eqref{equ-deli-sum} holds. This completes the proof.
\end{proof}

We proceed to show that $\pi$ satisfies the restriction in $D_{M,N}(n)$. The key procedure is to prove that $\Pi^{i,j}$ is a distinct partition for all $i,j$ and this is the content of Lemma \ref{lem-Pidistinct}. To this end, we first give a property on $f_0^j(k_i)$ in Lemma \ref{lem-513-1}. Then we use Lemma \ref{lem-513-1} to investigate the length of $\Pi^{i,j}$, which is Corollary \ref{cor-sright}. Next we give a proof of Lemma \ref{lem-Pidistinct} with the aid of Corollary \ref{cor-sright} and Lemma \ref{lem-pialpha}. Finally, Lemma \ref{lem-pimupro} shows that Lemma \ref{lem-Pidistinct} guarantees $\pi$ satisfies the restriction of $D_{M,N}(n)$. Together with Lemma \ref{lem-21-weight}, this confirms that $\psi$ maps $C_{M,N}(n)$ into $D_{M,N}(n)$.
    
\begin{lem} \label{lem-513-1}
For all integers $0\le i\le j\le q_\delta-1$, we have $f^{j}_0(k_i) > f^{j}_{0}(k_{i}+1)$. Consequently $k_i>k_{i+1}$. 
\end{lem}
\begin{proof}
If $j=i$, in this case, we argue by induction on $i$. From Proposition \ref{prop-1}, we see that $f^{i-1}_0(k_i)=k_i$ for all integer $0\le i\le q_{\delta}-1$. Moreover,  by   the induction hypothesis $k_{i-1}>k_i$, we have $f^{i-2}_0(k_i+1)=k_i+1$. Thus from \eqref{eq-def-fi}, we deduce that for $1\le s\le i-1$.
\begin{equation}
    \gamma^{s}_{k_i}=\gamma^s_{k_i+1}=N-k_{s-1}.
\end{equation}
Thus by \eqref{eq-fi-delta}, we have
\begin{equation}\label{ine-deltai-1-delta-i-1}
    \delta^{i-1}_{k_i}-\delta^{i-1}_{k_i+1}=\delta_{k_i}-\delta_{k_i+1}\le M+N-k_i.
\end{equation}
The last inequality follows from the definition of $C_{M,N}(n)$.

On the one hand, since $f^{i-1}(k_i)=k_i$, from \eqref{eq-def-fi} we see that $\gamma^{i}_{k_i}=N-k_{i-1}$. Thus by \eqref{eq-fi-delta},
\begin{equation}
    \delta^i_{k_i}=\delta^{i-1}_{k_i}-(M+1)-(N-k_{i-1}).
\end{equation}

On the other hand, \eqref{eq-fi-delta} and \eqref{eq-def-fi} yields that
\begin{equation}
    \delta^i_{f^{i-1}(k_i+1)}=\delta^{i-1}_{k_i+1}-(M+1)-\gamma_{k_i+1}^i=\delta^{i-1}_{k_i+1}-(M+1)-(N-f^{i-1}(k_i+1)-k_{i-1}+k_i+1).
\end{equation}

Thus, by the definition of $\overline{\delta}^i$ and $\tilde{\delta}^i$, we have
\begin{equation}\label{eq-overlinedelta}
    \tilde{\delta}^i_{k_i}=\delta^i_{k_i}-(M+1)=\delta^{i-1}_{k_i}-2(M+1)-(N-k_{i-1}),
\end{equation}
and
\begin{align}
    \overline{\delta}^i_{f^{i-1}(k_i+1)-k_i}=&\delta^i_{f^{i-1}(k_i+1)}\nonumber\\
    =&\delta^{i-1}_{k_i+1}-(M+1)-(N+1-f^{i-1}(k_i+1))+k_{i-1}-k_i.\label{eq-tildedelta}
\end{align}  

Combining \eqref{ine-deltai-1-delta-i-1}, \eqref{eq-overlinedelta} and \eqref{eq-tildedelta}, we deduce that
\begin{align}
\tilde{\delta}^i_{k_i}-(N-f^{i-1}(k_i+1))=&\delta^{i-1}_{k_i}-2(M+1)-(N-k_{i-1})-N+f^{i-1}(k_i+1)\nonumber\\
\le& \delta^{i-1}_{k_i+1}+M+N-k_i-2M-2-N+k_{i-1}-N+f^{i-1}(k_i+1)\nonumber\\
=&\delta^{i-1}_{k_i+1}-M-2-N+f^{i-1}(k_i+1)+k_{i-1}-k_i\nonumber\\
=& \overline{\delta}^i_{f^{i-1}(k_i+1)-k_i}.
\end{align}
By Remark \ref{rem20}, we see that $\gamma^{i+1}_{k_i}\le N-f^{i-1}(k_i+1)$. Using \eqref{eq-def-fi}, we have
\begin{equation}\label{ine-fi-ki-ge}
    f^i(k_i)\ge  f^{i-1}(k_i+1).
\end{equation}
Moreover, from $\gamma^{i+1}_{k_i}\le N-f^{i-1}(k_i+1)$ we find that $\gamma^{i+1'}_{N-f^{i-1}(k_i+1)+1}\le k_i-1$. Combining \eqref{eq-def-fi}, we deduce that
\begin{equation}\label{ine-fi-fi-1-le}
    f^i(f^{i-1}(k_i+1))\le f^{i-1}(k_i+1)-1.
\end{equation}
From \eqref{ine-fi-ki-ge} and \eqref{ine-fi-fi-1-le}, we deduce that
$f^i(k_i)>f^i(f^{i-1}(k_i+1)),$ which is equivalent to $f_0^i(k_i)>f_0^i(k_i+1)$.

It remains to show that $k_{i+1}<k_i$. In fact, from $f^i(k_i)>f^i(f^{i-1}(k_i+1))$, we see that
\[\delta^{i+1}_{f^i(k_i)}\le \delta^{i+1}_{f^i_{i-1}(k_i+1)}=\delta^{i}_{f^{i-1}(k_i+1)}\le M.\]
Thus there are at most $k_i-1$ elements in $\delta^{i+1}$ is larger than $M$, namely $\delta^{i+1}_{f^i{(1)}}, \ldots, \delta^{i+1}_{f^i{(k_i-1)}}$. This yields $k_{i+1}<k_i$.

We now assume that $j>i$. From the definition of $f^i(t)$ in \eqref{eq-def-fi}, we deduce that $f^r(t_1)>f^r(t_2)>k_{r+1}$ for any $0\le r\le q_{\delta}-1$ and $t_1>t_2> k_r$.
Thus iteratively using this inequality, we have
\begin{equation}
    f^j_{i+1}(t_1)>f^j_{i+1}(t_2)
\end{equation}
for any $t_1>t_2>k_{i+1}$.
Moreover, since $f_0^i(k_i)>f_0^i(k_{i}+1)> k_{i+1}$, we see that
\begin{equation}
    f_0^j(k_i)=f^{j}_{i+1}(f^i_0(k_i))>f^j_{i+1}(f^i_0(k_i+1))=f^j_0(k_i+1).
\end{equation}
This completes the proof.
\end{proof}

\begin{cor}\label{cor-sright}
For any $1\le i\le q_\delta$ and $0\le s\le i-1$. Let $r_{i,s}$ denote the minimum integer such that $\ell(\Pi^{i,r_{i,s}})=s$, then there exists $j_i>r_{i,s}$ such that $\ell(\Pi^{i,j_i})=s+1$. 
\end{cor}

\begin{proof}
     Set
    \begin{equation}
        j_i=f^{i-1}_{s+1}(f^{s}_0(k_s)) \quad \text{and}\quad r_{i,s}=f^{i-1}_0(k_s+1)=f^{i-1}_{s+1}(f^{s}_{0}(k_s+1)).
    \end{equation}
    From Proposition \ref{prop-2}, we know that $\ell(\Pi^{s+1, f^{s}_{0}(k_s)})=s+1$ and \[\ell(\Pi^{s+1,f^{s}_{0}(k_s+1)})=\ell(\Pi^{s,f^{s-1}_0(k_{s}+1)})=s.\]
    Thus for all $i>s$, 
    \begin{equation*}
        \ell(\Pi^{i,r_{i,s}})=\ell(\Pi^{i,f^{i-1}_{s+1}(f^{s}_0(k_s+1))})=s,
    \end{equation*}
    and
    \begin{equation*}
    \ell(\Pi^{i,j_i})=\ell(\Pi^{i,f^{i-1}_{s+1}(f^{s}_0(k_s))})=\ell(\Pi^{s+1,f^{s}_0(k_s)})=s+1.
    \end{equation*}
    Moreover, from Lemma \ref{lem-513-1}, we see that for all $0\le s\le q_{\delta}-1$, $f^{s}_0(k_s)>f^{s}_0(k_s+1)$. Furthermore, since $f^i(t)$ remains the relative position when $t\ge k_i$. We deduce that 
    $$j_i=f_{s+1}^{i-1}(f^{s}_0(k_s))> f_{s+1}^{i-1}(f^{s}_0(k_s+1))=r_{i,s}$$
\end{proof}
		
The following lemma gives a direct description of $\Pi^{i,w}$.

\begin{lem}\label{lem-pialpha}
    For any $1\le i \le q_\delta$ and $1\le w\le N$. Let $\Pi^{i,w}=(\Pi^{i,w}_1,\ldots,\Pi^{i,w}_{\ell(\Pi^{i,w})})$, where $0\le \ell(\Pi^{i,w})\le i$.  For $1\le s\le \ell(\Pi^{i,w})$, we have 
    \begin{equation}\label{eq-624.1}
        \Pi^{i,w}_s=M+1+\#\{1\le j\le N\colon \ell(\Pi^{i,j})\le \ell(\Pi^{i,w})-1-s\}+\#\{w+1\le j\le N\colon \ell(\Pi^{i,j})= \ell(\Pi^{i,w})-s\}.
    \end{equation}   
\end{lem}
\begin{proof}
    We prove this lemma by induction on $i$.  
    
    When $i=1$, if $w=f^{0}(t)$ where $t\ge k_0+1$, we know $\Pi^{1,w}=\Pi^{0,t}=\emptyset$. If $w=f^{0}(t)$ where $1\le t\le k_0$, By \eqref{equ-pi-i+1-w} and \eqref{eq-def-overlinepi},  we know 
    \begin{equation}
        \Pi^{1,w}=\overline{\Pi}^{0,t}=(\overline{\Pi}^{0,t}_1)=(M+1+\#\{N-\gamma^{1}_t+1\le k\le N\colon\ell(\Pi^{i,k})=0\}).
    \end{equation}
    Since $w=N-\gamma^{1}_t$, \eqref{eq-624.1} holds for $i=1$.

    Suppose that \eqref{eq-624.1} holds for $k\le i$, we now  consider $\Pi^{i+1,w}$ with the following two cases.

Case 1. If $w=f^{i}(t)$ where $k_i+1\le t\le N$, then by \eqref{equ-pi-i+1-w}, $\Pi^{i+1,w}=\Pi^{i,t}=(\Pi^{i,t}_1,\ldots,\Pi^{i,t}_{\ell(\Pi^{i,t})})$ with $\ell(\Pi^{i,t})=\ell(\Pi^{i+1,w})\le i$. By hypothesis, 
\begin{equation}\label{eq-Piit}
\Pi^{i,t}_{s}=M+1+\#\{1\le j\le N\colon\ell(\Pi^{i,j})\le \ell(\Pi^{i,t})-1-s\}+\#\{t+1\le j\le N\colon\ell(\Pi^{i,j})=\ell(\Pi^{i,t})-s\}.
\end{equation}
By Proposition \ref{prop-2}, it is clear that for any $q_\delta\ge x>y\ge 0$, we have
\begin{equation}\label{eq-set-1lejN}
   \#\{1\le j\le N\colon\ell(\Pi^{x,j})\le y\}=N-k_y.
\end{equation}
Set $y=\ell(\Pi^{i+1,w})-1-s$ and $x=i,i+1$ respectively, we deduce that
\begin{equation}\label{equ-1lejleNcolon}
    \#\{1\le j\le N\colon\ell(\Pi^{i+1,j})\le \ell(\Pi^{i+1,w})-1-s\}=\#\{1\le j\le N\colon\ell(\Pi^{i,j})\le \ell(\Pi^{i,t})-1-s\}.
\end{equation}
Moreover, by \eqref{eq-def-fi} we know $f^i(t_1)>f^i(t_2)$ for $N\ge t_1>t_2\ge k_i+1$. Therefore
\begin{equation}\label{equ-w+1lejlen}
    \{w+1\le j\le N\colon\ell(\Pi^{i+1,j})=\ell(\Pi^{i+1,w})-s\}=\{f^i(j)\colon t+1\le j\le N,\ell(\Pi^{i,j})=\ell(\Pi^{i,t})-s\}.
\end{equation}
From \eqref{eq-Piit}, \eqref{equ-1lejleNcolon} and \eqref{equ-w+1lejlen}, we deduce that
\begin{align*}
    \Pi^{i+1,w}_s=\Pi^{i,t}_{s}=M+1&+\#\{1\le j\le N\colon\ell(\Pi^{i+1,j})\le \ell(\Pi^{i+1,w})-1-s\}\\
    &+\#\{w+1\le j\le N\colon\ell(\Pi^{i+1,j})=\ell(\Pi^{i+1,w})-s\}.
\end{align*}

Case 2. If $w=f^{i}(t)$ where $1\le t\le k_i$,  by \eqref{equ-pi-i+1-w},
\begin{equation*}
\Pi^{i+1,w}=\overline{\Pi}^{i,t}=(\overline{\Pi}^{i,t}_1,\ldots,\overline{\Pi}^{i,t}_{i+1}).   
\end{equation*}
Using \eqref{eq-def-overlinepi}, we have 
\begin{equation}\label{eq-PIi+1}
    \Pi^{i+1,w}_s=\begin{cases}
        \Pi^{i,t}_s+\#\{N-\gamma^{i+1}_t+1\le j\le N\colon \ell(\Pi^{i,j})=i-s+1\}, &\text{if }1\le s\le i;\\
     M+1+   \#\{N-\gamma^{i+1}_t+1\le j\le N\colon \ell(\Pi^{i,j})=0\}, &\text{if }s=i+1.
    \end{cases}
\end{equation}
We claim that $N\ge j> N-\gamma^{i+1}_t$ if and only if $N\ge f^i(j)>w$ and $j> k_i$. On the one hand, if $j>N-\gamma^{i+1}_t$, then from $\gamma^{i+1}_t\le N-k_i$ we deduce that $j>k_i$. Using \eqref{eq-def-fi}, we have
\begin{equation*}
    f^i(j)=j-k_i+\gamma^{i+1}_{N-j+1}{}'\ge j-k_i+\gamma^{i+1}_{\gamma^{i+1}_t}{}'\ge j-k_i+t>N-\gamma^{i+1}_t-k_i+t=w.
\end{equation*}
On the other hand, assume the contrary, if $j>k_i$ and $j\le N-\gamma^{i+1}_t$, then we have
\[f^i(j)=j-k_i+\gamma^{i+1}_{N-j+1}{}'\le j-k_i+\gamma^{i+1}_{\gamma^{i+1}_t+1}{}'<j-k_i+t\le N-\gamma^{i+1}_t-k_i+t=w,\]
a contradiction. This yields our claim.

From the above claim and the observation that, for $j>k_i$,  $\ell(\Pi^{i,j})=\ell(\Pi^{i+1,f^i(j)})$, it is clear that
\begin{equation}\label{eq-num-N-gamma-i+1}
    \#\{N-\gamma^{i+1}_t+1\le j\le N\colon \ell(\Pi^{i,j})=i-s+1\}=\#\{ w+1\le j\le N\colon \ell(\Pi^{i+1,j})=i-s+1\}
\end{equation}
for any $1\le s\le i+1$. Moreover, note that $t\le k_i$ implies $\ell(\Pi^{i,j})=i$ for any $1\le j\le t$, we have
\begin{equation}\label{eq-set-t+1-le-j}
    \{t+1\le j\le N\colon \ell(\Pi^{i,j})= i-s\}=\{1\le j\le N\colon \ell(\Pi^{i,j})= i-s\}
\end{equation}
holds for all $1\le s\le i$.

Substituting \eqref{eq-num-N-gamma-i+1}, \eqref{eq-set-t+1-le-j} and the induction hypothesis into \eqref{eq-PIi+1}, we deduce
    \begin{align*}
        \Pi^{i+1,w}_s=M+1+&\#\{1\le j\le N\colon \ell(\Pi^{i+1,j})\le \ell(\Pi^{i+1,w})-1-s\}\\
        +&\#\{w+1\le j\le N\colon \ell(\Pi^{i+1,j})= \ell(\Pi^{i+1,w})-s\}.
    \end{align*}
    This completes the proof.
\end{proof}

\begin{lem}\label{lem-Pidistinct}
    For any $q_\delta\ge i\ge 0$, $1\le t\le N$, the partition $\Pi^{i,t}$ is a distinct partition.
\end{lem}
\begin{proof}
    Let $\Pi^{i,t}=(\Pi^{i,t}_1, \ldots, \Pi^{i,t}_{\ell(\Pi^{i,t})})$ where $0\le \ell(\Pi^{i,t})\le i$.  For $1\le s\le \ell(\Pi^{i,t})$. Set 
    \[
    u_{s,1}=\#\{1\le j\le N\colon \ell(\Pi^{i,j})\le \ell(\Pi^{i,t})-1-s\}
    \]
    and 
    \[
    u_{s,2}=\#\{t+1\le j\le N\colon \ell(\Pi^{i,j})= \ell(\Pi^{i,t})-s\}.
    \]
    By Lemma \ref{lem-pialpha},
    \begin{equation*}
        \Pi^{i,t}_s= M+1+u_{s,1}+u_{s,2}
    \end{equation*}
    and    
    \begin{equation*}
        \Pi^{i,t}_{s+1}=M+1+u_{s+1,1}+u_{s+1,2}.
    \end{equation*}
    Clearly, 
    \begin{align*}
        &\{1\le j\le N\colon \ell(\Pi^{i,j})\le \ell(\Pi^{i,t})-1-(s+1)\}\cup \{t+1\le j\le N\colon \ell(\Pi^{i,j})= \ell(\Pi^{i,t})-(s+1)\} \\
        &\subseteq \{1\le j\le N\colon \ell(\Pi^{i,j})\le \ell(\Pi^{i,t})-1-s\},
    \end{align*}
    which means $u_{s+1,1}+u_{s+1,2}\le u_{s,1}$. We consider the following two cases.
    
    \textbf{Case 1:} If $u_{s+1,1}+u_{s+1,2} < u_{s,1}$, then $\Pi^{i,t}_{s+1} < \Pi^{i,t}_s$. 
    
    \textbf{Case 2:} If $u_{s+1,1}+u_{s+1,2} = u_{s,1}$, then by definition, $\ell(\Pi^{i,j})\ne \ell(\Pi^{i,t})-s-1$ for all $1\le j\le t$. 
    By Proposition \ref{prop-2} there exists a smallest integer $j> t$ such that $\ell(\Pi^{i,j})=\ell(\Pi^{i,t})-s-1$. By Corollary \ref{cor-sright}, there exists an integer $k> j$ such that $\ell(\Pi^{i,k})=\ell(\Pi^{i,t})-s$.
    Hence, $u_{s,2}>0$, and the corollary follows.
\end{proof}

We are now in a position to show that for any $\delta\in C_{M,N}(n)$, the image $\psi(\delta)$ lies in $D_{M,N}(n)$; that is,  $\psi(\delta)=(\pi,\mu)\in D_{M,N}(n)$.
    
\begin{lem}\label{lem-pimupro}
    For any $\delta\in C_{M,N}(n)$, we have $\psi(\delta)=(\pi,\mu)\in D_{M,N}(n)$; that is, $\pi$ and $\mu$ satisfy the following two restrictions.
    \begin{itemize}
        \item[(1)] $\pi$ is a partition with all parts lying in $[M+1, M+N-1]$ and the number of occurrences of $M+i$ does not exceed $N-i$;
        \item[(2)] $\mu$ is the partition with at most $N$ parts and each part is less than or equal to $M$.
    \end{itemize}	
     \end{lem}
\begin{proof}
    We first show that for any $1\le t\le N$,
         \begin{equation}\label{ine-pi1}
             \Pi^{q_\delta,f^{q_\delta-1}_0(t)}_1\le M+N-t.
         \end{equation}
Let $w_t$ denote $f^{q_\delta-1}_0(t)$, and assume $k_s+1\le t\le k_{s-1}$, with the conventions $k_{-1}=N$ and $k_{q_\delta}=0$. By Proposition \ref{prop-2}, 
\begin{equation}
    \ell(\Pi^{q_\delta,w_t})=s.
\end{equation}
Moreover, from Lemma \ref{lem-pialpha} we have
\begin{equation}\label{equ-Pi-q-delta}
    \Pi^{q_\delta,w_t}_1=M+1+\#\{1\le j\le N\colon \ell(\Pi^{q_\delta,j})\le s-2\}+\#\{w_t+1\le j\le N\colon \ell(\Pi^{i,j})= s-1\}.
\end{equation}
By Proposition \ref{prop-2}, we see that  \[
\{1\le j\le N\colon \ell(\Pi^{q_\delta,j})\le s-2\}=\{f_0^{q_\delta-1}(j)\colon k_{s-2}+1\le j\le N\}.
\]
Therefore
\begin{equation}\label{eq-1lejleN}
    \#\{1\le j\le N\colon \ell(\Pi^{q_\delta,j})\le s-2\}=N-k_{s-2}.
\end{equation}
Next, we compute $\#\{w_t+1\le j\le N\colon \ell(\Pi^{i,j})= s-1\}$. On the one hand, if $t<k_{s-1}$, then $t+1\le k_{s-1}$. By Proposition \ref{prop-2}, we know $\ell(\Pi^{q_\delta,f^{q_\delta-1}_0(t+1)})\ge s$. Then
\[\{w_t+1\le j\le N\colon \ell(\Pi^{i,j})= s-1\}\subseteq \{f_0^{q_\delta-1}(j)\colon t+2\le j\le k_{s-2}\}.\]
On the other hand, if $t=k_{s-1}$, then $t+1=k_{s-1}+1$. From Lemma \ref{lem-513-1}, we find $w_t=f_0^{q_\delta-1}(t)>f_0^{q_\delta-1}(t+1)$. Thus we also deduce
\[\{w_t+1\le j\le N\colon \ell(\Pi^{i,j})= s-1\}\subseteq \{f_0^{q_\delta-1}(j)\colon t+2\le j\le k_{s-2}\}.\]

In both cases, we therefore have 
\begin{equation}\label{ine-wt+1}
    \#\{w_t+1\le j\le N\colon \ell(\Pi^{i,j})= s-1\}\le k_{s-2}-t-1.
\end{equation}
Substituting \eqref{eq-1lejleN} and \eqref{ine-wt+1} into \eqref{equ-Pi-q-delta} yeilds \eqref{ine-pi1}. Thus, every part of $\pi$ lies in $[M+1, M+N-1]$ unless $\pi=\emptyset$. Together with Lemma \ref{lem-Pidistinct} and \eqref{ine-pi1}, we see that for every $1\le i\le N-1$, $M+i$ can appear only among
\[
\Pi^{q_\delta,f^{q_{\delta-1}}_{0}(1)},\Pi^{q_\delta,f^{q_{\delta-1}}_{0}(2)},\ldots,\Pi^{q_\delta,f^{q_{\delta-1}}_{0}(N-i)},
\]
and at most once.  Therefore, the number of appearances of $M+i$ in $\pi$ does not exceed $N-i$. The restriction on $\mu$ is immediate from the construction of $\psi$, and the proof is complete.
\end{proof}

\subsection{The image set of $\psi^{-1}$}\label{subs-3}
In this subsection, we show that for any $(\pi,\mu)\in D_{M,N}(n)$, $\delta=\psi^{-1}(\pi,\mu)\in C_{M,N}(n)$ in the following lemma.
\begin{lem}\label{lem-psi-1inC}
Given $(\pi,\mu)\in D_{M,N}(n)$, let $\delta=(\delta_1,\ldots,\delta_N)=\psi^{-1}(\pi,\mu)$, then we have $0\le\delta_i-\delta_{i+1} \le M+N-i$ for $1\le i \le N$, with the convention that $\delta_{N+1}=0$. In other words, $\delta\in C_{M,N}(n)$.
\end{lem}
\begin{proof}
    First we observed that $\textit{``$F_{1}$"}$, $\textit{``$F_{2}$"}$, $\ldots$, $\textit{``$F_{i}$"}$ must appear in the first $i$th rows of the table \bm{$T$}. This follows that the number of appearances of $M+i$ in $\pi$ does not exceed $N-i$. 
    
    Next, we verify that $0 \le \delta_i-\delta_{i+1} \le M+N-i$ for $1\le i \le N$. To this end, we consider the following two cases.
    
    Case 1. If $\textit{``$F_{i}$"}$ and $\textit{``$F_{i+1}$"}$ are in the same row, we may assume that $\textit{``$F_{i}$"}$ is at the position of coordinate $(k+1,t)$ and $\textit{``$F_{i+1}$"}$ is at the position of coordinate $(k+1, s)$, where  $0\le k\le N-1$ and $1\le s<t\le N$. Let $a_j(1\le j\le k)$ denote the number with the coordinates $(j,t)$ and $b_j(1\le j\le k)$ denote the number with the coordinates $(j,s)$.
    Then we have 
    \[
    \delta_{N+1-i}=\sum_{j=1}^k a_j+\mu_t, \qquad \delta_{N-i}=\sum_{j=1}^k b_j+\mu_s.
    \]
    Moreover, from the construction of $\psi^{-1}$, we see that 
    \begin{equation}\label{equ-M+1-a-1}
        M+1\le a_1\le b_1\le a_2\le b_2\le \cdots \le a_k\le b_k \le M+i.
    \end{equation}
     On the one hand, by \eqref{equ-M+1-a-1} and noticing $\mu_s\ge \mu_t$, we have
     \[
     \delta_{N-i}-\delta_{N+1-i}=\sum_{j=1}^k(b_j-a_j)+\mu_s-\mu_t\ge 0.
     \]
    On the other hand, we have 
    \begin{equation}\label{eq-delta}
    \delta_{N-i}-\delta_{N+1-i}= (\mu_s-\mu_t)+\sum_{j=1}^{k-1}(b_j-a_{j+1})+b_k-a_1.
    \end{equation}
    By \eqref{equ-M+1-a-1}, we find that $\sum_{j=1}^{k-1}(b_j-a_{j+1})\le 0$ and $b_k-a_1\le i-1$. Moreover, from $0\le \mu_t\le\mu_s\le \mu_1\le M$ we deduce that $\mu_s-\mu_t\le M$. Hence \eqref{eq-delta} yields $\delta_{N-i}-\delta_{N+1-i}\le M+i-1$.

    Case 2. If $\textit{``$F_{i}$"}$ and $\textit{``$F_{i+1}$"}$ are not in the same row, we assume that $\textit{``$F_{i}$"}$ is at the position of coordinate $(k+1,t)$ and $\textit{``$F_{i+1}$"}$ is at the position of coordinate $(r+1,s)$. We first show that $r=k+1$ and $s>t$.

    From the construction of $\psi^{-1}$, we know the number between $\textit{``$F_{i}$"}$ and $\textit{``$F_{i+1}$"}$ must be $M+i+1$. After marking \textit{``$F_{1}$"}, \textit{``$F_{2}$"},$\ldots$, \textit{``$F_{i}$"}, the number of columns we can put $M+i+1$ at equals $N-i$. Moreover, by the definition of $\pi$, the number of appearances of $M+i+1$ does not exceed $N-i-1$. Thus we can place $M+i+1$ after $\textit{``$F_{i}$"}$ at most $N-i-1$ times. Together with $r>k$, we deduce that $\textit{``$F_{i+1}$"}$ must lie in the next row of $\textit{``$F_{i}$"}$ and also the right column of $\textit{``$F_{i}$"}$. Thus we have $s>t$ and $r=k+1$.

    Let $c_j(1\le j\le k)$ denote the number which lies at the coordinate $(j,t)$ and $d_j(1\le j\le k+1)$ denote the number which lies at the coordinate $(j,s)$. Then 
    \[
\delta_{N-i}=\sum_{j=1}^{k+1}d_j+\mu_s, \qquad \delta_{N+1-i}=\sum_{j=1}^{k}c_j+\mu_t.
\]
Similar with \eqref{equ-M+1-a-1}, we now have
    \begin{equation}\label{eq-M+1leak+1leM+i}
    M+1\le d_1\le c_1\le d_2\le c_2\le \cdots \le c_k\le d_{k+1}\le M+i.
    \end{equation}
    On the one hand, by \eqref{eq-M+1leak+1leM+i},
    \[
        \delta_{N-i} -\delta_{N+1-i}=\sum_{j=1}^{k}(d_j-c_j)+(\mu_s-\mu_t)+d_{k+1} \le M+i.
    \]
    On the other hand, by \eqref{eq-M+1leak+1leM+i},
    \[
    \delta_{N-i}-\delta_{N+1-i} =\mu_s-\mu_t+d_1+ \sum_{j=1}^{k}(d_{j+1}-c_j)\ge d_1-M\ge 1\ge 0.
    \]
    Let $(N-i)\rightarrow i$ in the above two cases, we deduce that $0\le \delta_i-\delta_{i+1} \le M+N-i$ for $1\le i \le N$. This completes the proof.
\end{proof}

\subsection{Proof of Theorem \ref{thm-CandD}}\label{subsec-2.4}
In this subsection, we conclude the proof of Theorem \ref{thm-CandD} by showing that $\psi^{-1}$ is indeed the inverse map of $\psi$. To this end, it suffices to verify that for any $\delta\in C_{M,N}(n)$, 
\begin{equation}\label{equ-psi-1-psi}
    \psi^{-1}(\psi(\delta))=\delta,
\end{equation}
and for any $(\pi,\mu)\in D_{M,N}(n)$, 
\begin{equation}\label{equ-psi-pais-1}
    \psi(\psi^{-1}(\pi,\mu))=(\pi,\mu).
\end{equation}

The key procedure to prove \eqref{equ-psi-1-psi} is the following lemma.

\begin{lem}\label{lem-psipsi-1}
    Given $0\le j\le N-2$, if  $f_0^{q_\delta-1}(N-j)<f_0^{q_\delta-1}(N-j+1)$, then    
    \[
    M+j+1\in \Pi^{q_\delta,k}\iff f_0^{q_\delta-1}(N-j)<k<f_0^{q_\delta-1}(N-j+1). 
    \]
    Moreover, if  $f_0^{q_\delta-1}(N-j)>f_0^{q_\delta-1}(N-j+1)$, then 
    \[
    M+j+1\in \Pi^{q_\delta,k}\iff k>f_0^{q_\delta-1}(N-j) \text{ or }k<f_0^{q_\delta-1}(N-j+1).
    \]
    Here we adopt the convention that $f_0^{q_\delta-1}(N+1)=+\infty$.  
\end{lem}
\begin{proof}
Let $r=\ell(\Pi^{q_\delta,k})$, from Lemma \ref{lem-pialpha}, we know
\begin{align}
    \Pi^{q_\delta,k}_{s}=M+1&+\#\{1\le \omega \le N : \ell(\Pi^{q_\delta,\omega})\le r-s-1\}\nonumber\\
    &+\#\{k+1\le \omega \le N : \ell(\Pi^{q_\delta,\omega})= r-s\}\label{eq-alphaM+1+jr-s-1r-s}.
\end{align}
Thus $M+j+1\in \Pi^{q_\delta,k}$ if and only if there exists $s$ such that
\begin{equation}\label{eq-j}
    j=\#\{1\le \omega \le N : \ell(\Pi^{q_\delta,\omega})\le r-s-1\}+\#\{k+1\le \omega \le N : \ell(\Pi^{q_\delta,\omega})= r-s\}.
\end{equation}
Using Proposition \ref{prop-2}, we find that 
\begin{equation}\label{eq-set-r-s-1}
    \{1\le \omega \le N : \ell(\Pi^{q_\delta,\omega})\le r-s-1\}=\{f_0^{q_\delta-1}(\omega)\colon k_{r-s-1}+1\le \omega\le N\}.
\end{equation}
Thus
\[
\#\{1\le \omega \le N : \ell(\Pi^{q_\delta,\omega})\le r-s-1\}=N-k_{r-s-1}.
\]
This implies \eqref{eq-j} is equivalent to the following identity
\begin{equation}\label{eq-num-k+1}
    \#\{k+1\le \omega \le N : \ell(\Pi^{q_\delta,\omega})= r-s\}=j-N+k_{r-s-1}.
\end{equation}
From the construction of $\psi$, it can be checked that
\[f_0^{q_\delta-1}(k_{r-s}+1)<f_0^{q_\delta-1}(k_{r-s}+2)<\cdots<f_0^{q_\delta-1}(k_{r-s-1}).\]
Let $c$ be the minimum positive integer such that $k<f_0^{q_\delta-1}(k_{r-s}+c)$, then we see that
\begin{equation}\label{eq-k+1-le-omega-1}
    \{k+1\le \omega \le N : \ell(\Pi^{q_\delta,\omega})= r-s\}
    =\{f_0^{q_\delta-1}(\omega)\colon k_{r-s}+c\le \omega\le k_{r-s-1}\}.
\end{equation}
Therefore
\eqref{eq-num-k+1} implies
\begin{equation}\label{eq-k+1-le-omega}
    \{k+1\le \omega \le N : \ell(\Pi^{q_\delta,\omega})= r-s\}
    =\{f_0^{q_\delta-1}(v)\colon N-j+1\le v\le k_{r-s-1}\}.
\end{equation}
We claim that \eqref{eq-k+1-le-omega} holds if and only if one of the following three cases must hold:
\begin{itemize}
    \item[Case 1.] $\ell(\Pi^{q_\delta,f_0^{q_\delta-1}(N-j+1)})= r-s-1$, $\ell(\Pi^{q_\delta,f_0^{q_\delta-1}(N-j)})= r-s$ and $f_0^{q_\delta-1}(N-j)< k$;
    \item[Case 2.] $\ell(\Pi^{q_\delta,f_0^{q_\delta-1}(N-j+1)})= r-s=\ell(\Pi^{q_\delta,f_0^{q_\delta-1}(N-j)})$ and $f_0^{q_\delta-1}(N-j)< k<f_0^{q_\delta-1}(N-j+1)$;
    \item[Case 3.] $\ell(\Pi^{q_\delta,f_0^{q_\delta-1}(N-j+1)})= r-s$, $\ell(\Pi^{q_\delta,f_0^{q_\delta-1}(N-j)})=r-s+1$ and $k<f_0^{q_\delta-1}(N-j+1)$.
\end{itemize}

On the one hand, from  \eqref{eq-k+1-le-omega}, we see that if $\ell(\Pi^{q_\delta,f_0^{q_\delta-1}(N-j+1)})= r-s-1$, then both sides of \eqref{eq-k+1-le-omega} are empty. Thus $N+j+1>k_{r-s-1}$. Together with \eqref{eq-num-k+1} we deduce that $N-j=k_{r-s-1}$. Using the fact $f_0^{q_\delta-1}(N-j)\not\in\{k+1\le \omega \le N : \ell(\Pi^{q_\delta,\omega})= r-s\}$, we find that $f_0^{q_\delta-1}(N-j)\le k$.  Moreover, by $\ell(\Pi^{q_\delta,k})=r$ and $\ell(\Pi^{q_\delta,f_0^{q_\delta-1}(N-j)})=r-s$, we deduce that $f_0^{q_\delta-1}(N-j)\ne k$. This yields Case 1. 

Otherwise we have $\ell(\Pi^{q_\delta,f_0^{q_\delta-1}(N-j+1)})= r-s$ and $k<f_0^{q_\delta-1}(N-j+1)$. Since $f_0^{q_\delta-1}(N-j)\not\in\{f_0^{q_\delta-1}(\omega)\colon N-j+1\le \omega\le k_{r-s-1}\}$, we see that either $\ell(\Pi^{q_\delta,f_0^{q_\delta-1}(N-j)})= r-s$ and $f_0^{q_\delta-1}(N-j)\le k$ or $\ell(\Pi^{q_\delta,f_0^{q_\delta-1}(N-j)})=r-s+1$. Using the same argument as in Case 1, we find that $f_0^{q_\delta-1}(N-j)\ne k$. This yields \eqref{eq-k+1-le-omega} implies that one of the above three cases holds.

On the other hand, it is trivial to check that if one of Case 1, Case 2 or Case 3 holds, then \eqref{eq-k+1-le-omega} is valid. This completes the proof of our claim.

Now we use the claim to prove Lemma \ref{lem-psipsi-1}. Clearly, Case 2 directly implies the first assertion of Lemma \ref{lem-psipsi-1} directly. 
Moreover, in Case 1,  from Proposition \ref{prop-2}  we know $N-j=k_{r-s-1}$. Thus, by Lemma \ref{lem-513-1}, we have $f^{q_{\delta}-1}_{0}(N-j+1)<f^{q_{\delta}-1}_{0}(N-j)< k$. In Case 3, using the same argument as in Case 1, we have $N-j=k_{r-s}$ and $k<f^{q_{\delta}-1}_{0}(N-j+1)<f^{q_{\delta}-1}_{0}(N-j)$. This establishes the second assertion of Lemma \ref{lem-psipsi-1}.
\end{proof}

To establish \eqref{equ-psi-pais-1}, we begin by introducing some notations. Given $(\pi,\mu)\in D_{M,N}(n)$, let $c(F_i)$ denote the column number of $``F_i"$, $r(F_i)$ denote the row number of $``F_i"$, where $``F_1",``F_2",\ldots,``F_N"$ are the entries in the construction of $\psi^{-1}(\pi,\mu)$.  Moreover, we let $o(i,j)$ denote the integer at coordinate $(i,j)$. (If the coordinate is $``F"$ or deleted, we adopt the convention $o(i,j)=0$.) Assume $r(F_N)=s$, and let $N=r_0>r_1>\cdots>r_{s}=0$ be the positive integers defined as follows:
\begin{equation}\label{eq-def-ri}
\{t\colon r(F_{N-t+1})=i+1\}=\{t\colon r_{i+1}+1\le t\le r_i\},
\end{equation}
 where $0\le i\le s-1$.  From the construction of $\psi^{-1}$, it is clear that for $1\le i\le s$,
 \[
 1\le c(F_{N-r_i})<c(F_{N-r_i-1})<\cdots<c(F_{N-r_{i-1}+1})\le N.
 \]
The above inequality enables us to define $\theta_j^i$ as follows.

\begin{defi}
    For any $1\le i\le s$ and $1\le j\le r_{i}$, we define $\theta^i_j$ as follows:
 \begin{equation}\label{eq-def-theta}
     \theta_j^i=\begin{cases}
         0,&\text{if }c(F_{N-j+1})\in [1,c(F_{N-r_i})];\\
         t,&\text{if }c(F_{N-j+1})\in [c(F_{N-r_i-t+1})+1,c(F_{N-r_i-t})] \\
         &\text{for some }1\le t\le r_{i-1}-r_i-1;\\
         r_{i-1}-r_i,& \text{if }c(F_{N-j+1})\in [c(F_{N-r_{i-1}+1})+1,N].
     \end{cases}
 \end{equation}
\end{defi}
From the construction of $\psi^{-1}$, it is easy to see that for every $1\le t<r(F_{N-j+1})$ we have
\begin{equation}\label{eq-o-theta}
    o(t,c(F_{N-j+1}))=N-r_t-\theta^t_j+M+1.
\end{equation}
We obtain the following result.

\begin{lem}\label{lem-215}
   Given $(\pi,\mu)\in D_{M,N}(n)$, let $\delta=\psi^{-1}(\pi,\mu)$ and apply $\psi$ on $\delta$, given $1\le i\le s$ and $r_{i+1}+1\le j\le r_i$, we have the following results.
   \begin{itemize}
       \item[(1)] We have
       \begin{align}
           r_{i}&=k_{i-1};\label{eq-riki-1}\\
           \delta^i_{f_0^{i-1}(j)}&=\mu_{c(F_{N+1-j})};\label{eq-f0-i-1}\\
           f_0^{i-1}(j)&=f^{i-1}(j)=j+\sum_{t=1}^i \theta^{t}_j.\label{eq-f0i-1-j}
       \end{align}
\item[(2)] For all   $1\le t\le i$,
\begin{equation}\label{eq-Pi-i-f0-i-1}
    \Pi^{i,f_0^{i-1}(j)}_t=o(i+1-t,c(F_{N-j+1}))=M+1+N-r_{i+1-t}-\theta^{i+1-t}_j.
\end{equation}
\item[(3)] Let $r_{i+1}+1\le j_1,j_2\le N$,
\begin{equation}\label{ine-order}
\big(c(F_{N+1-j_1})-c(F_{N+1-j_2})\big)\big(f_0^{i-1}(j_1)-f_0^{i-1}(j_2)\big)>0.
\end{equation}
   \end{itemize}
\end{lem}
\begin{proof}
We prove this lemma by induction on $i$. When $i=0$, we see that \[
\delta_N\le\delta_{N-1}\le\cdots\le \delta_{r_1+1}=\mu_{c(F_{N-r_1})}\le M.
\]
Together with $\delta_{r_1}\ge M+1$, we deduce that $k_0$ coincides with $r_1$ when we perform $\psi$ on the $\delta$ obtained from $\psi^{-1}$, which implies that \eqref{eq-riki-1} holds when $i=0$. From the construction of $\psi^{-1}$, for all $r_1+1\le j\le N$, we have
\[
\delta_j=\delta^0_j=\mu_{c(F_{N-j+1})}.
\]
Thus \eqref{eq-f0-i-1} holds when $i=0$. Moreover, \eqref{eq-f0i-1-j} holds since $\theta^{t}_{j}=0$ and \eqref{eq-Pi-i-f0-i-1} also holds since there is no $t$ with $1\le t\le 0$. Furthermore, from the construction of $\psi^{-1}$, we see that
\[
N\ge c(F_1)>c(F_2)>\cdots>c(F_{N-r_1})\ge 1,
\]
which implies \eqref{ine-order} holds when $i=0$.

Assuming that  \eqref{eq-riki-1}, \eqref{eq-f0-i-1}, \eqref{eq-f0i-1-j}, \eqref{eq-Pi-i-f0-i-1} and \eqref{ine-order} all hold for the case $i-1$, we first verify that the first assertion of Lemma \ref{lem-215} holds for the case $i$. For each $r_{i+1}+1\le j\le r_i$, from the construction of $\psi^{-1}$, we see that  
\begin{equation}\label{equ-delta-j}
    \delta_j=\mu_{c(F_{N-j+1})}+\sum_{t=1}^{i}o(t,c(F_{N-j+1})).
\end{equation}
Using \eqref{eq-def-theta} and \eqref{eq-o-theta}, we deduce that
\begin{equation}
   M+ \sum_{t=1}^{i}(N-r_{t}+M+1) \ge \delta_j\ge \sum_{t=1}^{i}(N-r_{t-1}+M+1).
\end{equation}
From the induction hypothesis, we see that $r_{t}=k_{t-1}$ for $1\le t\le i-1$. By the construction of $\psi$, it is clear that for any $1\le t\le i-1$, we have
\begin{equation}\label{eq-ft0-j=j}
    f^0(j)=\cdots=f^{i-2}(j)=j,\quad \delta^t_{f^{t-1}_0(j)}=\delta_j-\sum_{s=1}^t(M+1+N-r_{s}),
\end{equation}
and
\begin{equation}\label{eq-Pi-t-j}
    \Pi^{t,j}=(M+1+N-r_{t},M+1+N-r_{t-1},\ldots,M+1+N-r_{1}).
\end{equation}
Combining \eqref{eq-o-theta}, \eqref{equ-delta-j} and \eqref{eq-ft0-j=j}, we deduce that 
\begin{align}\label{eq-exp-delta-i-1}
    \delta^{i-1}_j=&\delta_j-\sum_{t=1}^{i-1}(M+1+N-r_t)\nonumber\\
    =&\mu_{c(F_{N-j+1})}+\sum_{t=1}^{i}o(t,c(F_{N-j+1}))-\sum_{t=1}^{i-1}(M+1+N-r_{t})\nonumber\\
    =&\mu_{c(F_{N-j+1})}+o(1,c(F_{N-j+1}))+\sum_{t=2}^{i}(o(t,c(F_{N-j+1}))-(M+1+N-r_{t-1}))\nonumber\\
    =&\mu_{c(F_{N-j+1})}+o(1,c(F_{N-j+1}))+\sum_{t=2}^{i}(r_{t-1}-r_t-\theta_j^t).
\end{align}
From \eqref{eq-exp-delta-i-1}, it is clear that $\delta^{i-1}_j\ge M+1$ for $r_{i+1}+1\le j\le r_i$. Moreover, from the induction hypothesis \eqref{eq-f0-i-1}, we see that 
\[
\delta^{i-1}_{r_i+1}=\delta^{i-1}_{f^{i-2}_0(p)}=\mu_{c(F_{N-p+1})}\le M
\]
for some $p\ge r_i+1$. Thus by the definition of $k_{i-1}$, we deduce that $k_{i-1}=r_i$.

Next, we verify that \eqref{eq-f0-i-1} and \eqref{eq-f0i-1-j} hold. First, we have
\begin{equation}
    \overline{\delta^{i-1}}=(\delta^{i-1}_{r_i+1},\ldots,\delta^{i-1}_N)
\end{equation}
and for all $1\le j\le r_i$,
\begin{align}
    \tilde{\delta}^{i-1}_j&=\delta^{i-1}_j-(M+1)\nonumber\\
    &=\mu_{c(F_{N-j+1})}+(M+1+N-r_1-\theta^1_j)+\sum_{t=2}^{i}(r_{t-1}-r_t-\theta_j^t)-(M+1)\nonumber\\
    &=\mu_{c(F_{N-j+1})}+\sum_{t=1}^{i}(r_{t-1}-r_t-\theta_j^t)\nonumber\\
    &=\mu_{c(F_{N-j+1})}+N-r_i-\sum_{t=1}^{i}\theta_j^t.
\end{align}
Thus
\begin{equation}\label{eq-tilde-delta-i-1-j-left}
    \tilde{\delta}^{i-1}_j-\left(N-r_i-\sum_{t=1}^{i}\theta_j^t\right)=\mu_{c(F_{N-j+1})}.
\end{equation}
From the definition of $\theta_j^t$ \eqref{eq-def-theta}, we see that for each $1\le t\le i$,
\begin{equation}
  \mu_{c(F_{N-r_t-\theta^t_j})} \le \mu_{c(F_{N-j+1})}\le \mu_{c(F_{N-r_t-\theta^t_j+1})}.
\end{equation}
Thus for any $r_{t-1}-r_t-1\ge s\ge \theta^t_j>r\ge 0$, we have
\begin{equation}\label{ine-mu-cfnrts}
   \mu_{c(F_{N-r_t-s})}\le \mu_{c(F_{N-r_t-\theta^t_j})}\le  \mu_{c(F_{N-j+1})}\le \mu_{c(F_{N-r_t-\theta^t_j+1})}\le \mu_{c(F_{N-r_t-r})}.
\end{equation}
Let $w_j=\sum_{t=1}^i\theta_j^t$ and define 
\begin{equation}\label{equ-def-A-B}
    A=\bigcup_{t=1}^i \{\mu_{c(F_{N-r_t-s})}\colon r_{t-1}-r_t-1\ge s\ge \theta^t_j\},\quad B=\bigcup_{t=1}^i \{\mu_{c(F_{N-r_t-r})}\colon \theta^t_j> r\ge 0\}.
\end{equation} 
Clearly $A$ contains $N-r_i-w_j$ elements and $B$ contains $w_j$ elements. 
Moreover, by \eqref{eq-def-theta} we see that for any $\mu_{c(F_a)}\in A$ and $\mu_{c(F_b)}\in B$, 
\[
c(F_a)\ge c(F_{N-j+1})\ge  c(F_b).
\]
Combining the induction hypothesis \eqref{ine-order},  we deduce
\begin{equation}\label{ine-cFa-cFb}
   f_0^{i-2}(N+1-a)\ge f_0^{i-2}(N+1-b).
\end{equation} 
  Furthermore, by the induction hypothesis \eqref{eq-f0-i-1}, we find that $A\cup B=\{\overline{\delta^{i-1}_{t}}\colon 1\le t\le N-r_i\}$. Together with \eqref{ine-cFa-cFb}, we deduce that
\begin{equation}\label{equ-oveline-set}
    \{\overline{\delta^{i-1}_{t}}\colon w_j+1\le t\le N-r_i\}=A
\end{equation}
and
\begin{equation}\label{eq-ove-del-i-1}
    \{\overline{\delta^{i-1}_{t}}\colon w_j\ge t\ge 1\}=B.
\end{equation}
Thus, $\overline{\delta^{i-1}_{w_j+1}}\in A$ and $\overline{\delta^{i-1}_{w_j}}\in B$. Combining \eqref{eq-tilde-delta-i-1-j-left}, \eqref{ine-mu-cfnrts} and \eqref{equ-def-A-B}, it is clear that
\begin{equation}
    \overline{\delta^{i-1}_{w_j}}\ge \tilde{\delta}^{i-1}_j-\left(N-r_i-w_j\right)\ge \overline{\delta^{i-1}_{w_j+1}}.
\end{equation}
Thus, by the construction of Algorithm Z, when $\psi$ is applied to $\psi^{-1}(\pi,\mu)$,  the partition $\gamma^i$ in Step 2 satisfies
\begin{equation}\label{eq-gammaij-n-ri}
    \gamma^i_j=N-r_i-w_{j}.
\end{equation}
Thus \eqref{eq-fi-delta} gives
\[
\delta^i_{f^{i-1}_0(j)}=\tilde{\delta}^{i-1}_j-\left(N-r_i-w_j\right)=\mu_{c(F_{N-j+1})}.
\]
Moreover, combining \eqref{eq-def-fi} and \eqref{eq-gammaij-n-ri}, we have
\begin{equation}\label{eq-f-i-1(j)}
    f^{i-1}(j)=N-r_i+j-(N-r_i-w_j)=j+w_j.
\end{equation}
This yields \eqref{eq-f0-i-1} and \eqref{eq-f0i-1-j} hold for the case $i$.

We proceed to verify the second assertion \eqref{eq-Pi-i-f0-i-1} holds for the case $i$. From the construction of $\psi$ in \eqref{equ-pi-i+1-w} and the hypothesis \eqref{eq-f0-i-1}, \eqref{eq-Pi-i-f0-i-1}, we deduce that for each $1\le t\le i-1$, $r_{t+1}+1\le j\le r_t$, 
\begin{equation}\label{eq-ell-pi-tf0-t-1}
    \ell(\Pi^{t,f_0^{t-1}(j)})=\ell(\Pi^{i-1,f_0^{i-2}(j)})=t.
\end{equation}
Moreover, using \eqref{eq-def-overlinepi} and \eqref{eq-Pi-i-f0-i-1}, we see that for any $r_{i+1}+1\le j\le r_i$, $1\le k\le i$, 
\begin{equation}\label{eq-Pi-ifi-1(j)}
    \Pi^{i,f^{i-1}(j)}_k=\begin{cases}
        \Pi^{i-1,j}_k+\#\{r_i+w_j+1\le t\le N\colon \ell(\Pi^{i-1,t})=i-k\},&\text{if }1\le k\le i-1;\\
        M+1+\#\{r_i+w_j+1\le t\le N\colon \ell(\Pi^{i-1,t})=0\},&\text{if }k=i.
    \end{cases}
\end{equation}
From \eqref{eq-ell-pi-tf0-t-1} it is clear that for any $0\le k\le i-1$,
\begin{equation}\label{eq-jing-ri-wj-1}
   \#\{r_i+w_j+1\le t\le N\colon \ell(\Pi^{i-1,t})=k\}=\#\{\delta^{i-1}_{f_0^{i-2}(t)}\colon r_i+w_j+1\le f_0^{i-2}(t)\le N, r_{k+1}+1\le t\le r_k\}.
\end{equation}

By \eqref{equ-def-A-B} and the hypothesis \eqref{eq-f0-i-1}, we find that
\[\{\delta^{i-1}_{f_0^{i-2}(t)}\colon r_{k}\ge t\ge \theta^{k+1}_j+r_{k+1}+1\}=A\cap\{\delta^{i-1}_{f_0^{i-2}(t)}\colon r_{k+1}+1\le t\le r_k\}.\]
Using \eqref{equ-oveline-set}, we deduce
\begin{equation*}
    \{\delta^{i-1}_{f_0^{i-2}(t)}\colon r_i+w_j+1\le f_0^{i-2}(t)\le N, r_{k+1}+1\le t\le r_k\}=A\cap\{\delta^{i-1}_{f_0^{i-2}(t)}\colon r_{k+1}+1\le t\le r_k\}.
\end{equation*}

Together with \eqref{eq-jing-ri-wj-1}, we have 
\begin{equation}\label{eq-t-colon}
    \#\{r_i+w_j+1\le t\le N\colon \ell(\Pi^{i-1,t})=k\}=r_k-r_{k+1}-\theta^{k+1}_j.
\end{equation}
Combining \eqref{eq-Pi-t-j}, \eqref{eq-Pi-ifi-1(j)} and \eqref{eq-t-colon}, we arrive at
\begin{equation}
    \Pi^{i,f_0^{i-1}(j)}_k=M+1+N-r_{i+1-k}-\theta^{i+1-k}_j,
\end{equation}
which implies that \eqref{eq-Pi-i-f0-i-1} holds.

Finally, we show that \eqref{ine-order} holds for the case $i$.   Given $r_{i+1}+1\le j_1,j_2\le N$, if $r_{i}+1\le j_1,j_2\le N$, then from the hypothesis we see that
\begin{equation*}
\big(c(F_{N+1-j_1})-c(F_{N+1-j_2})\big)\big(f_0^{i-2}(j_1)-f_0^{i-2}(j_2)\big)>0.
\end{equation*}
Since Algorithm Z preserves the relative order of $\overline{\delta^{i-1}}$ and $\tilde{\delta}^{i-1}$. $\delta_{f_0^{i-2}(j_1)}^{i-1}, \delta^{i-1}_{f_0^{i-2}(j_2)}$ are both in $\overline{\delta^{i-1}}$. We deduce that
\[
\big(f_0^{i-2}(j_1)-f_0^{i-2}(j_2)\big)\big(f_0^{i-1}(j_1)-f_0^{i-1}(j_2)\big)>0.
\]
This yields \eqref{ine-order} holds when $r_{i}+1\le j_1,j_2\le N$. When $r_{i}\ge j_1,j_2\ge r_{i+1}+1$, we can also derive that \eqref{ine-order} holds using the same argument. 

Now we consider the case $r_{i+1}+1\le j_1<r_{i}+1\le j_2\le N$. If $c(F_{N+1-j_1})<c(F_{N+1-j_2})$, then by \eqref{eq-def-theta} and \eqref{equ-oveline-set}, we know that $\overline{\delta^{i-1}_{f_0^{i-2}(j_2)-r_i}}$ is in the set $\{\overline{\delta^{i-1}_{t}}\colon w_{j_1}+1\le t\le N-r_i\}$. Therefore, we have
\begin{equation}
    f_0^{i-2}(j_2)\ge w_{j_1}+1+r_i.
\end{equation}
Combining \eqref{eq-def-fi}, we deduce
\begin{equation}\label{equ-gamma-i-N-f-i-2}
    f^{i-1}_0(j_2)=f_0^{i-2}(j_2)-r_i+{\gamma^{i'}_{N-f^{i-2}_0(j_2)+1}}\ge w_{j_1}+1+{\gamma^{i'}_{N-f^{i-2}_0(j_2)+1}}.
\end{equation}
Moreover, using \eqref{eq-gammaij-n-ri} we see that $\gamma^i_{j_1}=N-r_i-w_{j_1}\ge N-f^{i-2}_0(j_2)+1$, which implies
\begin{equation}\label{ine-est-gamma-i}
    {\gamma^{i'}_{N-f^{i-2}_0(j_2)+1}}\ge j_1.
\end{equation}
Combining \eqref{eq-f-i-1(j)}, \eqref{equ-gamma-i-N-f-i-2} and \eqref{ine-est-gamma-i}, we derive that \eqref{ine-order} holds when $c(F_{N+1-j_1})<c(F_{N+1-j_2})$.

Similarly, if $c(F_{N+1-j_1})>c(F_{N+1-j_2})$, then by \eqref{eq-def-theta} and \eqref{eq-ove-del-i-1}, using the same argument as above, we deduce that
\begin{equation}\label{eq-f-0i-2(j2)le}
    f_0^{i-2}(j_2)\le w_{j_1}+r_i.
\end{equation}
Again by \eqref{eq-gammaij-n-ri}, we know that $\gamma^i_{j_1}=N-r_i-w_{j_1}< N-f^{i-2}_0(j_2)+1$, which implies
\begin{equation}\label{ine-est-gamma-i-1}
    {\gamma^{i'}_{N-f^{i-2}_0(j_2)+1}}< j_1.
\end{equation}
Thus combining \eqref{eq-def-fi}, \eqref{eq-f-0i-2(j2)le} and \eqref{ine-est-gamma-i-1}, we derive that \eqref{ine-order} holds when $c(F_{N+1-j_1})>c(F_{N+1-j_2})$. This completes the proof of this lemma.
\end{proof}

We are now in a position to prove Theorem \ref{thm-CandD}.

{\noindent \it Proof of Theorem \ref{thm-CandD}. } From Lemma \ref{lem-pimupro} and Lemma \ref{lem-psi-1inC}, we know for all $\delta \in C_{M,N}(n)$, $\psi(\delta)=(\pi,\mu)\in D_{M,N}(n)$ and for all $(\pi,\mu)\in D_{M,N}(n)$, $\psi^{-1}(\pi,\mu)=\delta \in C_{M,N}(n)$. 

On the one hand, for any $\delta\in C_{M,N}(n)$, let $\psi(\delta)=(\pi,\mu)$. Setting $j=0$ in Lemma \ref{lem-psipsi-1}, we see that $M+1\in \Pi^{q_\delta,k}$ if and only if $ f_0^{q_\delta-1}(N)<k\le N$. From the construction of $\psi^{-1}$, we know that when we apply $\psi^{-1}$ to $(\pi,\mu)$, the column number of $``F_1"$ is $f_0^{q_\delta-1}(N)$. Similarly, the column number of $``F_i"$ is $f_0^{q_\delta-1}(N-i+1)$ for $1\le i\le N$. This yields \eqref{equ-psi-1-psi}.

On the other hand, utilizing \eqref{eq-o-theta} and \eqref{eq-Pi-i-f0-i-1}, it is clear to see that \eqref{equ-psi-pais-1} holds. This completes the entire proof.\qed

\section{Proof of Theorem \ref{thm-suanfa}}\label{sec-thm-suanfa}
In this section, we present a proof of Theorem \ref{thm-suanfa} using Theorem \ref{thm-CandD}.
We outline the main idea of the bijection $\phi_M$. Given a partition pair $(\alpha, \beta)$ in $A_{M,N}(n)$, we first select certain parts of $\beta$ to form a new partition $\epsilon$. We then apply the inverse bijection $\psi^{-1}$ to the pair $(\epsilon, \alpha)$ to obtain a partition $\eta$. The desired partition $\gamma$ is constructed from $\eta$ together with the remaining parts of $\beta$. The inverse of the map $\phi_M$ is also explicitly described.

\noindent\textit{ Proof of Theorem \ref{thm-suanfa}. } 
Given $(\alpha, \beta)\in A_{M,N}(n)$, by definition we may assume that $\beta=((M+1)^{g_1}, (M+2)^{g_2}, \ldots, (M+N)^{g_N})$ and $|\alpha|+|\beta|=n$. Let $d_i$ and $h_i$ be nonnegative integers satisfying
\begin{equation}\label{eq-gicihi}
    g_{i}=d_{i}(N-i+1)+h_i,
\end{equation}
where $0\le h_i\le N-i$ for $i=1,2,\ldots,N$. Define
\begin{equation}\label{eq-237}
    \epsilon=((M+1)^{h_1}, (M+2)^{h_2}, \ldots, (M+N)^{h_N}).
\end{equation}
Applying the inverse map $\psi^{-1} $ in Theorem \ref{thm-CandD} to $(\epsilon, \alpha)$ we obtain $\eta=(\eta_{1}, \eta_{2}, \ldots, \eta_{N})$ with $0\le \eta_{i}-\eta_{i+1}\le M+N-i$ for $i=1,2,\ldots,N$. Here we also adopt the convention $\eta_{N+1}=0$. 
Next, for $i=1,2,\ldots,N$ define 
  \begin{equation}\label{eq-240}
      \gamma_{i}=\eta_i+\sum_{j=i}^{N}d_{N+1-j}(M+N+1-j),
  \end{equation} 
where $d_{N+1-j}$ is given by \eqref{eq-gicihi}. Consequently, $\gamma=(\gamma_1, \ldots, \gamma_N)$ has $N$ nonnegative parts. We next calculate $|\gamma|$ as follows:
\begin{align}           |\gamma|=\sum_{i=1}^{N}\gamma_i&=|\eta|+\sum_{i=1}^{N}\sum_{j=i}^{N}d_{N+1-j}(M+N+1-j)\nonumber\\
    &=|\epsilon|+|\alpha|+\sum_{j=1}^{N}\sum_{i=1}^{j}d_{N+1-j}(M+N+1-j)\label{eq-mid-epsilon-alpha}\\
    &=\sum_{j=1}^{N}jd_{N+1-j}(M+N+1-j)+|\epsilon|+|\alpha|\nonumber\\
    &=\sum_{i=1}^{N}d_{i}(N+1-i)(M+i)+|\epsilon|+|\alpha|.\label{eq-0909-4}
\end{align} 
Moreover, from the definition of $\beta$, we have 
\begin{align}
    |\beta|&=\sum_{i=1}^{N}g_i(M+i)\nonumber\\
    &=\sum_{i=1}^{N}(d_i(N-i+1)+h_i)(M+i)\nonumber\\
    &=\sum_{i=1}^{N}d_{i}(N-i+1)(M+i)+\sum_{i=1}^{N}h_i(M+i)\nonumber\\
    &=\sum_{i=1}^{N}d_{i}(N-i+1)(M+i)+|\epsilon|.\label{eq-0909-2}
\end{align}
Combining \eqref{eq-0909-4} and \eqref{eq-0909-2}, we deduce $|\gamma|=|\alpha|+|\beta|$, thus $\phi_{M}(\alpha,\beta)=\gamma\in B_{M,N}(n)$.

Conversely, given $\gamma=(\gamma_1, \gamma_2, \ldots, \gamma_N)\in B_{N}(n)$ where $\gamma_1\ge \gamma_2\ge  \cdots\ge \gamma_N\ge 0$, we now construct the inverse map $\phi_M^{-1}$. Suppose for $1\le i\le N$,
\begin{equation}\label{eq-232}
    \gamma_{i}-\gamma_{i+1}=c_{N+1-i}(M+N+1-i)+r_{N+1-i},
\end{equation}
where $c_i$, $r_i$ are nonnegative integers and $0\le r_i\le M+i-1$. We also adopt the convention $\gamma_{N+1}=0$. Define 
$\delta_i=\sum_{j=1}^{N+1-i}r_j$ and let
\begin{equation}\label{eq-235}
        \delta=(\delta_1,\delta_2,\ldots,\delta_{N}).
\end{equation}
Thus, $\delta_{i}-\delta_{i+1}=r_{N+1-i}\le M+N-i$. Here we assume $\delta_{N+1}=0$. Therefore, we may apply
 the map $\psi$ in Theorem \ref{thm-CandD} on $\delta$ and let $(\epsilon, \alpha)=\psi(\delta)$ where 
\begin{equation*}
\epsilon=((M+1)^{t_1}, (M+2)^{t_2},\ldots, (M+N-1)^{t_{N-1}}, (M+N)^{t_{N}})
\end{equation*}
with $0\le t_i\le N-i$ for $i=1,2,\ldots, N$ and $\alpha=(\alpha_1, \alpha_2, \ldots, \alpha_{N})$ with $M\ge \alpha_1\ge \cdots \alpha_{N}\ge 0$. Then we define 
\begin{equation}\label{eq-0910-1}
   {\beta}=((M+1)^{c_1N+t_1},(M+2)^{c_2(N-1)+t_2},\ldots,(M+N)^{c_N+t_N}),
\end{equation}
where $c_i$ is obtained in \eqref{eq-232}. Now define $\phi^{-1}_M(\gamma)=(\alpha,\beta)$. To show that $(\alpha, \beta)\in A_{M,N}(n)$, clearly we only need to verify that $|\alpha|+|\beta|=n$.

On the one hand, we have  
\begin{align}
   n=\sum_{i=1}^{N}\gamma_i&=\sum_{i=1}^{N}i(\gamma_i-\gamma_{i+1})\nonumber\\
    &=\sum_{i=1}^{N}i\left(c_{N+1-i}(M+N+1-i)+r_{N+1-i}\right)\nonumber\\
    &=\sum_{i=1}^{N}ic_{N+1-i}(M+N+1-i)+\sum_{i=1}^{N}ir_{N+1-i}.\label{eq-233}
\end{align}

On the other hand, by the definition of $\delta$,
\begin{equation}\label{eq-0909-7}
|\delta|=\sum_{i=1}^{N}\delta_i=\sum_{i=1}^{N}\sum_{j=1}^{N+1-i}r_{j}=\sum_{i=1}^{N}ir_{N+1-i}.
\end{equation}
Thus we have
\begin{align}\label{eq-te-alpha-beta-alpha-sum}
    |\alpha|+|\beta|&=|\alpha|+\sum_{i=1}^{N}(M+i)(c_{i}(N+1-i)+t_i)\nonumber\\
    &=|\alpha|+\sum_{i=1}^{N}(M+i)t_i+\sum_{i=1}^{N}(M+i)c_i(N+1-i)\nonumber\\
    &=|\alpha|+|\epsilon|+\sum_{i=1}^{N}(M+i)c_i(N+1-i).
    \end{align}
From Theorem \ref{thm-CandD}, we find that $|\alpha|+|\epsilon|=|\delta|$. Combining with \eqref{eq-233} and \eqref{eq-0909-7}, we deduce that
\begin{align}
    |\alpha|+|\beta|&=|\delta|+\sum_{i=1}^{N}ic_{N+1-i}(M+N+1-i)\nonumber\\
    &=\sum_{i=1}^{N}ir_{N+1-i}+ \sum_{i=1}^{N}ic_{N+1-i}(M+N+1-i)\nonumber\\
    &=n.
\end{align}
Thus we deduce $\phi^{-1}_{M}(\gamma)=(\alpha, \beta)\in A_{M,N}(n)$.

By Theorem \ref{thm-suanfa}, it is routine to check that $\phi_M^{-1}$ is indeed the inverse map of $\phi_M$. 
\qed

For example, let $N=10$, $M=2$.  given $(\alpha, \beta)\in A_{2,10}(330)$ where 
\[
\alpha=(2,1,1,1,1,1,1,1,0,0)
\]
and
\[
\beta=(3^2,4^3, 5^3,6^2,  7^1,8^{10},  9^7, 10^8, 11^2,12^2).
\]
By \eqref{eq-gicihi} we know $d_{10}=2$, $d_{9}=1$, $d_{8}=2$, $d_{7}=1$, $d_{6}=2$, $d_{5}=d_{4}=d_{3}=d_{2}=d_{1}=0$ and $h_1=2$, $h_2=3$, $h_3=3$, $h_4=2$, $h_5=1$, $h_6=0$, $h_7=3$, $h_8=2$, $h_9=0$, $h_{10}=0$. Thus by \eqref{eq-237} we obtain
\[
\epsilon=(3^2,4^3, 5^3,6^2, 7^1, 9^3, 10^2).
\]
 Applying $\psi^{-1}$ on $(\epsilon, \alpha)$, we have 
\[
\eta=(28,26,20,12,6,6,5,3,1,1)\in C_{2,10}(108).
\]
From \eqref{eq-240}, we deduce
\[
\gamma=(108,82,65,37,22,6,5,3,1,1)
\]
and it can be easily checked that $\gamma\in B_{N}(330)$.

Conversely, given
\[
\gamma=(108,82,65,37,22,6,5,3,1,1)\in B_{N}(330),
\]
by \eqref{eq-232} we deduce that $c_{10}=2$, $c_{9}=1$, $c_{8}=2$, $c_{7}=1$, $c_{6}=2$, $c_{5}=c_{4}=c_{3}=c_{2}=c_{1}=0$, and $r_{10}=2$, $r_{9}=6$, $r_{8}=8$, $r_{7}=6$, $r_{6}=0$, $r_{5}=1$, $r_{4}=2$, $r_{3}=2$, $r_{2}=0$, $r_{1}=1$. 
Then by the definition of $\delta$ we have 
\[
\delta=(28,26,20,12,6,6,5,3,1,1)\in C_{2,10}(108).
\]
Using the injection $\psi$ in Theorem \ref{thm-CandD}, we get 
\[
\alpha=\delta^4=(2,1,1,1,1,1,1,1,0,0),
\]
and
\[
\epsilon=(3^2,4^3, 5^3,6^2, 7^1, 9^3, 10^2).
\]
 By \eqref{eq-0910-1}
\[
\beta=(3^2,4^3, 5^3,6^2,  7^1,8^{10},  9^7, 10^8, 11^2,12^2).
\]
It is easy to check that $(\alpha, \beta)\in A_{2,10}(330)$.  

\section{Proof of Theorem \ref{thm-B1989-lem3.4}}\label{sec-thm-B1989-lem3.4}
\begin{figure}
    \centering
    \begin{tikzpicture}[scale=.3]
    \draw[thick] (0,0) -- (14,0) -- (14,-12) -- (0,-12) -- cycle;
    \draw[thick] (14,0)--(19,0)--(19,-8)--(14,-12)--cycle;
    \draw[thick] (19,0)--(22,0)--(22,-6)--(19,-7)--cycle;
     \foreach \x/\y in {22.5/-2, 23.5/-2, 24.5/-2} {
    \fill[black] (\x,\y) circle (2pt);
    }
    \draw[thick] (25,0)--(28,0)--(28,-3)--(25,-4)--cycle;
    \draw[thick] (28,0)--(32,0)--(28,-2)--cycle;
    \draw[thick] (0,-12)--(9,-12)--(9,-19)--(0,-19)--cycle;
    \foreach \x/\y in {3/-19.5, 3/-20.5, 3/-21.5} {
    \fill[black] (\x,\y) circle (2pt);
    }
    \draw[thick] (0,-22)--(6,-22)--(6,-26)--(0,-26)--cycle;
    \draw[thick] (0,-26)--(4,-26)--(4,-28)--(0,-28)--cycle;
    \draw[thick] (0,-28)--(4,-28)--(0,-32)--cycle;
    \draw[thick] (9,-12)--(14,-12)--(9,-19)--cycle;
    \draw[thick] (4,-26)--(6,-26)--(4,-28)--cycle;
    \node[] at (7,-6) {$R^1$};
    \node[] at (4.5,-15.5) {$R^2$};
    \node[] at (3,-24) {$R^{k-1}$};
    \node[] at (2,-27) {$R^k$};
    \node[] at (11,-14) {$b^{1}$};
    \node[scale=0.6] at (4.8,-26.5) {$b^{k-1}$};
    \node[] at (16.4, -3.5) {$r^1$};
    \node[] at (20.5, -3) {$r^2$};
    \node[] at (26.5,-1.5) {$r^{k-1}$};
    \node[] at (29.5, -0.6) {$r^k$};
    \node[] at (1, -29) {$b^k$};
\end{tikzpicture}
    \caption{Decompose $\lambda$ into $3k$ partitions.}
    \label{fig-2mdurfee}
 \end{figure}
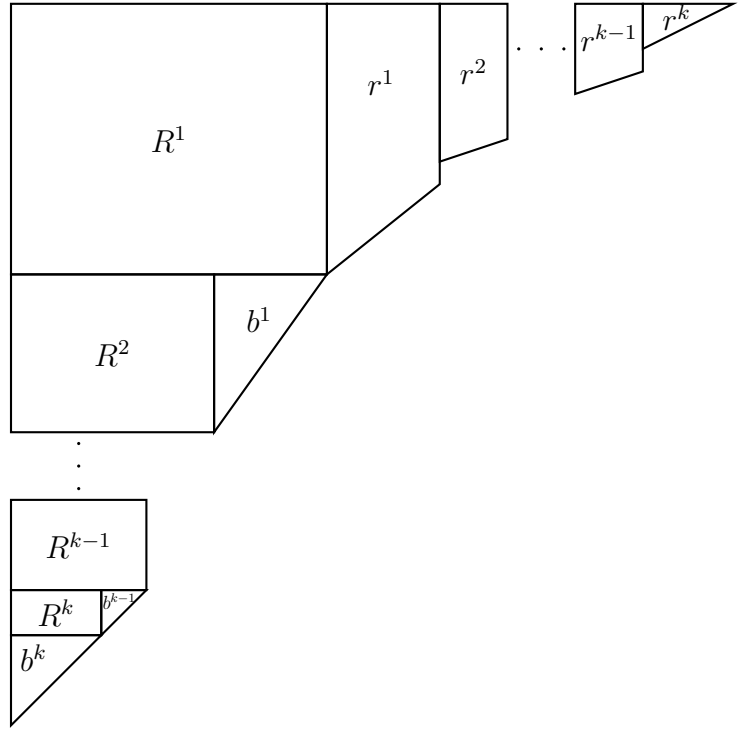

In this section, we provide an alternative proof of Theorem \ref{thm-B1989-lem3.4}. To this end, for any $(\lambda, \delta) \in R_{k,m}(n)$, we begin with the case $m \ge 0$ and decompose $\lambda$ into $3k$ partitions, namely $r^1,\ldots, r^{k-1},r^{k},b^1,\ldots,b^{k-1},b^{k},R^1,\ldots,R^k$, as illustrated in Figure \ref{fig-2mdurfee}. We then apply $\phi_{n_{i+1}}$ to the pair $({b^i}', r^i)$ to obtain $\nu^i$ for $1\le i\le k-1$. Next, using the inverse map $\phi^{-1}_{n_{i+1}+m}$ on $\nu^i$, we transform $\nu^i$ into a pair of partitions $(\overline{b}^i, \overline{r}^i)$. Furthermore, we apply $\phi^{-1}_{n_k + 2m}$ to $r^k$ to yield another pair of partitions $(r^{k,1}, r^{k,2})$. The remainder of the argument follows essentially the same procedure as in \cite{Bressound-1989-Generlaized}; we briefly outline the steps for completeness. The case $m < 0$ is handled similarly to the case $m \ge 0$ and we omit the details. 

{\noindent\textit{Proof of Theorem \ref{thm-B1989-lem3.4}. }}
There are two cases. 

\noindent \textbf{Case I}. $m\ge 0$. Given $(\lambda, \delta) \in R_{k,m}(n)$, recall that the $m$-Durfee rectangle of a partition $\lambda$, introduced by Gordon and Houten \cite{Gordon-Houten-1968}, is defined as the largest $(m + j) \times j$ rectangle contained in the Ferrers diagram of $\lambda$. Note that an $m$-Durfee rectangle reduces to a Durfee square when $m = 0$. 

Let $(n_1 + 2m)\cdot n_1$ denote the $2m$-Durfee rectangle of $\lambda$, labeled as $R^1$ in Figure \ref{fig-2mdurfee}. Iteratively, for $2 \le i \le k$, let $(n_i + 2m) \cdot n_i$ be the $2m$-Durfee rectangle of the subpartition of $\lambda$ consisting of all parts not exceeding $n_{i-1} + 2m$; this rectangle is marked as $R^i$ in Figure \ref{fig-2mdurfee}.

As illustrated in Figure \ref{fig-2mdurfee}, for $1 \le i \le k - 1$, let $b^i$ denote the partition situated below $R^i$ and to the right of $R^{i+1}$. Furthermore, let $b^k$ denote the subpartition of $\lambda$ consisting of all parts not exceeding $n_k + 2m$. For the right part of $R^1$, denoted by $\mathcal{T}$, we further divide $\mathcal{T}$ into $k$ parts as follows: define $r^i$ to be the partition consisting of all columns in the conjugate of $\mathcal{T}$ whose lengths lie in the interval $[n_{i+1} + 1, n_i]$ for $1 \le i \le k - 1$. Additionally, let $r^{k}$ denote the partition consisting of all columns with length not exceeding $n_k$.

For example, if $\lambda=(11,10,10,9,8,7,6,5,5,4,3,1,1)$, $m=1$, $k=3$, then by the definition of $\delta$ in \eqref{eq-defofdelta}, we see that  $\delta=(4)$. Thus $(\lambda,\delta)\in R_{3,1}(84)$. Moreover we know $n_1=5$, $n_2=3$, $n_3=2$, $b^1=(2,1)$, $b^2=(1)$, $b^3=(3,1,1)$. Furthermore, $\mathcal{T}=(4,3,3,2,1)$, and $\mathcal{T}'=(5,4,3,1)$. Thus $r^1=(5,4)$, $r^2=(3)$ and  $r^3=(1)$ (See Figure \ref{fig-mdur}). 
\begin{figure}
    \centering
     \begin{tikzpicture}[scale=.4][font=\footnotesize]
 \fill foreach \z [count=\y] in {11,10,10,9,8,7,6,5,5,4,3,1,1}
  {foreach \x in {1,...,\z} 
  {(\x,-\y) circle[radius=3pt]}};
\draw[black,thick] (0.6,-0.6) -- (7.4,-0.6) -- (7.4,-5.4) -- (0.6,-5.4) -- cycle;
\draw[black,thick] (0.6,-5.6) -- (5.4,-5.6) -- (5.4,-8.4) -- (0.6,-8.4) -- cycle;
\draw[black,thick] (0.6,-8.6) -- (4.4,-8.6) -- (4.4,-10.4) -- (0.6,-10.4) -- cycle;
\draw[black,thick] (0.6,-10.6) -- (4,-10.6) -- (0.6,-14) -- cycle;
\draw[black, thick] (7.6, -0.6) -- (9.4, -0.6) -- (9.4, -4.4) -- (7.6, -5.4)--cycle;
\draw[black, thick] (9.6, -0.6) -- (10.4, -0.6) -- (10.4, -3.4) -- (9.6, -3.6)--cycle;
\draw[black,thick] (10.6,-0.6) -- (13,-0.6) -- (10.6,-1.6) -- cycle;
\draw[black,thick] (5.6,-5.6) -- (7.4,-5.6) -- (5.6,-8.4) -- cycle;
\draw[black,thick] (4.6,-8.6) -- (5.4,-8.6) -- (4.6,-10) -- cycle;
\end{tikzpicture}
    \caption{Illustration of $2m$-Durfee rectangle of $\lambda$.}
    \label{fig-mdur}
\end{figure}
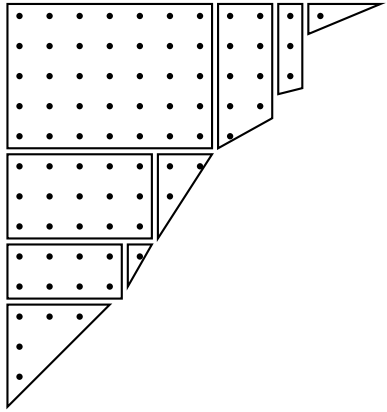

\begin{figure}
    \centering
    \begin{tikzpicture}[scale=.6][font=\footnotesize]
 \fill foreach \z [count=\y] in {11,9,9,8,7,7,6,6,5,5,3,3,3}
  {foreach \x in {1,...,\z} 
  {(\x,-\y) circle[radius=3pt]}};
\draw[black,thick] (0.6,-0.6) -- (6.4,-0.6) -- (6.4,-6.4) -- (0.6,-6.4) -- cycle;
\draw[black,thick] (0.6,-6.6) -- (4.4,-6.6) -- (4.4,-10.4) -- (0.6,-10.4) -- cycle;
\draw[black,thick] (0.6,-10.6) -- (3.4,-10.6) -- (3.4,-13.4) -- (0.6,-13.4) -- cycle;
\draw[black,thick] (6.6,-0.6) -- (12,-0.6) -- (6.6,-7) -- cycle;
\node[thick] at (0, -3.5) {$D^1$};
\node[thick] at (0, -8.5) {$D^2$};
\node[thick] at (0, -12) {$D^3$};
\node[thick] at (10, -4) {$\mathcal{R}'$};
\end{tikzpicture}
    \caption{{\bf Step 8} to construct $\overline{\alpha}$.}
    \label{fig-overlinealpha}
\end{figure}

From the above decomposition, it is readily seen that for $1\le i\le k-1$, $b^i$ is a partition with at most $n_{i+1}$ parts, each of size at most $n_i-n_{i+1}$; $b^k$ is a partition with each part not exceeding $n_k+2m$; for $1\le i\le k$, $r^i$  is a partition with each part lying in the interval $[n_{i+1}+1,n_i]$, here we adopt the convention that $n_{k+1}=0$.

We now describe the map $\chi$, which consists of the following ten steps.

{\bf \noindent Step 1.} For $1\le i\le k-1$, apply $\phi_{n_{i+1}}$ to the pair $({b^i}',r^i)$ to obtain a partition $\nu^i$. By Theorem \ref{thm-suanfa},  $\nu^i$ is a partition with at most $n_i-n_{i+1}$ parts. In this  example, we obtain $\nu^1=(10,2)$, $\nu^2=(4)$.

{\bf \noindent Step 2.} For $1\le i\le k-1$, apply $\phi^{-1}_{n_{i+1}+m}$ to $\nu^i$ to obtain a pair of partitions $(\overline{b}^i,\overline{r}^i)$. Again, by Theorem \ref{thm-suanfa},  $\overline{b}^i$ is a partition with at most $n_i-n_{i+1}$ parts, each not exceeding $n_{i+1}+m$, and $\overline{r}^i$ is a partition with all parts lying in $[n_{i+1}+m+1,n_{i}+m]$.
In this example, we obtain $\overline{b}^1=(4,2)$, $\overline{b}^2=\emptyset$, $\overline{r}^1=(6)$, $\overline{r}^2=(4)$.

{\bf \noindent Step 3.} Apply $\phi^{-1}_{n_k+2m}$ to $r^k$ to obtain a partition pair $(r^{k,1},r^{k,2})$, where $r^{k,1}$ is a partition with at most $n_k$ parts, each at most $n_k+2m$; and $r^{k,2}$ is a partition with each part lying in $[n_k+2m+1,2n_k+2m]$. In this example, we obtain $r^{k,1}=(1)$, $r^{k,2}=\emptyset$.

{\bf \noindent Step 4.} We divide the partition $b^k$ into two partitions $b^{k,1}$ and $b^{k,2}$ where $b^{k,1}$ is the partition with all parts not exceeding $n_k+m$ and $b^{k,2}$ is the partition with all parts lying in $[n_k+m+1, n_k+2m]$. In this example, we obtain $b^{k,1}=(3,1,1)$, $b^{k,2}=\emptyset$.

{\bf \noindent Step 5.} For the partition $\delta$, we define
\begin{equation}\label{eq-defthm1.4gammai}
\gamma_i=\delta_i-(2(m-i)+1)k
\end{equation}
where $1\le i\le m$ to get $\gamma=(m,m-1, \ldots, 1)$. In this example, we obtain $\gamma=(1)$.

{\bf \noindent Step 6.} For each $(n_{i}+2m) \cdot n_{i}$ ($1\le i\le k$) $2m$-Durfee rectangle in $\lambda$, we add $m^2$ to reshape it into Durfee Square $D^i$ with length $n_i+m$. In this example, we obtain $D^1=(6^6)$, $D^2=(4^4)$, $D^3=(3^3)$.

{\bf \noindent Step 7.} For $1\le i\le k-1$, we put $\overline{r}^{i+1}$ under $\overline{r}^{i}$ and put $b^{k,1}$ under $\overline{r}^{k-1}$ to get a new partition $\mathcal{R}$ with all parts not exceeding $n_1+m$. In this example, we obtain $\mathcal{R}=(6,4,3,1,1)$.

{\bf \noindent Step 8.} For $1\le i\le k-1$, we put $(\overline{b}^{i})'$ under $D^i$ and to the right of $D^{i+1}$ and put $\mathcal{R}'$ to the right of $D^1$ to get partition $\overline{\alpha}$. In this example, we obtain $\overline{\alpha}=(11,9,9,8,7,7,6,6,5,5,3,3,3)$ (See Figure \ref{fig-overlinealpha}).

{\bf \noindent Step 9.} Put $b^{k,2}$ under $r^{k,2}$ to get a new partition $\overline{\beta}$ with all parts lying in $[n_k+m+1, 2n_k+2m]$. In this example, we obtain $\overline{\beta}=\emptyset$.

{\bf \noindent Step 10.} Finally, conjugate $r^{k,1}$ to get partition $\overline{\xi}$ with all parts not exceeding $n_k$ and at most $n_k+2m$ parts. In this example, we obtain $\overline{\xi}=(1)$.

After the above steps, for $1\le i\le k$, let $s_i=n_i+m$, we get 
\begin{equation*}
    (\alpha, \beta, \gamma, \xi):=\chi(\lambda, \delta)=(\overline{\alpha}, \overline{\beta}, \gamma, \overline{\xi})
\end{equation*} 
in which $\alpha=\overline{\alpha}$ is a partition with all parts not less than $s_k$ and the length of the $i$-th Durfee squares in $\alpha$ equals $s_i$, $\beta=\overline{\beta}$ with all parts lying in $[s_k+1, 2s_k]$ and $\xi=\overline{\xi}$ with all parts not exceeding $s_k-m$ and at most $s_k+m$ parts.   Moreover, from \eqref{eq-defthm1.4gammai}, we have 
\[
|\gamma|=|\delta|-\sum_{i=1}^{m}2(m-i)+1=|\delta|-km^2,
\]
and
\[
|\alpha| + |\beta| + |\gamma| + |\xi|=|\delta|-km^2+km^2+|\lambda|=n.
\]

Thus, $(\alpha, \beta, \gamma, \xi)\in S_{k,m}(n)$. Since each step above is  reversible, we see that $\chi$ is a bijection between $R_{k,m}(n)$ and $S_{k,m}(n)$. In this example, we obtain 
\[
(\alpha, \beta, \gamma, \xi)=((11,9,9,8,7,7,6,6,5,5,3,3,3), \emptyset, (1), (1))
\]
and it can be checked that $(\alpha, \beta, \gamma, \xi)\in S_{3,1}(84)$.

\textbf{Case II}. $m<0$. This case is closely analogous to \textbf{Case I}. Instead of considering a $2m$-Durfee rectangle $R^i$, we take $R^i$ to be an $(n_i-2m) \cdot n_i$ rectangle, i.e. a $-2m$-Durfee rectangle. Moreover, in Step 6, we define $s_i=n_i-m$ rather than $s_i=n_i+m$. The remainder of the proof follows exactly as in \textbf{Case I}, and we omit the details.
\qed

\begin{rem}
Note that only Steps 1--3 make use of Theorem \ref{thm-suanfa}, which distinguishes this approach from the original proof in \cite{Bressound-1989-Generlaized}. The rest of the argument essentially follows the same lines as in \cite{Bressound-1989-Generlaized}.
\end{rem}

\section{Proof of Theorem \ref{thm-yingyong}}\label{sec-yingyongNkmn}

This section is aimed to give a combinatorial proof of Theorem \ref{thm-yingyong}.
We first recall the combinatorial interpretation of $N_k(m,n)$ which was first introduced by Garvan \cite{Garvan1994generalizations} and let $Q_k(m,n)$ denote the set of partitions counted by $N_k(m,n)$. Then we introduce the definition of $P_k(m,n)$, which is the set of  $(2k-1)$-tuple of partitions of $n$. Using Theorem \ref{thm-suanfa}, we show that there is a bijection $\eta$ between $P_k(m,n)$ and $Q_k(m,n)$. We then partition the set of $P_k(m,n)$ into $16$ disjoint subsets, namely $P^{i}_{k}(m,n)$ ($1\le i\le 16$). Consequently, $15$ disjoint subsets of $P_k(m,n+1)$, namely $P_{k,i}(m,n+1)$ ($1\le i\le 15$), will be listed. Then we build $15$ injections from $P^{i}_{k}(m,n)$ to $P_{k,i}(m,n+1) (1\le i\le 15)$ show that $P^{16}_k(m,n)$ is empty except when $n= |m|+k-1$ or $(m,k,n)= (0,3,8)$. This yields a combinatorial proof of Theorem \ref{thm-yingyong}.

We now recall the combinatorial interpretation of $N_k(m,n)$ given by Garvan \cite{Garvan1994generalizations}. For a partition $\pi$, let $d_1, d_2, \ldots$ to be the sizes of the successive Durfee squares of $\pi$. We denote $d_{\ell}=0$ if the number of successive Durfee squares of $\pi$ is less than $\ell$. The $k$-rank, $r_{k}(\pi)$, is the number of columns in the Ferrers graph of $\pi$ which lie to the right of the first Durfee square and whose length $\le d_{k-1}$ minus the number of parts of $\pi$ that lie below the $(k-1)$-th Durfee square.  To be specific, let $\alpha^{T}$ to denote the partition lying in the right of the first Durfee square with $\alpha_1 \le d_{k-1}$. Let $\beta$ denote the partition below the $(k-1)$-th Durfee square. Then $r_{k}(\pi)=\ell(\alpha)-\ell(\beta)$. We denote that $r_k(\pi)=0$ if $d_{k-1}=0$. 

\begin{figure}
    \centering
    \begin{tikzpicture}[scale=.6][font=\footnotesize]
 \fill foreach \z [count=\y] in {8,6,6,4,3,3,2,2,1}
  {foreach \x in {1,...,\z} 
  {(\x,-\y) circle[radius=3pt]}};
\draw[black,thick] (0.6,-0.6) -- (4.4,-0.6) -- (4.4,-4.4) -- (0.6,-4.4) -- cycle;
\draw[black,thick] (0.6,-4.6) -- (2.4,-4.6) -- (2.4,-6.4) -- (0.6,-6.4) -- cycle;
\draw[black,thick] (0.6,-6.6) -- (2.4,-6.6) -- (2.4,-8.4) -- (0.6,-8.4) -- cycle;
\draw[black,thick] (0.6,-8.6) -- (1.4,-8.6) -- (1.4,-9.4) -- (0.6,-9.4) -- cycle;
\end{tikzpicture}
    \caption{Illustration of successive Durfee squares of $\pi$.}
    \label{fig-rkpi}
\end{figure}
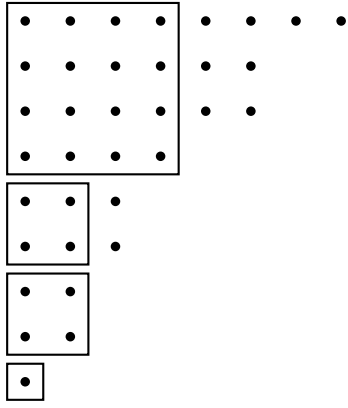

For example, given partition $\pi=(8,6,6,4,3,3,2,2,1)$, we know $d_1=4$, $d_2=d_3=2$, $d_4=1$(See Fig. \ref{fig-rkpi}). By simply calculating, we have $r_2(\pi)=-1$,  $r_3(\pi)=-1$, $r_4(\pi)=1$, $r_5(\pi)=2$, $r_i(\pi)=0$ for $i\ge 6$. 

Let $Q_k(m,n)$ denote the set of partitions $\pi$ of $n$ that have at least $k-1$ successive Durfee squares and $r_k(\pi)=m$. Garvan \cite{Garvan1994generalizations} showed that
\[\# Q_k(m,n)=N_k(m,n).\]

We next introduce the definition of $P_k(m,n)$.

\begin{defi}\label{def-Pkmn}
   Given $k\ge 3$, $n\ge |m|+k-1$ and $m\in\mathbb{Z}$, let $P_k(m,n)$ denote the set of $(2k-1)$-tuple of partitions 
    \[
    \Delta=(\alpha, \beta, \gamma^1, \gamma^2, \ldots, \gamma^{k-2}, \varpi^{1}, \ldots, \varpi^{k-1}),
    \]
    which satisfies the following restrictions:
    \begin{itemize}
        \item[(1)] $\varpi^{i}$($1\le i\le k-1$) is a partition such that each part equals $\ell(\varpi^{i})$. In other words, $\varpi^i=(d_i^{d_i})$, where $d_i=\ell(\varpi^{i})$;
        \item[(2)] $\ell(\varpi^{i})\ge \ell(\varpi^{i+1})\ge 1$ for $1\le i\le k-2$;
        \item[(3)] $\alpha$ is a partition such that all parts $\le \ell(\varpi^{k-1})$;
        \item[(4)] $\beta$ is a partition such that all parts $\le \ell(\varpi^{k-1})$;
        \item[(5)] $\ell(\alpha)-\ell(\beta)=m$;
        \item[(6)] $\gamma^{i}$($1\le i \le k-2$) are the partitions whose length $\le \ell(\varpi^{i}) - \ell(\varpi^{i+1})$;
        \item[(7)] $|\alpha|+|\beta|+\sum_{i=1}^{k-2}|\gamma^{i}|+\sum_{i=1}^{k-1}|\varpi^{i}|=n$.
    \end{itemize}
\end{defi}

For the sake of simplicity,  here and throughout this section, for 
\[\Delta=(\alpha, \beta, \gamma^1, \gamma^2, \ldots, \gamma^{k-2}, \varpi^{1}, \ldots, \varpi^{{k-1}})\in P_k(m,n),\]
we will use $d_i$ to denote $\ell(\varpi^{i})$. We will also frequently use $$\tilde{\Delta}=(\tilde{\alpha}, \tilde{\beta}, \tilde{\gamma}^1, \tilde{\gamma}^2, \ldots, \tilde{\gamma}^{k-2}, \tilde{\varpi}^{1}, \ldots, \tilde{\varpi}^{{k-1}})$$
to denote another partition tuples in $P_k(m,n)$. We write $\tilde{d}_i=\ell(\tilde{\varpi}^{i})$.

We are now in a position to describe the bijection between $P_k(m,n)$ and $Q_k(m,n)$.

\begin{thm}\label{thm-QtoP}
    There is a one-to-one correspondence $\eta$ between the set $P_k(m,n)$ and the set $Q_k(m,n)$.
\end{thm}
\begin{proof}
    For any $\pi\in Q_{k}(m,n)$ and $1\le j\le k-1$, let $d_j$ be the size of the $j$-th successive Durfee squares of $\pi$. From the definition of $Q_k(m,n)$, we see that $d_{k-1}\ge 1$. For $1\le i \le k-2$, let $\pi^{b_{i}}$ denote the partition consisting of columns below the $i$-th Durfee square and the right of the $(i+1)$-th Durfee square.
    Moreover, let $\pi^{r_i}$ be the partition consisting of the columns to the right of the first Durfee square such that the length of each column lies between $d_{i}$ and $ d_{i+1}+1$. Similarly, we use $\alpha$ to denote the partition consisting of the columns to the right of the first Durfee square with the length of each column not exceeding $d_{k-1}$, and use $\beta$ to denote the partition that contains the rows below the $k-1$-th Durfee square. We also denote these $k-1$ successive Durfee squares by $\varpi^1,\ldots,\varpi^{k-1}$ respectively (see Figure \ref{fig-pi} for an illustration). It is clear that $r_k(\pi)=\ell(\alpha)-\ell(\beta)=m$.

Applying the map $\phi_{d^{i+1}}$ defined in Theorem \ref{thm-suanfa} on $(\pi^{b_{i}}, \pi^{r_{i}})$ we can obtain a partition $\gamma^{i}$. Then $\gamma^{i}$ is a partition with at most $d_{i}-d_{i+1}$ parts. Now define 
    \[\eta(\pi)=(\alpha,\beta,\gamma^1,\ldots,\gamma^{k-2},\varpi^1,\ldots,\varpi^{k-1}).\]
    It is routine to check that $\eta(\pi)\in P_k(m,n)$. Moreover, it is easy to check that the above map is revertible since $\phi_{d^{i+1}}$ is revertible. 
\end{proof}

For example, let $k=3$, $m=3$, $n=48$ $\pi=(12,9,8,6,5,4,3,1)\in Q_3(3,48)$, then we know $d_{1}=5$, $d_{2}=2$, $\pi^{b_1}=(2,1,0)$, $\pi^{r_1}=(4,3,3)$, $\alpha=(2,1,1,1)$ and $\beta=(1)$ (see Figure \ref{fig-exam}). Using the map $\phi_2$ on $(\pi^{b_1}, \pi^{r_1})$, we get the $\gamma^{1}=(7,4,2)$. Thus by Theorem \ref{thm-QtoP}, we deduce $\Delta:=\eta(\pi)=((2,1,1,1), (1), (7,4,2), 5^5, 2^2)\in P_{3}(3,48)$.

\vspace{2em}
\begin{figure}[htbp]
  \centering
  \begin{minipage}[t]{0.48\textwidth}
    \centering
\begin{tikzpicture}[scale=.4]
    \draw[thick] (0,0) -- (7,0) -- (7,-7) -- (0,-7) -- cycle;
    \draw[thick] (7,0) -- (10,0) -- (10,-6) -- (7,-7) -- cycle;
    \foreach \x/\y in {10.5/-1, 11.5/-1, 12.5/-1} {
    \fill[black] (\x,\y) circle (2pt);
    }
    \draw[thick] (13,0) -- (16,0) -- (16,-3) -- (13,-4) -- cycle;
    \draw[thick] (16,0) -- (19,0) -- (16,-2) -- cycle;
    \draw[thick] (0,-7) -- (5,-7) -- (5,-12) -- (0,-12) -- cycle;
    \draw[thick] (5,-7) -- (7,-7) -- (5,-12) -- cycle;
     \foreach \x/\y in {1/-12.5, 1/-13.5, 1/-14.5} {
    \fill[black] (\x,\y) circle (2pt);
    }
    \draw[thick] (0,-15) -- (4,-15) -- (4,-19) -- (0,-19) -- cycle;
    \draw[thick] (0,-19) -- (2,-19) -- (2,-21) -- (0,-21) -- cycle;
    \draw[thick] (2,-19) -- (4,-19) -- (2,-21)-- cycle;
    \draw[thick] (0,-21) -- (2,-21) -- (0,-24) -- cycle;
    \draw[decorate, decoration={brace, amplitude=5pt, mirror}, thick] (0,0) -- (0,-7);
    \node[] at (-1, -3.5) {$d_{1}$};
    \draw[decorate, decoration={brace, amplitude=5pt, mirror}, thick] (0,-7) -- (0,-12);
    \node[] at (-1, -9.5) {$d_{2}$};
    \draw[decorate, decoration={brace, amplitude=5pt, mirror}, thick] (0,-15) -- (0,-19);
    \node[] at (-1.5, -17) {$d_{k-2}$};
      \draw[decorate, decoration={brace, amplitude=5pt, mirror}, thick] (0,-19) -- (0,-21);
    \node[] at (-1.5, -20) {$d_{k-1}$};
     \draw[decorate, decoration={brace, amplitude=5pt, mirror}, thick] (10,0) -- (10,-6);
     \node[scale=0.6] at (8.4, -3) {$\ge d_2+1$};
     \draw[decorate, decoration={brace, amplitude=5pt, mirror}, thick] (13,0) -- (13,-4);
     \node[scale=0.5] at (11.7, -2) {$\le d_{k-2}$};
     \draw[decorate, decoration={brace, amplitude=5pt, mirror}, thick] (16,0) -- (16,-3);
     \node[scale=0.4] at (14.5, -1.5) {$\ge d_{k-1}+1$};
    \draw[-latex, thick] (19,-2)--(17.5,-0.5);
    \node[] at (19.5,-2.5) {$\alpha$};
    \draw[-latex, thick] (16,-5)--(14.5,-3);
    \node[] at (16.5,-5.5) {$\pi^{r_{k-2}}$};
    \draw[-latex, thick] (10,-8)--(8.5,-5);
    \node[] at (10.5,-8.5) {$\pi^{r_{1}}$};
    \draw[-latex, thick] (7,-11)--(5.5,-8.5);
    \node[] at (7.5,-11.5) {$\pi^{b_{1}}$};
    \draw[-latex, thick] (4,-21)--(2.5,-19.5);
    \node[] at (4.5,-21.5) {$\pi^{b_{k-2}}$};
    \draw[-latex, thick] (2,-23)--(0.5,-22);
    \node[] at (2.5,-23.5) {$\beta$};
\end{tikzpicture}
    \caption{Example of $k-1$ successive Durfee Squares.}
    \label{fig-pi}
  \end{minipage}
  \hfill
  \begin{minipage}[t]{0.48\textwidth}
    \centering
\begin{tikzpicture}[scale=.6][font=\footnotesize]
 \fill foreach \z [count=\y] in {12,9,8,6,5,4,3,1}
  {foreach \x in {1,...,\z} 
  {(\x,-\y) circle[radius=3pt]}};
\draw[black,thick] (1,-1) -- (5,-1) -- (5,-5) -- (1,-5) -- cycle;
\draw[black,thick] (1,-6) -- (2,-6) -- (2,-7) -- (1,-7) -- cycle;
\draw[-latex,black,thick] (6,-1) -- (6,-4.5);
\node[] at (6, 0) {$\pi^{r_1}_1$};
\draw[-latex,black,thick] (7,-1) -- (7,-3.5);
\node[] at (7, 0) {$\pi^{r_1}_2$};
\draw[-latex,black,thick] (8,-1) -- (8,-3.5);
\node[] at (8, 0) {$\pi^{r_1}_3$};
\draw[-latex,black,thick] (9,-1) -- (9,-2.5);
\node[] at (9, 0) {$\alpha_1$};
\draw[-latex,black,thick] (10,-1) -- (10,-1.5);
\node[] at (10, 0) {$\alpha_2$};
\draw[-latex,black,thick] (11,-1) -- (11,-1.5);
\node[] at (11, 0) {$\alpha_3$};
\draw[-latex,black,thick] (12,-1) -- (12,-1.5);
\node[] at (12, 0) {$\alpha_4$};
\draw[-latex, black, thick] (3,-6) -- (3, -7.5);
\node[] at (3, -8) {$\pi^{b_1}_{1}$};
\draw[-latex, black, thick] (4,-6) -- (4, -6.5);
\node[] at (4, -7) {$\pi^{b_1}_{2}$};
\draw[-latex, black, thick] (5,-6) -- (5, -6.5);
\node[] at (5, -7) {$\pi^{b_1}_{3}$};
\draw[-latex, black, thick] (1,-8) -- (1.5,-8);
\node[] at (0,-8) {$\beta_1$};
\end{tikzpicture}
    \caption{Example of $2$ successive Durfee Squares.}
    \label{fig-exam}
  \end{minipage}
\end{figure}

Since $N_k(m,n)=N_k(-m,n)$, using Theorem \ref{thm-QtoP} we see that Theorem \ref{thm-yingyong} is a direct consequence of the following theorem.

\begin{thm}\label{thm-PnPn+1}
    For $k\ge 3$, $n\ge k-1$, $m\ge 0$,  there is an injection $\sigma$ from the set $P_{k}(m,n)$ to $P_{k}(m,n+1)$, except when $(m,k,n)=(m,k,m+k-1)$ or $(0,3,8)$.
\end{thm}

To construct this injection, we partition $P_{k}(m,n)$ into $16$ disjoint subsets $P^{i}_{k}(m,n)$ ($1\le i \le 16 $) as follows.
\begin{itemize}
    \item[(1)] $P^{1}_{k}(m,n)=\{\Delta \in P_{k}(m,n)\colon d_{1}\ne d_{k-1}\}$;
    \item[(2)] $P^{2}_{k}(m,n)=\{\Delta \in P_{k}(m,n)\colon d_{1}=d_{k-1}=1, n\ge m+k\}$;
    \item[(3)] $P^{3}_{k}(m,n)=\{\Delta \in P_{k}(m,n)\colon k\ge 4, d_{1}=d_{k-1} \ge 2, \alpha  =\emptyset \}$;
    \item[(4)] $P^{4}_{k}(m,n)=\{\Delta \in P_{k}(m,n)\colon k=3, d_{1}=d_{2}\ge 3, \alpha  =\emptyset\}$;
    \item[(5)] $P^{5}_{k}(m,n)=\{\Delta \in P_{k}(m,n)\colon d_{1}=d_{k-1}\ge 2, 1\le \alpha_{1}<d_{1}\}$;
    \item[(6)] $P^{6}_{k}(m,n)=\{\Delta \in P_{k}(m,n)\colon d_{1}=d_{k-1}\ge 2, \alpha_{1}=d_{1}, \alpha_2=0, \beta_1<d_{1}\}$;
    \item[(7)] $P^{7}_{k}(m,n)=\{\Delta \in P_{k}(m,n)\colon d_{1}=d_{k-1}\ge 2, \alpha_{1}=d_{1}, \alpha_2=0, \beta_1=d_{1}\}$; 
    \item[(8)] $P^{8}_{k}(m,n)=\{\Delta \in P_{k}(m,n)\colon d_{1}=d_{k-1}\ge 2,  \alpha_{1} = d_{1} > \alpha_2 \ge 1, \alpha_1-\alpha_2 \text{ is odd}\}$;
    \item[(9)] $P^{9}_{k}(m,n)=\{\Delta \in P_{k}(m,n)\colon d_{1}=d_{k-1}\ge 2, \alpha_{1} = d_{1} > \alpha_2 \ge 1, \alpha_1-\alpha_2 \text{ is even}\}$;
    \item[(10)] $P^{10}_{k}(m,n)=\{\Delta \in P_{k}(m,n)\colon d_{1}=d_{k-1}\ge 2, \alpha_{1} = d_{1} =\alpha_2, 1\le \beta_1<d_1\}$;
    \item[(11)] $P^{11}_{k}(m,n)=\{\Delta \in P_{k}(m,n)\colon d_{1}=d_{k-1}\ge 2, \alpha_{1} = d_{1} =\alpha_2, \beta_1=d_1\}$;
    \item[(12)] $P^{12}_{k}(m,n)=\{\Delta \in P_{k}(m,n)\colon d_{1}=d_{k-1}\ge 2, \alpha_{1} = d_{1} =\alpha_2=\alpha_3, \beta=\emptyset\}$;
    \item[(13)] $P^{13}_{k}(m,n)=\{\Delta \in P_{k}(m,n)\colon d_{1}=d_{k-1}\ge 2, \alpha_{1} = d_{1} =\alpha_2>\alpha_3, \beta=\emptyset,  k\ge 4\}$;
    \item[(14)] $P^{14}_{k}(m,n)=\{\Delta \in P_{k}(m,n)\colon k=3, d_{1}=d_{2}\ge 3, \alpha_{1} = d_{1} =\alpha_2>\alpha_3, \beta=\emptyset\}$;
    \item[(15)] $P^{15}_{k}(m,n)=\{\Delta \in P_{k}(m,n)\colon k=3, d_{1}=d_{2}= 2, \alpha_{1} = d_{1} =\alpha_2>\alpha_3, \beta=\emptyset\}$;
    \item[(16)] $P^{16}_{k}(m,n)=\{\Delta \in P_{k}(m,n)\colon \text{either }d_1=d_{k-1}=1, n=m+k-1 \text{ or } k=3, d_1=d_{2}=2, \alpha=\emptyset\}$
\end{itemize}

We now divide the set $P_{k}(m,n+1)$ into $15$ disjoint subsets $P_{k,i}(m,n+1)$($1\le i \le 15 $) as follows.
\begin{itemize}
    \item[(1)] $P_{k,1}(m,n+1)=\{\Delta \in P_{k}(m,n+1)\colon \gamma^i_1>\gamma^{i}_2, \text{ where } i=\min\{j\colon d_j>d_{j+1}\}\}$;
    \item[(2)] $P_{k,2}(m,n+1)=\{\Delta \in P_{k}(m,n+1)\colon \gamma^{1}=\cdots =\gamma^{k-2}=\emptyset, d_1=2, d_{2}=d_{k-1}=1\}$;
    \item[(3)] $P_{k,3}(m,n+1)=\{\Delta \in P_{k}(m,n+1)\colon k\ge 4, \gamma^{1}=\cdots =\gamma^{k-2}=\emptyset,  d_{1}=d_{k-2}=d_{k-1}+1\ge 2, \alpha=1^{d_1}, \beta=1^{d_1}\}$;
    \item[(4)] $P_{k,4}(m,n+1)=\{\Delta \in P_{k}(m,n+1)\colon  k=3, \gamma^{1}=\emptyset, d_1=d_2+1\ge 3, \alpha=\beta=1^{d_1}\}$;
    \item[(5)] $P_{k,5}(m,n+1)=\{\Delta \in P_{k}(m,n+1)\colon   d_{1}=d_{k-1}\ge 2, \alpha_1 >\alpha_2, \alpha_1\ge 2\}$;
    \item[(6)] $P_{k,6}(m,n+1)=\{\Delta \in P_{k}(m,n+1)\colon   d_{1}=d_{2}+1=\cdots=d_{k-2}+1=d_{k-1}+2\ge 3, \gamma_1=\cdots=\gamma_{k-2}=\emptyset, \alpha_1 =d_1-2, \alpha_2=0\}$;
    \item[(7)] $P_{k,7}(m,n+1)=\{\Delta \in P_{k}(m,n+1)\colon   d_{1}=d_{2}+1=\cdots=d_{k-1}+1\ge 3, \gamma_1=\cdots=\gamma_{k-2}=\emptyset, \alpha=\emptyset\}$;   
    \item[(8)] $P_{k,8}(m,n+1)=\{\Delta \in P_{k}(m,n+1)\colon   d_{1}=d_{k-1}\ge 2, 1\le \alpha_1=\alpha_2< d_{1}, f_{1}(\beta)\ge \frac{d_1-\alpha_1+1}{2}, d_1-\alpha_1 \text{ is odd} \}$;
    \item[(9)] $P_{k,9}(m,n+1)=\{\Delta \in P_{k}(m,n+1)\colon   d_{1}=d_{k-1}\ge 2, 1\le \alpha_1=\alpha_2< d_{1}, d_{1}-\alpha_1 \text{ is even}\}$;
    \item[(10)] $P_{k,10}(m,n+1)=\{\Delta \in P_{k}(m,n+1)\colon   d_{1}=d_{k-1}\ge 2,  \alpha_1=\alpha_2 = d_{1}, \beta_1 >\beta_2, \beta_{1}\ge 2 \}$;
    \item[(11)] $P_{k,11}(m,n+1)=\{\Delta \in P_{k}(m,n+1)\colon   d_{1}=d_{2}+1=\cdots=d_{k-1}+1\ge 3, \gamma_1=\cdots=\gamma_{k-2}=\emptyset, \alpha_1=d_1-1\}$;
    \item[(12)] $P_{k,12}(m,n+1)=\{\Delta \in P_{k}(m,n+1)\colon   d_{1}=d_{k-1}\ge 2,  d_1 = \alpha_{1}=\alpha_{2}, f_{d_1-1}(\alpha)\ge 1, \beta=(1) \}$;
    \item[(13)] $P_{k,13}(m,n+1)=\{\Delta \in P_{k}(m,n+1)\colon k\ge 4, d_1=d_2+1=\cdots=d_{k-2}+1=d_{k-1}+2\ge 3, \gamma^{1}=\cdots=\gamma^{k-3}=\emptyset, \gamma^{k-2}=(1), \alpha_1=\alpha_2=d_{k-1}, \beta=\emptyset\}$;
    \item[(14)] $P_{k,14}(m,n+1)=\{\Delta \in P_{k}(m,n+1)\colon k=3, d_1=d_2+2\ge 4,  \gamma^{1}=(1,1), \alpha_1=d_2, f_{\alpha_1-1}(\alpha)\ge 1\}$;
    \item[(15)] $P_{k,15}(m,n+1)=\{\Delta \in P_{k}(m,n+1)\colon k=3, d_{1}=d_{2}=2, \alpha_1=\alpha_2=\alpha_3=1,\ell(\alpha)\ge 3, \beta=(2)\}$.
\end{itemize}
It is easy to see that the set $P^{i}_{k}(m,n)\in P_{k}(m,n)$ ($1\le i\le 16$) are pairwise disjoint  and $P_{k}(m,n)=\bigcup_{i=1}^{16}P^{i}_{k}(m,n)$. Moreover, the set $P_{k,i}(m,n+1)\in P_{k}(m,n+1)$ ($1\le i\le 15$) are not intersected.

We proceed to show that for $1\le i\le 15$ there exist injections $\sigma^{i}$ from $P^{i}_{k}(m,n)$ to $P_{k,i}(m,n+1)$. In fact, three of these $\sigma^i$ are injections, namely $\sigma^1$, $\sigma^8$ and $\sigma^9$, and all the other $\sigma^i$ are bijections.

We next describe $\sigma^1 \sim \sigma^{15}$ in Lemma \ref{lem-pp1} $\sim$ Lemma \ref{lem-pp15} respectively.

\begin{lem}\label{lem-pp1}
    For $k\ge 3$, $n\ge k-1$ and $m\ge 0$, there exists an injection $\sigma^{1}$ from $P^{1}_{k}(m,n)$ to $P_{k,1}(m,n+1)$.
\end{lem}
\begin{proof}
    Given $\Delta=(\alpha, \beta, \gamma^1, \gamma^2, \ldots, \gamma^{k-2}, \varpi^{1}, \ldots, \varpi^{{k-1}}) \in P^{1}_{k}(m,n)$, by definition, we know $d_{1}\ne d_{k-1}$. Thus let $i$ denote the minimum integer such that $d_{i}>d_{i+1}$, $1\le i\le k-2$. Then we know $\ell(\gamma^i)\le d_{i}-d_{i+1}$. Define
    \begin{align*}
   \tilde{\Delta}:=& \sigma^{1}(\Delta)\\
    =&(\tilde{\alpha}, \tilde{\beta}, \tilde{\gamma}^{1}, \ldots, \tilde{\gamma}^{k-2}, \tilde{\varpi}^{1}, \ldots, \tilde{\varpi}^{k-1})\\
    =&(\alpha, \beta, \gamma^1, \ldots, \gamma^{i-1},(\gamma^i_1+1, \gamma^{i}_{2}, \ldots),\gamma^{i+1},\ldots, \gamma^{k-2}, \varpi^{1}, \ldots, \varpi^{k-1}).
    \end{align*}
    It is obvious that $\tilde{\gamma}^{i}_1=\gamma^i_1+1>\gamma^{i}_2=\tilde{\gamma}_2$ and $|\sigma^{1}(\Delta)|=|\Delta|+1=n+1$. Hence $\sigma^{1}(\Delta)\in P_{k,1}(m,n+1)$. To prove that the map $\sigma$ is an injection, let $
    H_{k,1}(m,n+1)$ be the image set of $\sigma^1$, which has been already shown to be a subset of $P_{k,1}(m,n+1)$. For any 
    $$\tilde{\Delta}=(\tilde{\alpha}, \tilde{\beta}, \tilde{\gamma}^{1}, \ldots, \tilde{\gamma}^{k-2}, \tilde{\varpi}^{1}, \ldots, \tilde{\varpi}^{k-1})\in H_{k,1}(m,n+1),$$ 
    by definition we see that there exists $i$ such that $\tilde{\gamma}^i\ne \emptyset$. We may choose such $i$ to be minimum. 
    Moreover, by the construction of $\sigma^1$, we find that $\tilde{\gamma}_1^i>\tilde{\gamma}^i_2$, and $\tilde{d}_1=\cdots=\tilde{d}_i>\tilde{d}_{i+1}$. Define 
    \begin{align*}
        \tilde{\Delta}:=&\zeta^{1}(\tilde{\Delta})\\
        =&(\alpha, \beta, \gamma^1, \gamma^2, \ldots, \gamma^{k-2}, \varpi^{1}, \ldots, \varpi^{k-1})\\
        =&(\tilde{\alpha}, \tilde{\beta}, \tilde{\gamma}^{1}, \ldots, \tilde{\gamma}^{i-1},(\tilde{\gamma}^{i}_1-1, \tilde{\gamma}^{i}_{2}, \ldots), \tilde{\gamma}^{i+1},\ldots, \tilde{\gamma}^{k-2},\tilde{\varpi}^{1}, \ldots, \tilde{\varpi}^{k-1}).
    \end{align*}
    It can be verified that $d_i=\tilde{d}_i>\tilde{d}_{i+1}=d_{i+1}$ and $|\zeta^{1}(\tilde{\Delta})|=n+1-1=n$, thus $\zeta^{1}(\tilde{\Delta}) \in P^{1}_{k}(m,n)$ and $\zeta^{1}(\sigma^{1}(\Delta))=\Delta$ for any $\Delta\in P^{1}_{k}(m,n)$. Hence the map $\sigma^{1}$ is an injection from $P^{1}_{k}(m,n)$ to $P_{k,1}(m,n+1)$.
\end{proof}
For example, let $\Delta=((2,1,1), (1, 1), (2), \emptyset, (3^3), (2^2), (2^2))\in P^{1}_{4}(1,25)$, using the $\sigma^{1}$ on $\Delta$, we get $\tilde{\Delta}=((2,1,1), (1,1), (3), \emptyset, (3^3), (2^2), (2^2))\in P_{4,1}(1,26)$. Applying $\zeta^1$ on  $\tilde{\Delta}$, we recover $\Delta$.

\begin{lem}\label{lem-pp2}
For $k\ge 3$, $n\ge m+k$ and $m\ge 0$, there exists a bijection $\sigma^{2}$ between $P^{2}_{k}(m,n)$ and $P_{k,2}(m,n+1)$.
\end{lem}
\begin{proof}
Let $\Delta=(\alpha, \beta, \gamma^1, \gamma^2, \ldots, \gamma^{k-2}, \varpi^{1}, \ldots, \varpi^{k-1}) \in P^{2}_{k}(m,n)$. By definition, we know $d_{1}=d_{k-1}=1$. Thus by definition we have $\gamma^i=\emptyset$ for $1\le i\le k-2$ and $\alpha_1\le 1$, $\beta_1\le 1$. Since $m=\ell(\alpha)-\ell(\beta)$, Assume $\beta=1^{t}$, then from $m=\ell(\alpha)-\ell(\beta)$, we see that $\alpha=1^{m+t}$. Since $n=|\alpha|+|\beta|+1\times (k-1)=2t+m+k-1\ge m+k$, we deduce that $t\ge 1$.

Define
\begin{align*}
    \tilde{\Delta}:=&\sigma^{2}(\Delta)\\
    =&(\tilde{\alpha}, \tilde{\beta}, \tilde{\gamma}^{1}, \ldots, \tilde{\gamma}^{k-2}, \tilde{\varpi}^{1}, \ldots, \tilde{\varpi}^{k-1})\\
    =&((\underbrace{1,\ldots,1}_{m+t-1}), (\underbrace{1,\ldots,1}_{t-1}), \underbrace{\emptyset, \ldots, \emptyset}_{k-2}, (2,2),\underbrace{(1), \ldots, (1)}_{k-2}).
\end{align*}

It is easy to check that $\tilde{d}_1=\ell(\tilde{\varpi}^{1})=2$ and $\tilde{d_2}=\ell(\tilde{\varpi^{2}})=\cdots=\ell(\tilde{\varpi}^{k-1})=\tilde{d_{k-1}}=1$. Moreover, $|\sigma^{2}(\Delta)|=n-2+4-1=n+1$. Thus $\sigma^{2}(\Delta)\in P_{k,2}(m,n+1)$. 

To show that $\sigma^{2}$ is a bijection, we construct the inverse map $\zeta^{2}$. Let 
\[
\tilde{\Delta}=(\tilde{\alpha}, \tilde{\beta}, \tilde{\gamma}^{1}, \ldots, \tilde{\gamma}^{k-2}, \tilde{\varpi}^{1}, \ldots, \tilde{\varpi}^{k-1})\in P_{k,2}(m,n).
\]
By definition, we have   $\tilde{d}_1=\ell(\tilde{\varpi}^{1})=2$, $\tilde{d}_{i}=\ell(\tilde{\varpi}^{i})=1$ and $\gamma^{i-1}=\emptyset$ for $2\le i\le k-1$. By Definition \ref{def-Pkmn}, we find that $\alpha_1\le 1$ and $\beta_1\le 1$. Assume $\ell(\beta)=s$, where $s\ge 0$. Then $\beta=1^s$ and $\alpha=1^{m+s}$. Define $\zeta^{2}(\tilde{\Delta})$ to be
\begin{align*}
    \zeta^{2}(\tilde{\Delta})=&(\alpha, \beta, \gamma^1, \gamma^2, \ldots, \gamma^{k-2}, \varpi^{1}, \ldots, \varpi^{k-1})\\
    =&((\underbrace{1,\ldots,1}_{m+s+1}), (\underbrace{1,\ldots,1}_{s+1}), \underbrace{\emptyset, \ldots, \emptyset}_{k-2}, \underbrace{(1),\ldots,(1)}_{k-1}).
\end{align*}

It is easy to check that $\zeta^{2}(\tilde{\Delta})\in P^{2}_{k}(m,n)$ and $\zeta^2$ is the inverse map of $\sigma^2$. So we conclude that $\sigma^{2}$ is a bijection between $P^{2}_{k}(m,n)$ and $P_{k,2}(m,n+1)$.
\end{proof}	
For example, let $\Delta=((1,1,1), (1, 1), \emptyset, \emptyset, (1^1), (1^1), (1^1))\in P^{2}_{4}(1,8)$, using the $\sigma^{2}$ on $\Delta$, we get $\tilde{\Delta}=((1,1), (1), \emptyset, \emptyset, (2^2), (1^1), (1^1))\in P_{4,2}(1,9)$. Applying $\zeta^2$ on  $\tilde{\Delta}$, we recover $\Delta$.

\begin{lem}\label{lem-pp3}
For $k\ge 4$, $n\ge k-1$ and $m\ge 0$, there exists a bijection $\sigma^{3}$ between $P^{3}_{k}(m,n)$ and $P_{k,3}(m,n+1)$.
\end{lem}
\begin{proof}
Given $\Delta=(\alpha, \beta,  \gamma^1, \ldots, \gamma^{k-2}, \varpi^{1}, \ldots, \varpi^{k-1})\in P^{3}_{k}(m,n)$, by definition, we know that $d_{1}=d_{k-1}\ge 2$, $k\ge 4$, $\alpha=\emptyset$. Thus by Definition \ref{def-Pkmn}, $\gamma^{1}=\cdots=\gamma^{k-2}=\emptyset$. Moreover, since $m\ge 0$, we see that $0\le \ell(\beta)=\ell(\alpha)-m=-m\le 0$, which yields $m=0$ and $\alpha=\beta=\emptyset$. 
Define
\begin{align*}
    \tilde{\Delta}:=&\sigma^{3}(\Delta)\\
    =&(\tilde{\alpha}, \tilde{\beta}, \tilde{\gamma}^{1}, \ldots, \tilde{\gamma}^{k-2}, \tilde{\varpi}^{1}, \ldots, \tilde{\varpi}^{k-1})\\
    =&(1^{d_1}, 1^{d_1}, \underbrace{\emptyset, \ldots, \emptyset}_{k-2}, \varpi^{1}, \ldots, \varpi^{{k-2}}, ((d_{k-1}-1)^{d_{k-1}-1})).
\end{align*}
It is easy to see that $\tilde{d}_{i}=\ell(\tilde{\varpi}^{i})\ge 2$($1\le i \le k-2$), $\tilde{d}_{k-1}=\ell(\tilde{\varpi}^{k-1})=d_{k-1}-1=\ell(\tilde{\varpi}^{1})-1=\tilde{d}_{1}-1$ and $\tilde{\alpha}=\tilde{\beta}=(1^{\ell(\tilde{\varpi}^{1})})=(1^{\tilde{d}_{1}})$.    Moreover, $|\sigma^{3}(\Delta)|=n-(2d_{k-1}-1)+2d_{k-1}=n+1$. Thus we deduce that $\sigma^{3}(\Delta)\in P_{k,3}(m,n+1)$. To prove $\sigma^{3}$ is a bijection, we construct the inverse map $\zeta^{3}$ of $\sigma^{3}$. Let
\[
    \tilde{\Delta}=(\tilde{\alpha}, \tilde{\beta}, \tilde{\gamma}^{1}, \ldots, \tilde{\gamma}^{k-2}, \tilde{\varpi}^{1}, \ldots, \tilde{\varpi}^{k-1})
\]
be a $(2k-1)$-tuple partition in $P_{k,3}(m,n+1)$. By definition, we know $\tilde{d}_{1}=\cdots=\tilde{d}_{k-2}=\tilde{d}_{k-1}+1\ge 2$, and $\tilde{\alpha}=\tilde{\beta}=(1^{\tilde{d}_{1}})$. Define $\zeta^{3}(\tilde{\Delta})$ to be
\begin{align*}
    \zeta^{3}(\tilde{\Delta}):=&(\alpha, \beta,  \gamma^1, \gamma^2, \ldots, \gamma^{k-2}, \varpi^{1}, \ldots, \varpi^{k-1})\\
    =&(\underbrace{\emptyset, \ldots, \emptyset}_{k}, \tilde{\varpi}^{1}, \ldots, \tilde{\varpi}^{k-2}, (\tilde{d}_{k-1}+1)^{\tilde{d}_{k-1}+1}).
\end{align*}
Note that $d_{1}=\ell(\varpi^{1})=\ldots=\ell(\varpi^{k-2})=d_{k-2}=\tilde{d}_{1}\ge 2$, $d_{k-1}=\ell(\varpi^{k-1})=\tilde{d}_{k-1}+1=\tilde{d}_{1}=d_{1}$ and $\alpha=\beta=\emptyset$. Moreover, 
\[
|\zeta^{3}(\tilde{\Delta})|=n+1-2\tilde{d}_1+(2\tilde{d}_{k-1}+1)=n+1-2\tilde{d}_1+(2(\tilde{d}_{1}-1)+1)=n,
\]
thus we deduce that $\zeta^{3}(\tilde{\Delta})\in P^{3}_{k}(m,n)$ and it is easy to check that $\zeta^{3}$ is the inverse map of $\sigma^{3}$. Thus we conclude that $\sigma^{3}$ is a bijection between $P^{3}_{k}(m,n)$ and $P_{k,3}(m,n+1)$.
\end{proof}

For example, let $\Delta=(\emptyset, \emptyset, \emptyset, \emptyset, (2^2), (2^2), (2^2))\in P^{3}_{4}(0,12)$, using the $\sigma^{3}$ on $\Delta$, we get $\tilde{\Delta}=((1,1), (1,1), \emptyset, \emptyset, (2^2), (2^2), (1^1))\in P_{4,3}(0,13)$. Applying $\zeta^3$ on  $\tilde{\Delta}$, we recover $\Delta$.

\begin{lem}\label{lem-pp4}
For $k= 3$, $n\ge 2$, $m\ge 0$,  there exists a bijection $\sigma^{4}$ between $P^{4}_{3}(m,n)$ and $P_{3,4}(m,n+1)$.
\end{lem}
\begin{proof}
Given $\Delta=(\alpha, \beta,  \gamma^1,  \varpi^{1},\varpi^{2})\in P^{4}_{3}(m,n)$, by definition, we see that $\alpha=\emptyset$ and $d_1=d_2\ge 3$. Since $\ell(\alpha)-\ell(\beta)=-\ell(\beta)=m\ge 0$, we deduce that  $m=0$ and $\beta=\emptyset$. Thus 
\[
\Delta=(\emptyset, \emptyset, \emptyset, \varpi^{1}, \varpi^{2}).
\]
Define
\begin{align*}
\tilde{\Delta}:=&\sigma^{4}(\Delta)\\
=&(\tilde{\alpha}, \tilde{\beta}, \tilde{\gamma}^{1}, \tilde{\varpi}^{1},  \tilde{\varpi}^{2})\\
=&(1^{d_1}, 1^{d_1}, \emptyset, \varpi^{1}, (d_2-1)^{d_2-1}).  
\end{align*}

It is easy to check that 
\[
\tilde{d}_1= d_1=d_2= \tilde{d_2}+1\ge 3,
\]
$\tilde{\gamma}^{1}=\emptyset$ and $\tilde{\alpha}=\tilde{\beta}=1^{d_1}=1^{\tilde{d}_1}$. Moreover,  $|\sigma^{4}(\Delta)|=n-(2d_2-1)+2d_1=n+1$. Thus we confirm that $\sigma^{4}(\Delta) \in P_{3,4}(m,n+1)$. To show $\sigma^{4}$ is a bijection, we construct the inverse map $\zeta^{4}$ of $\sigma^{4}$. Let
\[
    \tilde{\Delta}=(\tilde{\alpha}, \tilde{\beta}, \tilde{\gamma}^{1}, \tilde{\varpi}^{1},   \tilde{\varpi}^{2})
\]
denote the $5$-tuple partition in $P_{3,4}(m,n+1)$ where $\tilde{d}_1=\ell(\varpi^{1})=\tilde{d}_2+1=\ell(\varpi^{2})+1 \ge 3$ and $\tilde{\alpha}=\tilde{\beta}=1^{\tilde{d_1}}$. Define
\begin{align*}
\Delta:=&\zeta^{4}(\tilde{\Delta})\\
=&(\alpha, \beta,  \gamma^1,  \varpi^{1}, \varpi^{2})\\
 =&(\emptyset, \emptyset, \emptyset, \tilde{\varpi}^{1}, (\tilde{d}_2+1)^{\tilde{d}_2+1}).   
\end{align*}
Note that
\[
d_1= \tilde{d}_1=\tilde{d}_2+1= d_2\ge 3,
\]
$\alpha=\beta=\emptyset$ and $|\zeta^{4}(\tilde{\Delta})|=n+1-(2\tilde{d}_1)+(2(\tilde{d}_2+1)-1)=n+1-(2\tilde{d}_1)+(2\tilde{d}_1-1)=n$. Thus we confirmed that $\zeta^{4}(\tilde{\Delta})\in P^{4}_{3}(m,n)$. It is easy to check that $\zeta^{4}$ is the inverse map of $\sigma^{4}$. Thus $\sigma^{4}$ is a bijection.
\end{proof}

For example, let $\Delta=(\emptyset, \emptyset, \emptyset, (3^3), (3^3))\in P^{4}_{3}(0,18)$, using the $\sigma^{4}$ on $\Delta$, we get $\tilde{\Delta}=((1^3), (1^3), \emptyset, (3^3), (2^2))\in P_{3,4}(0,19)$. Applying $\zeta^4$ on $\tilde{\Delta}$, we recover $\Delta$.

\begin{lem}\label{lem-pp5}
For $k\ge 3$, $n\ge k-1$, $m\ge 0$,  there exists a bijection $\sigma^{5}$ between $P^{5}_{k}(m,n)$ and $P_{k,5}(m,n+1)$.
\end{lem}
\begin{proof}
Given $\Delta=(\alpha, \beta, \gamma^1, \ldots, \gamma^{k-2}, \varpi^{1}, \ldots, \varpi^{k-1})\in P^{5}_{k}(m,n)$, by definition we know $d_1=d_{k-1}\ge 2$ and  $\gamma^{i}=\emptyset$ for $1\le i\le k-2$. Moreover $1\le \alpha_1 < d_{k-1}$. Define
\begin{align*}
    \tilde{\Delta}:=&\sigma^{5}(\Delta)\\
    =&(\tilde{\alpha}, \tilde{\beta}, \tilde{\gamma}^{1}, \ldots, \tilde{\gamma}^{k-2}, \tilde{\varpi}^{1}, \ldots, \tilde{\varpi}^{k-1})\\
    =&((\alpha_{1}+1, \alpha_2, \ldots), \beta, \underbrace{\emptyset, \ldots, \emptyset}_{k-2}, \varpi^{1}, \ldots, \varpi^{k-1}).
\end{align*}
It is easy to see that $\tilde{d}_{1}=d_1=d_{k-1}=\tilde{d}_{k-1}\ge 2$($1\le i\le k-1$),  $\tilde{\alpha}_{1}=\alpha_1+1>\alpha_2=\tilde{\alpha}_2$, and $\tilde{\alpha}_1=\alpha_1+1\ge 2$. Moreover, $|\sigma^{5}(\Delta)|=n+1$. Hence $\sigma^{5}(\Delta)\in P_{k,5}(m,n+1)$. We next show that $\sigma^{5}$ is a bijection between $P^{5}_{k}(m,n)$ and $P_{k,5}(m,n+1)$. 
For any
\[
\tilde{\Delta}=(\tilde{\alpha}, \tilde{\beta}, \tilde{\gamma}^{1}, \ldots, \tilde{\gamma}^{k-2}, \tilde{\varpi}^{1}, \ldots, \tilde{\varpi}^{k-1})\in P_{k,5}(m,n+1),
\]
by the definition of $P_{k,5}(m,n+1)$, we know $\tilde{\alpha}_1>\tilde{\alpha}_2$ , $\tilde{d}_{1}=\tilde{d}_{k-1}\ge 2$, $\tilde{\alpha}_1\ge 2$. Define
\begin{align*}
    \Delta:=&\zeta^{5}(\tilde{\Delta})\\
    =&(\alpha,\beta, \gamma^{1}, \ldots, \gamma^{k-2}, \varpi^{1}, \ldots, \varpi^{k-1})\\
    =&((\tilde{\alpha}_{1}-1, \tilde{\alpha}_2, \ldots), \tilde{\beta}, \underbrace{\emptyset, \ldots, \emptyset}_{k-2}, \tilde{\varpi}^{1}, \ldots, \tilde{\varpi}^{k-1}).
\end{align*}
It is easy to check that $\alpha_1=\tilde{\alpha}_{1}-1\ge 1$, $\alpha_1=\tilde{\alpha}_{1}-1\ge \tilde{\alpha}_2=\alpha_2$, $\alpha_1=\tilde{\alpha}_1-1< \tilde{d}_1=d_1$ and $d_1=\tilde{d}_1=\tilde{d}_{k-1}=d_{k-1}\ge 2$. Moreover, $|\zeta^{5}(\tilde{\Delta})|=n+1-1=n$. Thus we deduce that $\zeta^{5}(\tilde{\Delta})\in P^{5}_{k}(m,n)$ and it can be easily checked that $\zeta^{5}(\sigma^{5}(\Delta))=\Delta$. Thus $\sigma^{5}$ is a bijection between $P^{5}_{k}(m,n)$ and $P_{k,5}(m,n+1)$. .
\end{proof}

For example, let $\Delta=((2,1,1), (1), \emptyset, (3^3), (3^3))\in P^{5}_{3}(2,23)$, using the $\sigma^{5}$ on $\Delta$, we get $\tilde{\Delta}=((3,1,1), (1), \emptyset, (3^3), (3^3))\in P_{3,5}(2,24)$. Applying $\zeta^5$ on  $\tilde{\Delta}$, we recover $\Delta$.

\begin{lem}\label{lem-pp6}
  For $k\ge 3$, $n\ge k-1$, $m\ge 0$,  there exists a bijection $\sigma^{6}$ between $P^{6}_{k}(m,n)$ and $P_{k,6}(m,n+1)$.  
\end{lem}
\begin{proof}
 Given $\Delta=(\alpha, \beta, \gamma^{1}, \ldots, \gamma^{k-2}, \varpi^{1}, \ldots, \varpi^{k-1})\in P^{6}_{k}(m,n)$, by definition, we know that $d_1=d_{k-1}\ge 2$, $\alpha_1=d_1$, $\alpha_2=0$ and $\beta_1<d_{1}$. Define
 \begin{align*}
\tilde{\Delta}:=&\sigma^{6}(\Delta)\\
=&(\tilde{\alpha}, \tilde{\beta}, \tilde{\gamma}^{1}, \ldots, \gamma^{k-2}, \tilde{\varpi}^{1}, \ldots, \tilde{\varpi}^{k-1})\\
=&((\alpha_1-1), \beta, \underbrace{\emptyset, \ldots, \emptyset}_{k-2}, (d_1+1)^{d_1+1}, \varpi^{2}, \ldots, \varpi^{k-2}, (d_{k-1}-1)^{d_{k-1}-1}).
 \end{align*}
 Note that for $2\le i\le k-2$, 
 \[
 \tilde{d}_1=d_1+1=d_{i}+1=\tilde{d}_i+1=d_{k-1}-1+2=\tilde{d}_{k-1}+2\ge 3,
 \]
 and $\tilde{\alpha}_1=\alpha_1-1=d_1-1=\tilde{d}_1-2\ge 1$. Moreover, $\tilde{d}_{k-1}=d_{k-1}-1\ge \beta_1=\tilde{\beta}_1$ and $|\sigma^{6}(\Delta)|=n-1-(2d_{k-1}-1)+(2(d_{1}+1)-1)=n+1$. Thus we deduce that $\sigma^{6}(\Delta)\in P_{k,6}(m,n)$. To show that $\sigma^{6}$ is a bijection between $P^{6}_{k}(m,n)$ and $P_{k,6}(m,n)$, we now consider the inverse map of $\sigma^{6}$. Let 
 \[
\tilde{\Delta}=(\tilde{\alpha}, \tilde{\beta}, \tilde{\gamma}^{1}, \ldots, \tilde{\gamma}^{k-2}, \tilde{\varpi}^{1}, \ldots, \tilde{\varpi}^{k-1})\in P_{k,6}(m,n+1).
\]
By definition, we know that for $2\le i\le k-2$, $\tilde{d}_1=\tilde{d}_i+1=\tilde{d}_{k-1}+2\ge 3$, $\tilde{\gamma}^{i}=\emptyset$,  $\tilde{\alpha}_1=\tilde{d}_1-2 \ge 1$ and $\tilde{\alpha}_2=0$. Define
\begin{align*}
    \Delta:=\zeta^{6}(\tilde{\Delta})=&(\alpha, \beta, \gamma^{1}, \ldots, \gamma^{k-2}, \varpi^{1}, \ldots, \varpi^{k-1})\\
    =&((\tilde{\alpha}_1+1), \tilde{\beta}, \underbrace{\emptyset, \ldots, \emptyset}_{k-2}, (\tilde{d}_1-1)^{\tilde{d}_1-1}, \tilde{d}_2, \ldots, \tilde{d}_{k-2}, (\tilde{d}_{k-1}+1)^{\tilde{d}_{k-1}+1}).
\end{align*}
It is easy to see that for $2\le i\le k-2$, $d_1=\tilde{d}_1-1=\tilde{d}_i=d_i=\tilde{d}_{k-1}+1=d_{k-1}$,  $\alpha_1=\tilde{\alpha}_1+1=\tilde{d}_1-2+1=d_1$, $\alpha_2=0$ and $\beta_1=\tilde{\beta_1}\le \tilde{d}_{k-1}<\tilde{d}_{k-1}+1=d_{k-1}$. Moreover,  $\gamma^{i}=\emptyset$ for $1\le i\le k-2$ and $|\zeta^{6}(\tilde{\Delta})|=n+1+1-(2\tilde{d}_1-1)+(2(\tilde{d}_{k-1}+1)-1)=n$. Thus We confirm that $\zeta^{6}(\tilde{\Delta})\in P^{6}_{k}(m,n)$. Furthermore, it is easy to verify $\zeta^{6}$ is the inverse map of $\sigma^6$. Thus we deduce that $\sigma^{6}$ is a bijection between $P^{6}_{k}(m,n)$ and $P_{k,6}(m,n+1)$.
\end{proof}

For example, let $\Delta=((2), (1), \emptyset, \emptyset, (2^2), (2^2), (2^2))\in P^{6}_{4}(0,15)$, using the $\sigma^{6}$ on $\Delta$, we get $\tilde{\Delta}=((1), (1), \emptyset, \emptyset, (3^3), (2^2), (1^1))\in P_{4,6}(0,16)$. Applying $\zeta^6$ on $\tilde{\Delta}$, we recover $\Delta$.

\begin{lem}\label{lem-pp7}
  For $k\ge 3$, $n\ge k-1$, $m\ge 0$,  there exists a bijection $\sigma^{7}$ between $P^{7}_{k}(m,n)$ and $P_{k,7}(m,n+1)$.  
\end{lem}
\begin{proof}
 Given $\Delta=(\alpha, \beta, \gamma^{1}, \ldots, \gamma^{k-2}, \varpi^{1}, \ldots, \varpi^{k-1})\in P^{7}_{k}(m,n)$, by definition, we know that $d_1=d_{k-1}\ge 2$, $\alpha_1=d_1$, $\alpha_2=0$ and $\beta_1=d_{1}$. Moreover, $m\ge 0$ implies that $\ell(\beta)=1$. Define
 \begin{align*}
\tilde{\Delta}:=&\sigma^{7}(\Delta)\\
=&(\tilde{\alpha}, \tilde{\beta}, \tilde{\gamma}^{1}, \ldots, \gamma^{k-2}, \tilde{\varpi}^{1}, \ldots, \tilde{\varpi}^{k-1})\\
=&(\emptyset, \emptyset, \underbrace{\emptyset, \ldots, \emptyset}_{k-2}, (d_1+1)^{d_1+1}, \varpi^{2}, \ldots, \varpi^{k-1}).
 \end{align*}
 Note that for $2\le i\le k-1$, $\tilde{d}_1=d_1+1=d_{i}+1=\tilde{d}_i+1\ge 3$, and for $1\le i\le k-2$, $\tilde{\gamma}^{i}=\emptyset$. Moreover, $\tilde{\alpha}=\tilde{\beta}=\emptyset $ and $|\sigma^{7}(\Delta)|=n-2d_1+(2d_1+1)=n+1$. Thus we deduce that $\sigma^{7}(\Delta)\in P_{k,7}(m,n)$. To show that $\sigma^{7}$ is a bijection between $P^{7}_{k}(m,n)$ and $P_{k,7}(m,n)$, we now consider the inverse map of $\sigma^{7}$. Let 
 \[
\tilde{\Delta}=(\tilde{\alpha}, \tilde{\beta}, \tilde{\gamma}^{1}, \ldots, \tilde{\gamma}^{k-2}, \tilde{\varpi}^{1}, \ldots, \tilde{\varpi}^{k-1})\in P_{k,7}(m,n+1).
\]
By definition, we know that for $2\le i\le k-1$, $\tilde{d}_1=\tilde{d}_i+1\ge 3$, $\tilde{\gamma}^{i-1}=\emptyset$ and $\tilde{\alpha}=\emptyset$. Moreover, $m\ge0$ implies $\tilde{\beta}=\emptyset$. Define
\begin{align*}
    \Delta:=\zeta^{7}(\tilde{\Delta})=&(\alpha, \beta, \gamma^{1}, \ldots, \gamma^{k-2}, \varpi^{1}, \ldots, \varpi^{k-1})\\
    =&((\tilde{d}_1-1), (\tilde{d}_1-1), \underbrace{\emptyset, \ldots, \emptyset}_{k-2}, (\tilde{d}_1-1)^{\tilde{d}_1-1}, \tilde{\varpi}^{2}, \ldots, \tilde{\varpi}^{k-1}).
\end{align*}
It is easy to see that for $2\le i\le k-1$, $d_1=\tilde{d}_1-1=\tilde{d}_i=d_i$, $\gamma^{i-1}=\emptyset$, $\alpha_1=\tilde{d}_1-1=d_1$, $\alpha_2=0$ and $\beta_1=\tilde{d}_1-1=d_1$. Moreover, $|\zeta^{7}(\tilde{\Delta})|=n+1+(2\tilde{d}_1-2)-(2\tilde{d}_1-1)=n$. We confirm that $\zeta^{7}(\tilde{\Delta})\in P^{7}_{k}(m,n)$. Furthermore, it is easy to verify $\zeta^{7}$ is the inverse map of $\sigma^7$. Thus we deduce that $\sigma^{7}$ is a bijection between $P^{7}_{k}(m,n)$ and $P_{k,7}(m,n+1)$.
\end{proof}

For example, let $\Delta=((2), (2), \emptyset, \emptyset, (2^2), (2^2), (2^2))\in P^{7}_{4}(0,16)$, using the $\sigma^{7}$ on $\Delta$, we get $\tilde{\Delta}=(\emptyset, \emptyset, \emptyset, \emptyset, (3^3), (2^2), (2^2))\in P_{4,7}(0,17)$. Applying $\zeta^7$ on $\tilde{\Delta}$, we recover $\Delta$.

\begin{lem}\label{lem-pp8}
For $k\ge 3$, $n\ge k-1$, $m\ge 0$,  there exists an injection $\sigma^{8}$ from $P^{8}_{k}(m,n)$ to $P_{k,8}(m,n+1)$.
\end{lem}
\begin{proof}
Given $\Delta=(\alpha, \beta, \gamma^{1}, \ldots, \gamma^{k-2}, \varpi^{1}, \ldots, \varpi^{k-1})\in P^{8}_{k}(m,n)$, by definition, we know that $d_{1}=d_{k-1}\ge 2$, $\alpha_1=d_{1}>\alpha_2\ge 1$, $\gamma^{i}=\emptyset$ ($1\le i\le k-2$) and $\alpha_1-\alpha_{2}$ is odd. Define
\begin{align}
\tilde{\Delta}:=&\sigma^{8}(\Delta)\nonumber\\
=&(\tilde{\alpha}, \tilde{\beta}, \tilde{\gamma}^{1}, \ldots, \gamma^{k-2}, \tilde{\varpi}^{1}, \ldots, \tilde{\varpi}^{k-1})\nonumber\\
=&((\alpha_{2}, \alpha_{2}, \alpha_{3}, \ldots, \alpha_{\ell(\alpha)}, 1^{\frac{\alpha_1-\alpha_{2}+1}{2}}), \nonumber\\
&(\beta_{1}, \ldots, \beta_{\ell(\beta)}, 1^{\frac{\alpha_1-\alpha_{2}+1}{2}}),
\underbrace{\emptyset, \ldots, \emptyset}_{k-2}, \varpi^{1}, \ldots, \varpi^{k-1}).\label{eq-sigma8}
\end{align}
Note that $\tilde{d}_1=d_1=d_{k-1}=\tilde{d}_{k-1}\ge 2$ and
\[
1\le \tilde{\alpha}_{1}=\tilde{\alpha}_{2}=\alpha_2<d_1 =\tilde{d}_1.
\]
Moreover, $\tilde{d}_{1}-\tilde{\alpha}_{1}=d_{1}-\alpha_2=\alpha_1-\alpha_2$ is odd, 
\[
\ell(\tilde{\alpha})-\ell(\tilde{\beta})=\ell(\alpha)+\frac{\alpha_1-\alpha_{2}+1}{2}-\left(\ell(\beta)+\frac{\alpha_1-\alpha_{2}+1}{2}\right)=\ell(\alpha)-\ell(\beta)=m,
\]
$f_{1}(\tilde{\beta})\ge \frac{\tilde{d}_1-\tilde{\alpha}_1+1}{2}\ge 1$
and $|\sigma^{8}(\Delta)|=n-(\alpha_1-\alpha_{2})+2\cdot\frac{\alpha_1-\alpha_{2}+1}{2}=n+1$.
Thus $\sigma^{8}(\Delta)\in P_{k,8}(m,n+1)$. To show that $\sigma^{8}$ is an injection from $P^{8}_{k}(m,n)$ to $P_{k,8}(m,n+1)$, let $
H_{k,8}(m,n+1)=\{\sigma^{8}(\Delta)\colon\Delta\in P^{8}_{k}(m,n)\}
$
be the image set of $\sigma^8$, which has already been shown to be a subset of $P_{k,8}(m,n+1)$. Now we construct the inverse map $\zeta^{8}$ from $H_{k,8}(m,n+1)$ to $P^{8}_{k}(m,n)$. Given 
\[
\tilde{\Delta}=(\tilde{\alpha}, \tilde{\beta}, \tilde{\gamma}^{1}, \ldots, \tilde{\gamma}^{k-2}, \tilde{\varpi}^{1}, \ldots, \tilde{\varpi}^{k-1})\in H_{k,8}(m,n+1).
\]
By definition, we know that $\tilde{d}_1=\tilde{d}_{k-1}\ge 2$, $\tilde{\alpha}_1=\tilde{\alpha}_2<\tilde{d}_1$ and $\tilde{d}_1-\tilde{\alpha}_1$ is odd. Moreover, by the construction of $\sigma^{8}$ in \eqref{eq-sigma8}, we know $f_1(\tilde{\alpha})\ge \frac{\tilde{d}_1-\tilde{\alpha}_1+1}{2}$ and $f_1(\tilde{\beta})\ge \frac{\tilde{d}_1-\tilde{\alpha}_1+1}{2}\ge 1$. 

Define
\begin{align*}
    \Delta:=&\zeta^{8}(\tilde{\Delta})\\
    =&(\alpha, \beta, \gamma^{1}, \ldots, \gamma^{k-2}, \varpi^{1}, \ldots, \varpi^{k-1})\\
    =&((\tilde{d}_1, \tilde{\alpha}_{2}, \ldots, \tilde{\alpha}_{\ell(\tilde{\alpha})-\frac{\tilde{d}_1-\tilde{\alpha}_1+1}{2}}),\\
    &(\tilde{\beta}_1, \tilde{\beta}_{2}, \ldots, \tilde{\beta}_{\ell(\tilde{\alpha})-\frac{\tilde{d}_1-\tilde{\alpha}_1+1}{2}}), \underbrace{\emptyset, \ldots, \emptyset}_{k-2}, \tilde{\varpi}^{1}, \ldots, \tilde{\varpi}^{k-1}).
\end{align*}
Note that $d_{1}=\tilde{d}_1=\tilde{d}_{k-1}=d_{k-1}\ge 2$. Moreover, 
\[
\alpha_{1}= \tilde{d}_{1}>\tilde{\alpha}_2=\alpha_2\ge 1.
\]
Furthermore, $\alpha_1-\alpha_2=\tilde{d}_{1}-\tilde{\alpha}_2=\tilde{d}_1-\tilde{\alpha}_1$ is odd, 
\[
\ell(\alpha)-\ell(\beta)=\ell(\tilde{\alpha})-\frac{\tilde{d}_1-\tilde{\alpha}_1+1}{2}-\left(\ell(\tilde{\beta})-\frac{\tilde{d}_1-\tilde{\alpha}_1+1}{2}\right)=m
\]
and $|\zeta^{8}(\tilde{\Delta})|=n+1+\tilde{d}_1-\tilde{\alpha}_1-2\frac{\tilde{d}_1-\tilde{\alpha}_1+1}{2}=n$. Thus we deduce that $\zeta^{8}(\tilde{\Delta})\in P^{8}_{k}(m,n)$ and it is clear that $\zeta^{8}(\sigma^8(\Delta))=\Delta$ for any $\Delta\in P^8_k(m,n)$. Hence the map $\sigma^{8}$ is an injection from $P^{8}_{k}(m,n)$ to $P_{k,8}(m,n+1)$.  
\end{proof}

For example, let $k=4$, $\Delta=((2,1), (1), \emptyset, \emptyset, (2^2), (2^2), (2^2))\in P^{8}_{4}(1,16)$. Using the map $\sigma^{8}$ on $\Delta$ we get $\tilde{\Delta}=((1,1,1), (1,1), \emptyset, \emptyset, (2^2), (2^2), (2^2))\in P_{4,8}(1,17)$. 

Inversely, given $\tilde{\Delta}=((1,1,1), (1,1), \emptyset, \emptyset, (2^2), (2^2), (2^2))\in P_{4,8}(1,17)$, using the inverse map $\zeta^{8}$, we get $\Delta=((2,1), (1), \emptyset, \emptyset, (2^2), (2^2), (2^2))\in P^{8}_{4}(1,16)$.

\begin{lem}\label{lem-pp9}
For $k\ge 3$, $n\ge k-1$, $m\ge 0$,  there exists an injection $\sigma^{9}$ from $P^{9}_{k}(m,n)$ to $P_{k,9}(m,n+1)$.
\end{lem}
\begin{proof}
Given $\Delta=(\alpha, \beta, \gamma^{1}, \ldots, \gamma^{k-2}, \varpi^{1}, \ldots, \varpi^{k-1})\in P^{9}_{k}(m,n)$. By definition, we know that $d_{1}=d_{k-1}\ge 2$,  $\alpha_{1}=d_{1}>\alpha_{2}\ge 1$ and $\alpha_1-\alpha_{2}$ is even. Let $i$ be the maximum integer such that $\beta_i\ge 2$. Now we define $\sigma^{9}$ as follows.
\begin{align}
    \sigma^{9}(\Delta):=&(\tilde{\alpha}, \tilde{\beta}, \tilde{\gamma}^{1}, \ldots, \tilde{\gamma}^{k-2}, \tilde{\varpi}^{1}, \ldots, \tilde{\varpi}^{k-1})\nonumber\\
    =&((\alpha_{2}, \alpha_{2}, \alpha_{3}, \ldots, \alpha_{\ell(\alpha)}, 1^{\frac{\alpha_{1}-\alpha_2}{2}}), \nonumber\\
    &(\beta_{1}, \ldots, \beta_{i},2,\beta_{i+1},\ldots,\beta_{\ell(\beta)}, 1^{\frac{\alpha_{1}-\alpha_2}{2}-1}), \underbrace{\emptyset, \ldots, \emptyset}_{k-2}, \varpi^{1}, \ldots, \varpi^{k-1}).\label{eq-sigma9}
\end{align}
Note that $\tilde{d}_{1}=d_1=d_{k-1}=\tilde{d}_{k-1}\ge 2$.
Moreover, $\tilde{d}_1-\tilde{\alpha}_1=d_{1}-\alpha_{2}=\alpha_1-\alpha_2$ is even, $1\le \tilde{\alpha}_1=\alpha_2=\tilde{\alpha}_2 < d_1=\tilde{d}_1$ and $|\sigma^{9}(\Delta)|=n-(\alpha_1-\alpha_2)+\frac{\alpha_1-\alpha_2}{2}+\frac{\alpha_1-\alpha_2}{2}-1+2=n+1$. Thus we deduce that $\sigma^{9}(\Delta)\in P_{k,9}(m,n+1)$. To show that $\sigma^{9}$ is an injection from $P^{9}_{k}(m,n)$ to $P_{k,9}(m,n+1)$, let 
\[
H_{k,9}(m,n+1)=\{\sigma^{9}(\Delta)\colon \Delta\in P^{9}_{k}(m,n)\}.
\]
By the above analysis, we know that $H_{k,9}(m,n+1)\subseteq P_{k,9}(m,n+1)$. Now we construct the map $\zeta^9$ from $H_{k,9}(m,n+1)$ to $P^{9}_{k}(m,n)$. Let 
\[
\tilde{\Delta}=(\tilde{\alpha}, \tilde{\beta}, \tilde{\gamma}^{1}, \ldots, \tilde{\gamma}^{k-2}, \tilde{\varpi}^{1}, \ldots, \tilde{\varpi}^{k-1})\in H_{k,9}(m,n+1)
\]
be the $(2k-1)$-tuple partition in $H_{k,9}(m,n+1)$. We know that $\tilde{d}_1=\tilde{d}_{k-1}\ge 2$, $\tilde{d}_1-\tilde{\alpha}_{1}$ is even. Moreover, by the construction of $\sigma^7$ in \eqref{eq-sigma9}, we see that $1\le \tilde{\alpha}_1=\tilde{\alpha}_2<\tilde{d}_1$, $f_{1}(\tilde{\alpha})\ge \frac{\tilde{d}_1-\tilde{\alpha}_1}{2}$, $f_{1}(\tilde{\beta})\ge \frac{\tilde{d}_1-\tilde{\alpha}_1}{2}-1$ and there exists $i$ such that $\beta_i=2$. Define
\begin{align*}
    \zeta^{9}(\tilde{\Delta}):=&(\alpha, \beta, \gamma^{1}, \ldots, \gamma^{k-2}, \varpi^{1}, \ldots, \varpi^{k-1})\\
    =&((\tilde{d}_1, \tilde{\alpha}_2, \ldots, \tilde{\alpha}_{\ell(\tilde{\alpha})-\frac{\tilde{d}_1-\tilde{\alpha}_{1}}{2}}), \\
    &(\tilde{\beta}_{1}, \ldots, \tilde{\beta}_{i-1},\tilde{\beta}_{i+1},\ldots,\tilde{\beta}_{\ell(\tilde{\beta})-(\frac{\tilde{d}_1-\tilde{\alpha}_{1}}{2}-1)}), \underbrace{\emptyset, \ldots, \emptyset}_{k-2}, \tilde{\varpi}^{1}, \ldots, \tilde{\varpi}^{k-1}).
\end{align*}
It is easy to check that $d_{1}=d_{k-1}\ge 2$. Moreover, we have
\[
\alpha_1=\tilde{d}_{1}>\tilde{\alpha}_2=\alpha_2\ge 1,
\]
using the fact $1\le \tilde{\alpha}_1=\tilde{\alpha}_2<\tilde{d}_1$. Furthermore, $\alpha_1-\alpha_2=\tilde{d}_1-\tilde{\alpha}_1$ is even and
\[
\ell(\alpha)-\ell({\beta})=\ell(\tilde{\alpha})-\frac{\tilde{d}_1-\tilde{\alpha}_1}{2}-(\ell(\tilde{\beta})-(\frac{\tilde{d}_1-\tilde{\alpha}_1}{2}-1+1))=m.
\]
and it can be checked that 
\[
|\zeta^{9}(\tilde{\Delta})|=n+1+(\tilde{d}_1-\tilde{\alpha}_{1})-\frac{\tilde{d}_1-\tilde{\alpha}_{1}}{2}-(\frac{\tilde{d}_1-\tilde{\alpha}_{1}}{2}-1)-2=n.
\]
Thus we verified that $\zeta^{9}(\tilde{\Delta})\in P^{9}_{k}(m,n)$. It is clear that for any $\Delta\in P_k^9(m,n)$
\[
\zeta^9(\sigma^9(\Delta))=\Delta.
\]
This yields that $\sigma^{9}$ is an injection from $P^{9}_{k}(m,n)$ to $P_{k,9}(m,n+1)$.
\end{proof}

For example, let $\Delta=((3,1), (1), \emptyset, (3^3), (3^3))\in P^{9}_{3}(1,23)$, using the $\sigma^{9}$ on $\Delta$, we get $\tilde{\Delta}=((1,1,1), (2,1), \emptyset, (3^3), (3^3))\in P_{3,9}(1,24)$. Applying $\zeta^9$ on  $\tilde{\Delta}$, we recover $\Delta$.

\begin{lem}\label{lem-pp10}
For $k\ge 3$, $n\ge k-1$, $m\ge 0$,  there exists a bijection $\sigma^{10}$ between $P^{10}_{k}(m,n)$ and $P_{k,10}(m,n+1)$.
\end{lem}
\begin{proof}
Given $\Delta=(\alpha, \beta, \gamma^{1}, \ldots, \gamma^{k-2}, \varpi^{1}, \ldots, \varpi^{k-1})\in P^{10}_{k}(m,n)$, by definition, we know $d_{1}=d_{k-1}\ge 2$,  $\alpha_{1}=\alpha_{2}=d_{1}>\beta_1\ge 1$. Define
\begin{align*}
    \sigma^{10}(\Delta):=&(\tilde{\alpha}, \tilde{\beta}, \tilde{\gamma}^{1}, \ldots, \tilde{\gamma}^{k-2}, \tilde{\varpi}^{1}, \ldots, \tilde{\varpi}^{k-1})\\
    =&(\alpha, (\beta_{1}+1, \beta_{2}, \ldots), \underbrace{\emptyset, \ldots, \emptyset}_{k-2}, \varpi^{1}, \ldots, \varpi^{k-1}).
\end{align*}
It is obvious that $\sigma^{10}(\Delta)\in P_{k,10}(m,n+1)$ and $\sigma^{10}$ is a bijection. We omit the trivial verification steps.
\end{proof}

For example, let $\Delta=((3,3), (1,1), \emptyset, (3^3), (3^3))\in P^{10}_{3}(0,26)$, using the $\sigma^{10}$ on $\Delta$, we get $\tilde{\Delta}=((3,3), (2,1), \emptyset, (3^3), (3^3))\in P_{3,10}(0,27)$. Applying $\zeta^{10}$ on  $\tilde{\Delta}$, we recover $\Delta$.

\begin{lem}\label{lem-pp11}
For $k\ge 3$, $n\ge k-1$, $m\ge 0$,  there exists a bijection $\sigma^{11}$ between $P^{11}_{k}(m,n)$ and $P_{k,11}(m,n+1)$.
\end{lem}
\begin{proof}
$\Delta=(\alpha, \beta, \gamma^{1}, \ldots, \gamma^{k-2}, \varpi^{1}, \ldots, \varpi^{k-1})\in P^{11}_{k}(m,n)$, by definition, we know $d_{1}=d_{k-1}\ge 2$,  $\alpha_{1}=\alpha_{2}=\beta_{1}=d_{1}$. Define
\begin{align*}
    \sigma^{11}(\Delta):=&(\tilde{\alpha}, \tilde{\beta}, \tilde{\gamma}^{1}, \ldots, \tilde{\gamma}^{k-2}, \tilde{\varpi}^{1}, \ldots, \tilde{\varpi}^{k-1})\\
    =&((\alpha_{2}, \ldots), (\beta_{2}, \ldots), \underbrace{\emptyset, \ldots, \emptyset}_{k-2}, ((d_1+1)^{d_1+1}), \varpi^{2}, \ldots, \varpi^{k-1}).
\end{align*}
Note that for $1\le i\le k-2$, $\tilde{d}_1=d_1+1=d_{i+1}+1=\tilde{d}_{i+1}+1\ge 3$ and $\tilde{\gamma}^{i}=\emptyset$. $\tilde{\alpha}_1=\alpha_2=d_1=\tilde{d}_1-1$. Moreover, $|\sigma^{11}(\Delta)|=n-|\alpha_1|-|\beta_1|+2d_1+1=n-2d_1+2d_1+1=n+1$. Thus we deduce that $\sigma^{11}(\Delta)\in P_{k,11}(m,n+1)$. To show that $\sigma^{11}$ is a bijection, now we construct the inverse map of $\sigma^{11}$. Given 
\[
\tilde{\Delta}=(\tilde{\alpha}, \tilde{\beta}, \tilde{\gamma}^{1}, \ldots, \tilde{\gamma}^{k-2}, \tilde{\varpi}^{1}, \ldots, \tilde{\varpi}^{k-1})\in P_{k,11}(m,n+1).
\]
By definition, we know that for $1\le i\le k-2$, $\tilde{d}_1=\tilde{d}_{i+1}+1\ge 3$, $\tilde{\gamma}^{i}=\emptyset$ and $\tilde{\alpha}_1=\tilde{d}_1-1$. Let
\begin{align*}
    \zeta^{11}(\Delta):=&(\alpha, \beta, \gamma^{1}, \ldots, \gamma^{k-2}, \varpi^{1}, \ldots, \varpi^{k-1})\\
    =&(\tilde{d}_1-1, \tilde{\alpha}_{1}, \ldots), (\tilde{d}_1-1, \tilde{\beta}_{1}, \ldots), \underbrace{\emptyset, \ldots, \emptyset}_{k-2}, (\tilde{d}_1-1)^{\tilde{d}_1-1}, \tilde{\varpi}^{2}, \ldots, \tilde{\varpi}^{k-1}).
\end{align*}
From the construction of $\zeta^{11}$, we see that for $1\le i\le k-2$, 
\[
d_1=\tilde{d}_1-1=\tilde{d}_{i+1}=d_{i+1},
\]
and $\gamma^{i}=\emptyset$. Moreover, $\alpha_{1}=\beta_1=\tilde{d}_1-1=\tilde{\alpha}_{1}=\alpha_2$. Furthermore, $|\zeta^{11}(\Delta)|=n+1-(2\tilde{d}_1-1)+2(\tilde{d}_1-1)=n$. Hence we verified that $\zeta^{11}(\Delta)\in P^{11}_{k}(m,n)$. It is easy to check that $\zeta^{11}$ is the inverse map of $\sigma^{11}$. This completes the proof.
\end{proof}

For example, let $\Delta=((3,3,1), (3,2), \emptyset, (3^3), (3^3))\in P^{11}_{3}(1,30)$, using the $\sigma^{11}$ on $\Delta$, we get $\tilde{\Delta}=((3,1), (2), \emptyset, (4^4), (3^3))\in P_{3,11}(1,31)$. Applying $\zeta^{11}$ on  $\tilde{\Delta}$, we recover $\Delta$.

\begin{lem}\label{lem-pp12}
For $k\ge 3$, $n\ge k-1$, $m\ge 0$,  there exists a bijection $\sigma^{12}$ between $P^{12}_{k}(m,n)$ and $P_{k,12}(m,n+1)$.
\end{lem}
\begin{proof}
Let $\Delta=(\alpha, \beta, \gamma^{1}, \ldots, \gamma^{k-2}, \varpi^{1}, \ldots, \varpi^{k-1})\in P^{12}_{k}(m,n)$, by definition, we know $d_{1}=d_{k-1}\ge 2$,  $\alpha_{1}=\alpha_2=\alpha_{3}=d_1$ and $\beta=\emptyset$. Let $t$ be the maximum integer such that $\alpha_t=d_1$. Clearly, $t\ge 3$. Define
\begin{align}
    \sigma^{12}(\Delta):=&(\tilde{\alpha}, \tilde{\beta}, \tilde{\gamma}^{1}, \ldots, \tilde{\gamma}^{k-2}, \tilde{\varpi}^{1}, \ldots, \tilde{\varpi}^{k-1})\nonumber\\
    =&((\alpha_{1}, \ldots, \alpha_{t-1}, \alpha_{t}-1, \alpha_{t+1}, \ldots,\alpha_{\ell(\alpha)}, 1), (1), \underbrace{\emptyset, \ldots, \emptyset}_{k-2}, \varpi^{1}, \ldots, \varpi^{k-1})\label{eq-jt1}.
\end{align}
Note that $\tilde{d}_1=d_1=d_{k-1}=\tilde{d}_{k-1}\ge 2$, $\tilde{\gamma}^{i}=\emptyset$ ($1\le i\le k-2$), $\tilde{\alpha}_1=\alpha_1=d_1=\tilde{\alpha}_2=\tilde{d}_1$, $\tilde{\alpha}_{t}=\alpha_{t}-1=d_1-1=\tilde{d}_1-1$ and $\beta=(1)$.  Moreover, $|\sigma^{12}(\Delta)|=n+1$. Thus we deduce that $\sigma^{12}(\Delta)\in P_{k,12}(m,n+1)$. To show that $\sigma^{12}$ is a bijection, we construct the inverse map of $\sigma^{12}$. Given 
\[
\tilde{\Delta}=(\tilde{\alpha}, \tilde{\beta}, \tilde{\gamma}^{1}, \ldots, \tilde{\gamma}^{k-2}, \tilde{\varpi}^{1}, \ldots, \tilde{\varpi}^{k-1})\in P_{k,12}(m,n+1).
\]
By definition, we know that $\tilde{d}_1=\tilde{d}_{k-1}\ge 2$, $\tilde{\gamma}^{i}=\emptyset$ ($1\le i\le k-2$), $f_{\tilde{d}_1-1}(\tilde{\alpha})\ge 1$ and $\tilde{\beta}=(1)$. Moreover, let $j\ge 2$ denote the maximum integer such that $\tilde{\alpha}_j=\tilde{d}_1$, which means that $\tilde{d}_1=\tilde{\alpha}_1=\cdots=\tilde{\alpha}_{j}>\tilde{\alpha}_{j+1}=\tilde{d}_1-1$.
Let
\begin{align}
    \zeta^{12}(\Delta):=&(\alpha, \beta, \gamma^{1}, \ldots, \gamma^{k-2}, \varpi^{1}, \ldots, \varpi^{k-1})\nonumber\\
    =&((\tilde{\alpha}_1, \ldots, \tilde{\alpha}_j, \tilde{\alpha}_{j+1}+1, \tilde{\alpha}_{j+2}, \ldots, \tilde{\alpha}_{\ell(\tilde{\alpha})-1}), \emptyset, \underbrace{\emptyset, \ldots, \emptyset}_{k-2}, \tilde{\varpi}^{1}, \tilde{\varpi}^{2}, \ldots, \tilde{\varpi}^{k-1})\label{eq-jt2}.
\end{align}
From the construction of $\zeta^{12}$, we see that $d_1=\tilde{d}_1=\tilde{d}_{k-1}=d_{k-1}$, $\beta=\emptyset$ and for $1\le i\le k-2$, $\gamma^{i}=\emptyset$. Moreover, when $j=2$, we have $\alpha_{3}=\tilde{\alpha}_{3}+1=\tilde{d}_1-1+1=d_1$. For $j\ge 3$, we have $\alpha_{3}=\tilde{\alpha}_{3}=d_1$. So in either case, we conclude that $\alpha_3=d_1$. Furthermore, $|\zeta^{12}(\Delta)|=n+1-2+1=n$. Hence, we have verified that $\zeta^{12}(\Delta)\in P^{12}_{k}(m,n)$. It is straightforward to check that the number $j$ in \eqref{eq-jt2} coincides with the number $t-1$ in \eqref{eq-jt1}. Thus, $\zeta^{12}$ is the inverse map of $\sigma^{12}$. This completes the proof.
\end{proof}

For example, let $\Delta=((2,2,2,1,1),\emptyset, \emptyset, (2^2), (2^2))\in P^{12}_{3}(5,16)$, using the $\sigma^{12}$ on $\Delta$, we get $t=3$, $\tilde{\Delta}=((2,2,1,1,1,1), (1), \emptyset, (2^2), (2^2))\in P_{3,12}(5,17)$. Applying $\zeta^{12}$ on  $\tilde{\Delta}$, we deduce $j=2$ and $\zeta^{12}(\sigma^{12}(\Delta))=\Delta$.

\begin{lem}\label{lem-pp13}
For $k\ge 4$, $n\ge k-1$, $m\ge 0$,  there exists a bijection $\sigma^{13}$ between $P^{13}_{k}(m,n)$ and $P_{k,13}(m,n+1)$.
\end{lem}
\begin{proof}
Let $\Delta=(\alpha, \beta, \gamma^{1}, \ldots, \gamma^{k-2}, \varpi^{1}, \ldots, \varpi^{k-1})\in P^{13}_{k}(m,n)$, by definition, we know $k\ge 4$, $d_{1}=d_{k-1}\ge 2$,  $\alpha_{1}=\alpha_{2}=d_1>\alpha_3$ and $\beta=\emptyset$. Define
\begin{align*}
    \sigma^{13}(\Delta):=&(\tilde{\alpha}, \tilde{\beta}, \tilde{\gamma}^{1}, \ldots, \tilde{\gamma}^{k-2}, \tilde{\varpi}^{1}, \ldots, \tilde{\varpi}^{k-1})\\
    =&((\alpha_{1}-1, \alpha_2-1, \alpha_3, \ldots), \emptyset, \underbrace{\emptyset, \ldots, \emptyset}_{k-3},(1), (d_1+1)^{d_1+1}, \varpi^{2} \ldots, \varpi^{k-2}, (d_{k-1}-1)^{d_{k-1}-1}).
\end{align*}
Note that for $2\le i\le k-2$, $\tilde{d}_1=d_1+1=d_i+1=\tilde{d}_i+1=d_{k-1}-1+2=\tilde{d}_{k-1}+2\ge 3$, $\tilde{\gamma}^{i-1}=\emptyset$, $\tilde{\gamma}^{k-2}=(1)$ and $\tilde{\beta}=\emptyset$. Moreover, 
\[
\tilde{\alpha}_1=\alpha_1-1=d_1-1=\tilde{\alpha}_2=\tilde{d}_1-2=\tilde{d}_{k-1}\ge \tilde{\alpha}_3=\alpha_3,
\]
Furthermore, $|\sigma^{13}(\Delta)|=n-(2d_1-1)-2+(2(d_1+1)-1)+1=n+1$. Thus we deduce that $\sigma^{13}(\Delta)\in P_{k,13}(m,n+1)$. To show that $\sigma^{13}$ is a bijection, now we construct the inverse map of $\sigma^{13}$. Given 
\[
\tilde{\Delta}=(\tilde{\alpha}, \tilde{\beta}, \tilde{\gamma}^{1}, \ldots, \tilde{\gamma}^{k-2}, \tilde{\varpi}^{1}, \ldots, \tilde{\varpi}^{k-1})\in P_{k,13}(m,n+1).
\]
By definition, we know that for $2\le i\le k-2$, $\tilde{d}_1=\tilde{d}_i+1=\tilde{d}_{k-1}+2\ge 3$, $\tilde{\gamma}^{i-1}=\emptyset$, $\tilde{\gamma}^{k-2}=(1)$, $\tilde{\alpha}_1=\tilde{\alpha}_2=\tilde{d}_{k-1}$ and $\tilde{\beta}=\emptyset$. 
Let
\begin{align*}
    \zeta^{13}(\Delta):=&(\alpha, \beta, \gamma^{1}, \ldots, \gamma^{k-2}, \varpi^{1}, \ldots, \varpi^{k-1})\\
    =&((\tilde{\alpha}_1+1, \tilde{\alpha}_2+1, \tilde{\alpha}_3, \ldots), \emptyset, \underbrace{\emptyset, \ldots, \emptyset}_{k-2}, (\tilde{d}_1-1)^{\tilde{d}_1-1}, \tilde{\varpi}^{2}, \ldots, \tilde{\varpi}^{k-2}, (\tilde{d}_{k-1}+1)^{\tilde{d}_{k-1}+1}).
\end{align*}
From the construction of $\zeta^{13}$, we see for $2\le i\le k-2$, $d_1=\tilde{d}_1-1=\tilde{d}_{i}=d_i=\tilde{d}_{k-1}+1=d_{k-1}$, $\beta=\emptyset$ and for $1\le i\le k-2$, $\gamma^{i}=\emptyset$. Moreover, $\alpha_{1}=\tilde{\alpha}_1+1=\alpha_2=d_1=d_{k-1}$. Furthermore, $|\zeta^{13}(\Delta)|=n+1-(2\tilde{d}_1-1)+(2(\tilde{d}_1-1)-1)-1+2=n$. Hence we verified that $\zeta^{13}(\Delta)\in P^{13}_{k}(m,n)$. It is easy to check that $\zeta^{13}$ is the inverse map of $\sigma^{13}$. This completes the proof.
\end{proof}

For example, let $\Delta=((2,2,1,1),\emptyset, \emptyset, \emptyset, (2^2), (2^2), (2^2))\in P^{13}_{4}(4,18)$, using the $\sigma^{13}$ on $\Delta$, we get $\tilde{\Delta}=((1,1,1,1), \emptyset, \emptyset, (1), (3^3), (2^2), (1^1))\in P_{4,13}(4,19)$. Applying $\zeta^{13}$ on  $\tilde{\Delta}$, we recover $\Delta$.

\begin{lem}\label{lem-pp14}
For $k=3$, $n\ge 2$, $m\ge 0$,  there exists a bijection $\sigma^{14}$ between $P^{14}_{3}(m,n)$ and $P_{3,14}(m,n+1)$.
\end{lem}
\begin{proof}
Let $\Delta=(\alpha, \beta, \gamma^{1},  \varpi^{1},  \varpi^{2})\in P^{14}_{3}(m,n)$, by definition, we know $d_{1}=d_{2}\ge 3$,  $\alpha_{1}=\alpha_{2}=d_1>\alpha_3$ and $\beta=\emptyset$. Moreover, let $t$ be the maximum integer such that $\alpha_t\ge d_1-1$, clearly $t\ge 2$ and $\alpha_t> d_1-2\ge \alpha_{t+1}$. Define
\begin{align}
    \sigma^{14}(\Delta):=&(\tilde{\alpha}, \tilde{\beta}, \tilde{\gamma}^{1}, \tilde{\varpi}^{1}, \tilde{\varpi}^{2})\nonumber\\
    =&((\alpha_{1}-1,  \alpha_3, \ldots, \alpha_t, \alpha_2-2, \alpha_{t+1}, \ldots), \emptyset, (1, 1), (d_1+1)^{d_1+1}, (d_{2}-1)^{d_{2}-1})\label{eq-jjtt1}.
\end{align}
Note that $\tilde{d}_1=d_1+1=d_2-1+2=\tilde{d}_2+2\ge 4$, $\tilde{\gamma}^{1}=(1,1)$ and $\tilde{\alpha}_1=\alpha_1-1=d_1-1=\tilde{d}_2$. Moreover, $\tilde{\alpha}_{t}=\alpha_2-2=\alpha_1-2=\tilde{\alpha}_1-1$.  
Furthermore, $|\sigma^{14}(\Delta)|=n-(2d_2-1)-3+(2(d_1+1)-1)+2=n+1$. Thus we deduce that $\sigma^{14}(\Delta)\in P_{3,14}(m,n+1)$. To show that $\sigma^{14}$ is a bijection, now we construct the inverse map of $\sigma^{14}$. Given 
\[
\tilde{\Delta}=(\tilde{\alpha}, \tilde{\beta}, \tilde{\gamma}^{1}, \tilde{\varpi}^{1}, \tilde{\varpi}^{2})\in P_{3,14}(m,n+1).
\]
By definition, we know that $\tilde{d}_1=\tilde{d}_2+2\ge 4$, $\tilde{\gamma}_1=(1,1)$, $\tilde{\alpha}_1=\tilde{d}_2$, and there exists $j\ge 2$ such that $\tilde{\alpha}_j=\tilde{\alpha}_1-1$, we choose such $j$ to be minimum.
Let
\begin{align}
    \zeta^{14}(\Delta):=&(\alpha, \beta, \gamma^{1}, \varpi^{1}, \varpi^{2})\nonumber\\
    =&((\tilde{\alpha}_1+1, \tilde{\alpha}_{j}+2, \tilde{\alpha}_2, \ldots, \tilde{\alpha}_{j-1},  \tilde{\alpha}_{j+1}, \ldots), \emptyset, \emptyset, (\tilde{d}_1-1)^{\tilde{d}_1-1}, (\tilde{d}_{2}+1)^{\tilde{d}_{2}+1})\label{eq-jjtt2}.
\end{align}
From the construction of $\zeta^{14}$, we see $d_1=\tilde{d}_1-1=\tilde{d}_{2}+1=d_2\ge 3$, $\beta=\gamma^{1}=\emptyset$ and $\alpha_1=\tilde{\alpha}_1+1=\tilde{d}_2+1=d_2=d_1=\tilde{\alpha}_{j}+2=\alpha_2>\alpha_3=\tilde{\alpha}_2$ which follows $\tilde{\alpha}_2\le \tilde{\alpha}_1$. Furthermore, $|\zeta^{14}(\Delta)|=n+1-(2\tilde{d}_1-1)-2+(2(\tilde{d}_1-1)-1)+1+2=n$. Hence we verified that $\zeta^{14}(\Delta)\in P^{14}_{3}(m,n)$. It is easy to check that the number $j$ in \eqref{eq-jjtt2} coincides with the number $t$ in \eqref{eq-jjtt1}. Thus $\zeta^{14}$ is the inverse map of $\sigma^{14}$. This completes the proof.
\end{proof}

For example, let $\Delta=((3,3,2,2,1), \emptyset, \emptyset, (3^3), (3^3))\in P^{14}_{3}(5,29)$,  using the $\sigma^{14}$ on $\Delta$, we get $t=4$ and $\tilde{\Delta}=((2,2,2,1,1), \emptyset, (1,1), (4^4), (2^2))\in P_{3,14}(5,30)$. Applying $\zeta^{14}$ on  $\tilde{\Delta}$, we deduce $j=4$ and $\zeta^{14}(\sigma^{14}(\Delta))=\Delta$.

\begin{lem}\label{lem-pp15}
For $k= 3$, $n\ge 2$, $m\ge 0$,  there exists a bijection $\sigma^{15}$ between $P^{15}_{3}(m,n)$ and $P_{3,15}(m,n+1)$.
\end{lem}
\begin{proof}
Let $\Delta=(\alpha, \beta, \gamma^{1}, \varpi^{1}, \varpi^{2})\in P^{15}_{3}(m,n)$, by definition, we $d_{1}=d_{2}=2$,  $\alpha_{1}=\alpha_{2}=d_1=2>\alpha_3$ and $\beta=\emptyset$. Define
\begin{align*}
    \sigma^{15}(\Delta):=&(\tilde{\alpha}, \tilde{\beta}, \tilde{\gamma}^{1}, \tilde{\varpi}^{1}, \tilde{\varpi}^{2})\\
    =&((\alpha_{1}-1,  \alpha_2-1, 1, \alpha_3, \ldots, \alpha_{\ell(\alpha)}), (2), \emptyset, \varpi^{1}, \varpi^{2}).
\end{align*}
Note that $\tilde{d}_1=d_1=d_2=\tilde{d}_2=2$, $\tilde{\alpha}_1=\alpha_1-1=1=\tilde{\alpha}_2=\tilde{\alpha}_3$, $\ell(\tilde{\alpha})=\ell(\alpha)+1\ge 3$ and $\beta=(2)$. Moreover, $|\sigma^{15}(\Delta)|=n-2+2+1=n+1$. Thus we deduce that $\sigma^{15}(\Delta)\in P_{3,15}(m,n+1)$. To show that $\sigma^{15}$ is a bijection, now we construct the inverse map of $\sigma^{15}$. Given 
\[
\tilde{\Delta}=(\tilde{\alpha}, \tilde{\beta}, \tilde{\gamma}^{1}, \tilde{\varpi}^{1}, \tilde{\varpi}^{2})\in P_{3,15}(m,n+1).
\]
By definition, we know that $\tilde{d}_1=\tilde{d}_2=2$, $\tilde{\alpha}_1=1$, $\ell(\tilde{\alpha})\ge 3$ and $\tilde{\beta}=(2)$.
Let
\begin{align*}
    \zeta^{15}(\Delta):=&(\alpha, \beta, \gamma^{1}, \varpi^{1}, \varpi^{2})\\
    =&((\tilde{\alpha}_1+1, \tilde{\alpha}_{2}+1, \ldots, \tilde{\alpha}_{\ell(\tilde{\alpha})-1}), \emptyset, \emptyset, \tilde{\varpi}^1, \tilde{\varpi}^2).
\end{align*}
From the construction of $\zeta^{15}$, we see $d_1=\tilde{d}_1=\tilde{d}_2=d_2=2$, $\beta=\gamma^{1}=\emptyset$ and $\alpha_1=\tilde{\alpha}_1+1=\tilde{d}_2=d_2=2$. Moreover, $|\zeta^{15}(\Delta)|=n+1-3+2=n$. Hence we verified that $\zeta^{15}(\Delta)\in P^{15}_{3}(m,n)$. It is easy to check that $\zeta^{15}$ is the inverse map of $\sigma^{15}$. This completes the proof.
\end{proof}

For example, let $\Delta=((2,2), \emptyset, \emptyset, (2^2), (2^2))\in P^{15}_{3}(2,12)$, using the $\sigma^{15}$ on $\Delta$, we get $\tilde{\Delta}=((1,1,1), (2), \emptyset, (2^2), (2^2))\in P_{3,15}(2,13)$. Applying $\zeta^{15}$ on  $\tilde{\Delta}$, we recover $\Delta$.

{\noindent \it{Proof of Theorem \ref{thm-PnPn+1}. }} For $k\ge 3$, $m\ge 0$, $n\ge k-1$, let $\Delta$ be a partition in $P_{k}(m,n)$. If $\Delta\in P^{i}_{k}(m,n)$, define
\begin{equation*}
    \sigma(\Delta)=\sigma^{i}(\Delta), 
\end{equation*}
where $1\le i\le 15$.
If $\Delta\in P^{16}_{k}(m,n)$, there are only two cases.

Case 1: $d_1=d_{k-1}=1$ and $n=m+k-1$; then $(m,k,n)=(m,k,m+k-1)$.

Case 2: $k=3$, $d_1=d_2=2$, $\alpha=\emptyset$; then $\beta=\emptyset$, thus $\Delta=(\emptyset, \emptyset, \emptyset, (2^2), (2^2))$ and $(m,k,n)=(0,3,8)$.

By Lemmas \ref{lem-pp1} $\sim$ \ref{lem-pp15}, we conclude that $\sigma$ is an injection from $P_{k}(m,n)$ to $P_{k}(m,n+1)$ except for the two cases in $P^{16}_{k}(m,n)$.
\qed

\noindent{\bf Acknowledgments.}   This work was supported by the National Science Foundation of China grants 12171358 and 12371336.

\end{document}